%% file: HSWEMED.tex
\documentclass[10pt,a4paper]{article}
\usepackage[margin=0.7in]{geometry} 
\usepackage{setspace}
\usepackage{authblk}
\usepackage{amsmath}
\usepackage{amsthm}
\usepackage{amssymb}
\usepackage{nccmath}
\usepackage{bm}
\usepackage{svg}
\usepackage{multicol}
\usepackage{multirow}
\usepackage[makeroom]{cancel}
\usepackage{wrapfig}
\usepackage{booktabs}
\usepackage[titletoc]{appendix}
\usepackage{xcolor}
\usepackage{tikz}
\usetikzlibrary{arrows.meta,calc,decorations.pathreplacing,positioning,fit,backgrounds}
\definecolor{reconblue}{RGB}{0,84,180}
\definecolor{avgred}{RGB}{190,30,45}
\definecolor{gridgray}{RGB}{170,170,170}
\usepackage[normalem]{ulem}

\usepackage{titlesec}
\newcommand{\secref}[1]{\textcolor{blue}{\hyperref[#1]{Section~\ref*{#1}}}}

\makeatletter
\renewcommand{\theparagraph}{\thesubsubsection.\arabic{paragraph}}

\renewcommand{\paragraph}[1]{%
	\refstepcounter{paragraph}
	\par\noindent\textbf{\theparagraph~#1}\par
}

\newcommand{\parref}[1]{\hyperref[#1]{Section~\ref*{#1}}}
\makeatother

\titlespacing*{\paragraph}
{0pt}{1.25ex plus 1ex minus .2ex}{1.5ex plus .2ex}

\usepackage[utf8]{inputenc}
\usepackage{hyperref}
\hypersetup{
	unicode,
	colorlinks=true,
	linktoc=all,
	linkcolor=blue,
	breaklinks,
	urlcolor=cyan   
}

\newcommand{\eqrefc}[1]{\textcolor{red}{\hyperref[#1]{\eqref{#1}}}}

\usepackage[backend=biber,sorting=nyt,giveninits=true,style=numeric,maxnames=99]{biblatex}
\DeclareNameAlias{author}{given-family}

\newcommand{\theoremref}[1]{\textcolor{green}{\hyperref[#1]{Theorem~\ref*{#1}}}}

\newcommand{\appref}[1]{\textcolor{purple}{\hyperref[#1]{Appendix~\ref*{#1}}}}

\newcommand{\remref}[1]{\textcolor{orange}{\hyperref[#1]{Remark~\ref*{#1}}}}

\newcommand{\tabref}[1]{\textcolor{orange}{\hyperref[#1]{Table~\ref*{#1}}}}

\theoremstyle{plain}
\newtheorem{theorem}{Theorem}[section]

\newtheorem{proposition}[theorem]{Proposition}

\theoremstyle{definition}
\newtheorem{definition}{Definition}[section]

\newtheorem{remark}{Remark}[section]

\usepackage{graphicx, color}
\usepackage[skip=5pt]{caption}
\graphicspath{{fig/outputs/}}
\usepackage{svg}
\svgpath{{fig/outputs/}}
\usepackage{subcaption}

\usepackage[mathlines]{lineno} 

\modulolinenumbers[1]         

\newcommand{\figref}[1]{\textcolor{red}{\hyperref[#1]{Figure~\ref*{#1}}}}

\usepackage{algorithm, algpseudocode} 
\usepackage{mathrsfs} 
\usepackage{array}

\title{A Moment-Based Shallow Water Model for Coupled Bedload-Suspended Load Transport with Variable Density: Hyperbolicity and Numerical Validation}
\author[1,2]{Afroja Parvin \thanks{Corresponding author. \textit{Email address:}  \texttt{afroja.parvin@kuleuven.be}}}
\author[1]{Giovanni Samaey}
\author[3,4]{Julian Koellermeier}

\affil[1]{\small Department of Computer Science, KU Leuven, Belgium}
\affil[2]{\small School of Mathematical Sciences, Peking University, China}
\affil[3]{\small Bernoulli Institute, University of Groningen, Netherlands}
\affil[4]{\small Department of Mathematics, Computer Science and Statistics, Ghent University, Belgium}
\makeatletter
\renewcommand{\@date}{} 
\makeatother

\begin{document}
	  \setstretch{1.1}
	   \maketitle	
\begin{abstract}
		In this paper, we develop a Hyperbolic Shallow Water Exner Moment model with Erosion and Deposition (HSWEMED), which extends the recent shallow water moment framework to capture coupled morphodynamics with erosion and deposition. Extending existing moment models, HSWEMED introduces a sediment concentration equation for suspended load, couples a concentration-dependent sediment-water mixture density with the momentum equation and higher-order moments, and adds source terms arising from erosion and deposition. Starting from the incompressible Navier-Stokes equations for a water-sediment mixture, we derive a coupled system that integrates: (i) the shallow water equations, (ii) moment equations for polynomial velocity coefficients, (iii) a depth-averaged suspended-sediment equation, and (iv) an Exner equation for bedload transport with erosion-deposition coupling. The reduced model evolves only the depth-averaged concentration, while the vertical structure is retained for the horizontal velocity through the moment expansion. 
		For the regularized model, we derive the characteristic polynomial and identify a real-root condition for the sediment-coupled part. The analysis shows that strict hyperbolicity can hold for even moment orders, while odd moment orders lead to weak hyperbolicity because of a repeated eigenvalue.
		We also derive dissipative energy balance relations for the lower-order models. Numerical results are obtained with a second-order path-conservative MUSCL-Rusanov finite-volume scheme with a well-balanced correction for the water-at-rest equilibrium. We verify the well-balanced property, perform a reference-convergence study, and compare dam-break simulations with laboratory experiments from the literature. The results illustrate how the moment approximation and the erosion-deposition coupling affect the free surface, bed evolution, and suspended concentration.
		
    \vskip 0.5cm 
	\noindent\textbf{Keywords:} Shallow water Exner moment model; erosion-deposition; bedload and suspended load sediment transport; hyperbolicity; dam-break, wet/dry fronts.
	\end{abstract}
        \input{01-Introduction}

        \input{02-Mathematical_model}
        \input{03-Hyperbolicity_modified}
        \input{04-Numerical_simulation_modified}
        \input{05-Conclusion}

	\section*{Acknowledgments} 
		Afroja Parvin acknowledges the financial support from KU Leuven Global PhD Partnership fellowship for a joint PhD with Peking University (grant agreement GPPKU/21/009). This work is part of the HiWAVE project (file no. VI.Vidi.233.066) within the NWO Vidi ENW programme, partly funded by the Dutch Research Council (NWO; grant DOI: 10.61686/CBVAB59929).
	\section*{Declarations}
	
	\textbf{Competing interests:}  On behalf of all authors, the corresponding author states that there is no conflict of interest.
	
    \printbibliography
       \section*{Appendices}
         \begin{appendices}
            \input{appendix_reference_system}
            \input{appendix_depthaveraged_equation}
            \input{appendix_energy_equation}
	        \input{Proof_of_characteristic_polynomial}\label{app:3}
       \end{appendices}
\end{document}

%% file: 01-Introduction.tex
\section{Introduction}
Studying sediment transport in shallow water is an active research topic. It involves the movement of sediment particles influenced by gravity and friction \cite{diaz2008sediment}. This process is essential for evolving riverbeds, estuaries, and coastal areas. Understanding sediment transport dynamics is essential in civil engineering, river and coastal engineering, and flood risk management \cite{gonzalez2020robust}.
\par
\begin{figure}[H]
	\centering
	\includeinkscape[inkscapelatex=true, width=0.6\textwidth]{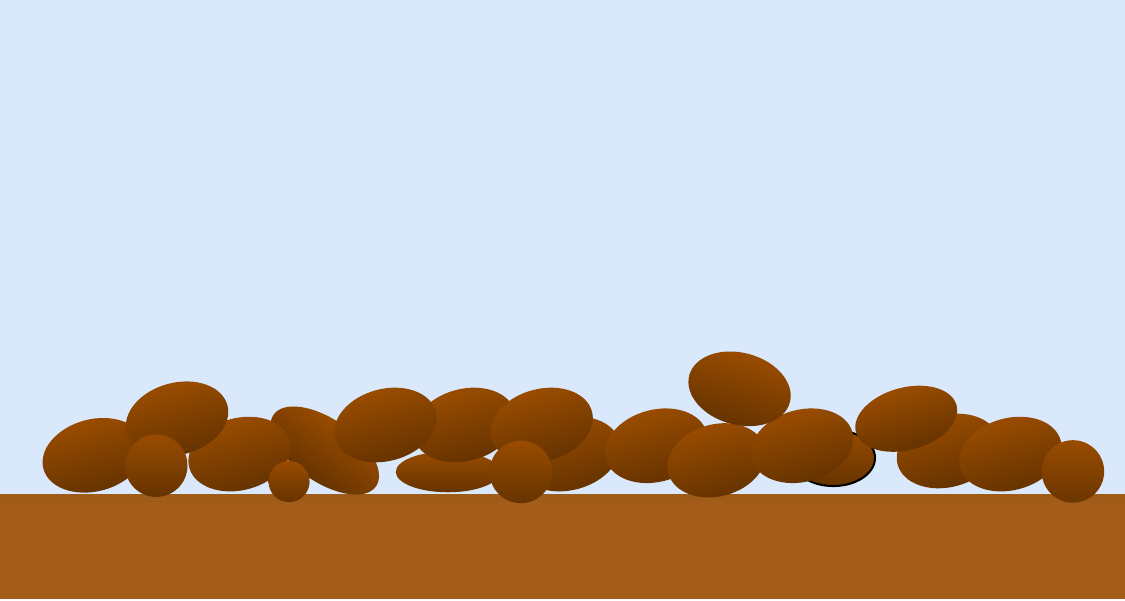}
	\caption{Sketch of sediment transport dynamics, illustrating bedload movement along the bed, suspended load within the water column, and zones of erosion and deposition controlled by variations in flow velocity and shear stress.}
	\label{model_sketch}
\end{figure}

 Sediments are usually transported as bedload along the riverbed and suspended loads in the flow. Bedload transport refers to coarser particles that move close to the bed by rolling, sliding, or jumping. In contrast, a suspended load consists of finer particles eroded from the bed, which remain suspended for a time before settling  \parencite{gonzalez2020robust,subhasish2014fluvial}(see \figref{model_sketch}).
\par The typical method for modeling bedload transport couples the Shallow Water model (SW)  \parencite{audusse2004fast,meng2020localized} with the bed evolution equation.  The result is known as the Shallow Water Exner model \cite{exner1925uber}. The Exner equation is based on mass conservation and describes the solid transport flux using a closure relation. In most existing works, this is integrated with the transport equation for suspended sediment, including erosion and deposition processes  \parencite{cantero2019vertically,del2023lagrange,gonzalez2023numerical,zhao2019depth}. There are many empirical formulas to quantify solid transport flux, erosion, and deposition at the bed surface  (see, e.g.,  \parencite{boittin2019modeling,cozzolino2014novel,el2011applicability,kubo2004experimental,morales2009shallow,nakajima2002laboratory,parker1986self,rijn1987mathematical,zhao2017comparison}, and many others).
\par  In classical shallow water sediment transport models, the near-bed velocity, essential for estimating bed shear stress, sediment flux, and erosion, is commonly approximated using a logarithmic velocity profile. Since these models do not resolve the vertical structure of the velocity, empirical friction laws such as Manning or Chézy are often introduced to parametrize shear stress, with their coefficients typically calibrated by integrating or fitting the logarithmic profile under idealized conditions, such as steady and uniform flow \parencite{afzalimehr2009determination,elgamal2021moment,lajeunesse2010bed,zordan2018entrainment}. Although tuning friction coefficients may yield reasonable results in steady-state regimes, this approach cannot capture the full dynamics of sediment transport in more general flow scenarios. Consequently, the model lacks the capability to dynamically resolve the vertical velocity structure, limiting its effectiveness in representing near-bed flow in complex or unsteady environments.
Given that velocities near the bed are critical for bedload transport, reliance on depth-averaged quantities and idealized profiles can lead to inaccurate predictions of sediment erosion under non-equilibrium conditions.
\par Recent research has aimed to improve fluid descriptions in the vertical direction using two-phase flows \cite{audusse2021asymptotic} or a multi-layer approach  \parencite{bonaventura2018multilayer,audusse2011multilayer,fernandez2013multilayer,rowan2020efficient}. However, this approach raises challenges, such as the lack of analytical expressions for eigenvalues, the requirement for many layers to capture complex velocity profiles, and increased computational cost. \cite{CiCP-25-669} developed an extended Shallow Water Moment model (SWM) to address these challenges. This model uses a Legendre polynomial expansion to represent horizontal velocity variations in the vertical direction. It applies a Galerkin projection of a transformed Navier-Stokes system for deriving evolution equations for all basis coefficients, called 'moments.' The model was recently extended for open curved shallow flows \cite{steldermann2023shallow}, non-hydrostatic pressure \cite{scholz2024dispersion}, a splitting scheme for dispersive shallow moment equations \cite{scholz2026splitting}, and multi-layer models \cite{garres2023general}. While the SWM performs well compared to the standard SW \cite{CiCP-25-669}, it lacks hyperbolicity for higher-order moments. Losing hyperbolicity can cause numerical instabilities. In \cite{koellermeier2020analysis}, the authors resolve this issue for the SWM using a hyperbolic regularization, resulting in the Hyperbolic Shallow Water Moment model (HSWM). Later, the HSWM model was extended by incorporating the Exner equation to simulate bedload transport, resulting in the Hyperbolic Shallow Water Exner Moment model (HSWEM)  \cite{garres2023general}. While this approach includes non-linear friction adapted to the moment framework, it neglects erosion, deposition, and suspension processes, and assumes a constant fluid density. These simplifications limit its ability to capture bottom evolution, since sediment-laden flows should be treated as a water-sediment mixture with variable density. To our best knowledge, there is no hyperbolic shallow water moment model that simultaneously simulates bedload and suspended load transport in a coupled framework.
 
 \par To address these limitations, this work extends the HSWEM \cite{garres2023general} to create a coupled framework for simulating bedload and suspended load sediment transport. The proposed framework incorporates vertical velocity structures and accounts for sediment mass changes at the bottom by including erosion and deposition processes. The derivation starts with the standard incompressible Navier-Stokes Equations for a water-sediment mixture with a divergence-free velocity field. Using a moment-based approach \cite{CiCP-25-669}, we derive a shallow water moment model by averaging over the vertical variable, incorporating volumetric sediment concentration for suspended load with erosion and deposition, and the Exner equation for bedload transport, resulting in the Shallow Water Exner Moment model with Erosion and Deposition (SWEMED). In our new SWEMED, incorporating sediment suspension, erosion, and deposition effects introduces extra terms in the mass conservation and Exner equations, which quantify the rate of bed deformation. Additionally, two more terms are added to the momentum conservation and higher-order averaged equations, including friction based on flow velocity at the bottom. These terms account for the influence of sediment concentration and the momentum transfer resulting from sediment exchange between the water column and the erodible bed. The higher-order averages include a sediment discharge term coupled with the Exner equation.

\par Therefore, the key novelty of this work is the extension of the shallow water moment framework to a coupled sediment-water system with erosion and deposition. Unlike previous moment models that neglect sediment exchange, our proposed model: (i) introduces a depth-averaged sediment concentration equation for suspended load, (ii) couples the variable mixture density with the momentum and higher-order moment equations, and (iii) adds new source terms arising from erosion and deposition, representing the transfer of mass and momentum between the flow and the erodible bed. The reduced model evolves only the depth-averaged concentration $c_m$. The near-bed concentration $c_b$ required in the deposition closure is prescribed algebraically from $c_m$ through an empirical relation. Accordingly, within the intended dilute to moderately dilute regime, the suspended concentration is not treated as a passive scalar: it feeds back on the hydrodynamics through the mixture density in the hydrostatic pressure and associated buoyancy terms in the depth-averaged momentum and moment equations.

\par The main contributions of this paper are the following. First, we derive the SWEMED model by coupling the shallow water moment model with a depth-averaged suspended-sediment equation and an Exner equation including erosion and deposition. Second, we analyze the hyperbolicity of the regularized system following the recent work on hyperbolic shallow water moment models \cite{koellermeier2020analysis,Garres}, resulting in the Hyperbolic Shallow Water Exner Moment model with Erosion and Deposition (HSWEMED). Third, we derive dissipative energy balance relations for the lower-order models SWEMED0 and SWEMED1. 
Fourth, we develop a finite-volume discretization for HSWEMED based on a second-order path-conservative MUSCL-Rusanov scheme \cite{vanLeer1979MUSCL}. We introduce a well-balanced correction for the water-at-rest equilibrium by modifying the Rusanov viscosity term. The corrected scheme preserves the water-at-rest equilibrium over smooth and discontinuous bottoms, and the reference-convergence study gives rates consistent with the expected second-order behaviour. Finally, we use dam-break simulations to compare HSWEMED with SWEED and HSWEM and to assess the model behaviour against laboratory experiments from the literature \cite{Spinewine01122007}. The numerical results show that, in the tested PVC case, HSWEMED predicts less bottom erosion than the depth-averaged SWEED due to the moment-based bottom velocity, while the differences are smaller for the sand case. The comparisons focus on the free surface, sediment-bottom evolution, and depth-averaged suspended concentration.

\par The remainder of this paper is organized as follows: Section 2 derives the coupled Shallow Water Moment model with Erosion and Deposition (SWEMED), incorporating the moment approach for shallow flows together with the sediment concentration equation for suspended load and the Exner equation for bedload transport. Section 3 analyzes the eigenvalue structure and hyperbolicity properties of the regularized HSWEMED transport matrix. Section 4 presents numerical results for dam-break problems and compares them with laboratory experiments from the literature. Finally, Section 5 provides concluding remarks. 
We summarize the model acronyms used throughout this work in \tabref{tab:model_acronyms}. Since several model names contain overlapping terminology, the table is included to make the comparison between the models easier to follow.
\begin{table}[h!]
	\centering
	\begin{tabular}{ll}
		\toprule
		\textbf{Acronym} & \textbf{Full Model Name} \\
		\midrule
		SW        & Shallow Water model \\
		SWEED     & Shallow Water Exner model with Erosion and Deposition \parencite{cantero2019vertically,del2023lagrange,gonzalez2023numerical,zhao2019depth} \\
		SWM       & Shallow Water Moment model \cite{CiCP-25-669} \\
		HSWM      & Hyperbolic Shallow Water Moment model \cite{koellermeier2020analysis} \\
		HSWEM     & Hyperbolic Shallow Water Exner Moment model \cite{Garres} \\
		SWEMED    & Shallow Water Exner Moment model with Erosion and Deposition (this work) \\
		HSWEMED   & Hyperbolic Shallow Water Exner Moment model with Erosion and Deposition (this work) \\
		\bottomrule
	\end{tabular}
	\caption{Summary of model acronyms and their corresponding full names used in this work.}
	\label{tab:model_acronyms}
\end{table}

%% file: 02-Mathematical_model.tex
\section{Shallow water moment model for sediment transport}
In this section we present our new model for sediment-laden shallow flows as extension of the HSWEM from \cite{Garres} by including a sediment concentration equation and erosion deposition effects. The resulting Shallow Water Exner Moment model with Erosion and Deposition (SWEMED) incorporates both hydrodynamic and morphodynamic evolution that interact through additional erosion and deposition source terms. The new system explicitly accounts for sediment concentration, suspension and erosion-deposition processes within the moment framework, leading to a set of coupled equations for the sediment-water mixture. Specifically, the hydrodynamic part is described by the depth-averaged Shallow Water Moment model (SWM) with erosion-deposition, while the morphodynamic part combines suspended load transport through the sediment concentration equation with bed evolution through the Exner equation.
\subsection{Hydrodynamic equations}\label{Hydrodynamic}
We consider the two-dimensional inhomogeneous Navier-Stokes equations
for simplicity to derive the coupled mathematical model for sediment
transport. However, the derivation can be extended to three spatial
dimensions \cite{CiCP-25-669}. The inhomogeneous Navier-Stokes equations play a central
role in the description of geophysical flows, including rivers and
shallow coastal currents. These flows are incompressible, but the
density varies due to the presence of suspended sediment. In the
present work we consider the following two-dimensional
inhomogeneous Navier-Stokes equations
\begin{align}
	\partial_x u + \partial_z w & = 0, \label{eq:2D NSE system_1}\\
	\partial_t (\rho u)+ \partial_x (\rho u^2) + \partial_z (\rho uw)
	&= -\partial_x p + \partial_x \sigma_{xx} + \partial_z \sigma_{xz},
	\label{eq:2D NSE system_2} \\
	\partial_t (\rho w)+ \partial_x (\rho uw) + \partial_z (\rho w^2)
	&= -\partial_z p + \partial_x \sigma_{zx} + \partial_z \sigma_{zz}
	- \rho g,
	\label{eq:2D NSE system_3} \\
	\partial_t \rho + \partial_x (\rho u) + \partial_z (\rho w)
	&= 0,
	\label{eq:2D NSE system_4}
\end{align}
where $u$ and $w$ denote the horizontal and vertical velocity
components, respectively, $g$ is the gravitational acceleration, and
$\sigma_{xx},\sigma_{xz},\sigma_{zx},\sigma_{zz}$ are the components of
the deviatoric stress tensor.

Due to the water-sediment mixture, the density $\rho$ is not constant and is
written as
\begin{equation}\label{density equation}
	\rho(t,x,z)= \rho_w + c(t,x,z)\,(\rho_s-\rho_w),
\end{equation}
where $c(t,x,z)\in[0,1]$ is the local volumetric sediment concentration,
and $\rho_w$ and $\rho_s$ are the densities of water and sediment,
respectively. The mixture is assumed to be incompressible and satisfies
$\nabla\cdot \boldsymbol{u}=0$. The density variations are caused by the variations of the concentration $c$.

In principle, the density $\rho$ and the concentration $c$ both depend on the vertical coordinate $z$, and this dependence should be taken into account when deriving any reduced model. In the present work, we restrict our attention to dilute to moderately dilute suspensions, where the vertical variations of the density $\rho$ are small compared with the overall vertical variation of the velocity $u$. Accordingly, in the shallow water reduction we let the density depend only on the depth-averaged concentration $c_m$ by using
\begin{equation}\label{density average equation}
	\rho(t,x) = \rho_w + c_m(t,x)\,(\rho_s-\rho_w),
\end{equation}
where $c_m(t,x)$ denotes the depth-averaged volumetric sediment concentration. With a slight abuse of notation, we use the same symbol $\rho$ for this averaged density in the reduced model. The near-bed concentration $c_b$ required in the deposition closure will be prescribed separately from $c_m$ through the empirical relation given in \secref{bedflux}.

\par Under the assumption of shallowness, we adopt the hydrostatic approximation for the pressure. With the density approximation \eqref{density average equation}, the two-dimensional inhomogeneous Navier-Stokes system
\eqref{eq:2D NSE system_1}-\eqref{eq:2D NSE system_4} reduces to
\begin{align}
	\partial_x u + \partial_z w & = 0, \label{reference system21}\\
	\partial_t u + \partial_x (u^2) + \partial_z (u w)
	&= - \frac{1}{\rho}\,\partial_x p
	+ \frac{1}{\rho}\,\partial_z \sigma_{xz}, \label{reference system22}\\
	\partial_t c + \partial_x (c u) + \partial_z (c w)
	&= 0.  \label{reference system23}
\end{align}
Consistently with the hydrostatic approximation, the pressure is
expressed as
\begin{equation}\label{pressure equation}
	p(t,x,z)
	= \rho(t,x)\,g\big(h(t,x)+ h_b(t,x) - z\big),
\end{equation}
where $h(t,x):=h_s(t,x)-h_b(t,x)$ is the water depth, $h_b$ denotes the
bed elevation and $h_s$ the free-surface elevation. The flow domain is
bounded below by the bottom topography $h_b(t,x)$ and above by the free
surface $h_s(t,x)$. In this formulation, the dominant feedback of the suspended sediment on the hydrodynamics occurs through the mixture density $\rho=\rho(c_m)$ in the hydrostatic pressure \eqref{pressure equation} and the buoyancy terms in the depth-averaged momentum and moment equations, while possible dependence of the viscous and turbulent stresses on the concentration $c$ is neglected. Note that the system does not transport the sediment concentration as a passive scalar only, but couples it back to the hydrodynamic equations through the density.

\par The boundary conditions for the system include kinematic boundaries applied both at the free surface and the bottom, which are defined as follows 
\begin{align}
	\partial_t h_s + \left.u \right\vert_{z=h_s}\partial_x h_s-\left. w \right\vert_{z=h_s}	& = 0, \label{bc_z_cord_1}\\
	\partial_t h_b + \left.u \right\vert_{z=h_b}\partial_x h_b-\left. w \right\vert_{z=h_b}& = -\, F_b, \label{bc_z_cord_2}
\end{align}
where the surface flux is zero, and the bed exchange rate $F_b$, which couples the model to the morphology, will be explained in \secref{bedflux}. Following classical erosion-deposition modelling, we represent the exchange by a net interfacial flux, as commonly used in depth-averaged formulations (e.g.\cite{armanini1988one, galappatti1985depth}). Existing shallow water moment models simplify the problem by neglecting the bed exchange rate, i.e., by taking $F_b
=0$ \cite{Garres,CiCP-25-669}.
\par Additionally, we assume that the remaining shear stress component of the deviatoric stress tensor, \( \sigma_{xz} \), vanishes at the free surface and follows Manning type friction \( \tau \), at the bottom \cite{Garres}. As a result, the boundary conditions for the shear stress are
\begin{equation}\label{shear_stress_z_cord}
	\left.\sigma_{xz}\right\vert_{z=h_s}=0 \quad \quad \text{and} \quad \quad \left.\sigma_{xz}\right\vert_{z=h_b}=\, \tau.
\end{equation}
While the system in equations \eqrefc{reference system21}–\eqrefc{reference system23} already defines our reference system, we aim to reduce its dimension and derive a reduced model by adopting the so-called moment approach outlined in \cite{CiCP-25-669}. This moment approach is based on two main ideas: first, mapping to a reference (\(\zeta-\)coordinate) system and second, expanding the horizontal velocity using a moment expansion. To derive our coupled model, we adopt the approach of \cite{CiCP-25-669} but adapt it to our specific setting, which incorporates coupling with the morphology discussed later in \secref{morphology}. 
\subsubsection{Mapped reference system and moment expansion}
\begin{figure}[H] 
	\centering
	\includegraphics[width=0.95\textwidth]{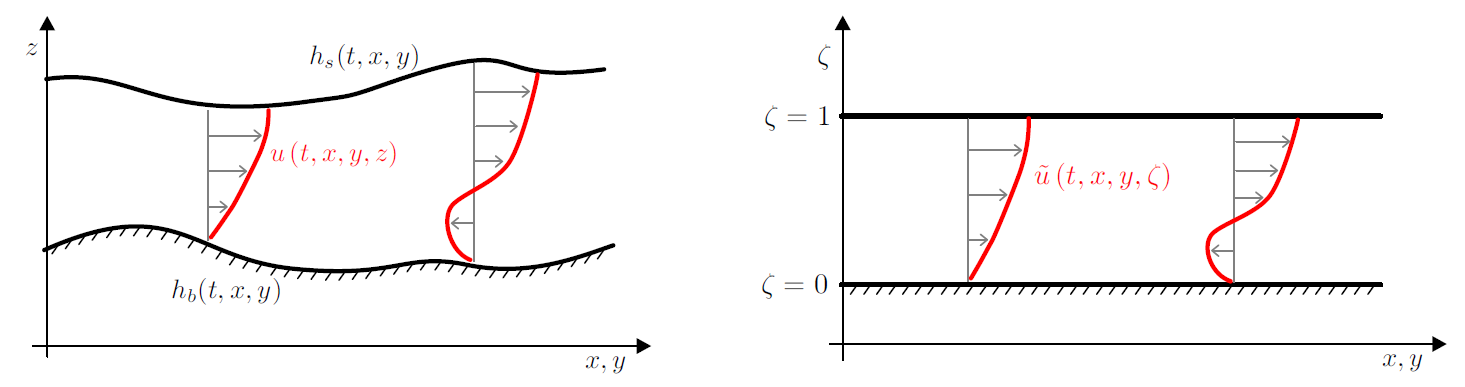}
	\caption{The mapping from physical $z-$space to transformed $\zeta-$space \cite{CiCP-25-669}.}
	\label{fig:mapping_graph}
\end{figure}

According to the first main idea of the Shallow Water Moment model (SWM) approximation described in \cite{CiCP-25-669} we introduce a scaled vertical variable $\zeta(t,x)$, which transforms the $z-$coordinate from the physical space $z\in \left[h_b,h_s\right]$ to the mapped space $\zeta\in \left[0,1\right]$, as shown in \figref{fig:mapping_graph}. The definition of the mapping to $\zeta-$coordinates is given by 
\begin{equation}\label{mapping}
	\zeta(t,x) = \frac{z-h_b(t,x)}{h_s(t,x)-h_b(t,x)}= \frac{z-h_b(t,x)}{h(t,x)}.
\end{equation}
Using \eqrefc{mapping}, for any function $\Phi(t,x,z)$, the corresponding mapped function in $\zeta-$coordinates is given by
\begin{equation}
	\Tilde{\Phi}(t,x,\zeta) = \Phi(t,x,z)\zeta+ h_b(t,x).
\end{equation}
The corresponding differential operators read
\begin{equation}\label{diffop1}
	\partial_\zeta \Tilde{\Phi} = h\partial_z\Phi \quad \text{and} \quad h\partial_s\Phi = \partial_s(h\Tilde{\Phi})-\partial_\zeta(\partial_s(\zeta h+h_b)\Tilde{\Phi}), \quad \text{for}\,\,\, s \in \left[t,x\right]. 
\end{equation}
Taking into account the mapping \eqrefc{mapping} and using the differential operators \eqrefc{diffop1}, \textcolor{black}{the complete vertically resolved system} can be written as
\begin{align}
	\partial_x \left(h \Tilde{u}\right)+\partial_\zeta \left[\Tilde{w}-\Tilde{u}\partial_x(\zeta h+h_b) \right]   &= 0, \label{eq:resolved  system1_1}\\ 
	\begin{split}
		\partial_t (h \Tilde{u})+ \partial_x(h \Tilde{u}^2) +  \partial_\zeta \left[\Tilde{u}\left(\Tilde{w} -\Tilde{u}\partial_x\left(\zeta h+h_b\right)- \partial_t\left(\zeta h+h_b\right)\right)\right] \\
		+gh\partial_x (h + h_b) + \frac{gh^2}{\rho}(1-\zeta)\partial_x \rho & = \frac{1}{\rho}\partial_\zeta \Tilde{\sigma}_{xz},
	\end{split} \label{eq:resolved  system1_2} \\
	\partial_t (h \Tilde{c})  + \partial_x \left(h \Tilde{c}\Tilde{u}\right) +\partial_\zeta \left[\Tilde{c} \left(\Tilde{w}- \partial_x(\zeta h+h_b)\Tilde{u}-\partial_t(\zeta h+h_b) \right) \right]&=0. \label{eq:resolved  system1_3}
\end{align}
The system of equations \eqrefc{eq:resolved  system1_1} to \eqrefc{eq:resolved  system1_3} is referred to as the vertically resolved system because it incorporates the dependence on the vertical variable $\zeta.$\\
Furthermore, the boundary conditions \eqrefc{bc_z_cord_1} and \eqrefc{bc_z_cord_2} and also the shear stress at the free surface and bottom are mapped to $\zeta-$coordinate as
\begin{align}
	& \partial_t h_s + \left.\Tilde{u} \right\vert_{\zeta=1} \partial_x h_s-\left. \Tilde{w} \right\vert_{\zeta=1} = 0, \quad \quad \quad \quad \quad \quad  \quad \left.\Tilde{\sigma}_{xz}\right\vert_{\zeta=1}=0, \label{bc_zeta_1}\\
	&  \partial_t h_b + \left.\Tilde{u}\right\vert_{\zeta=0}  \partial_x h_b-\left. \Tilde{w} \right\vert_{\zeta=0}= - F_b, \quad \quad \quad \quad \quad \quad \left.\Tilde{\sigma}_{xz}\right\vert_{\zeta=0}=\, \tau. \label{bc_zeta_2}
\end{align}
\par Following the second main idea of the Shallow Water Moment model in \cite{CiCP-25-669} we expand the horizontal velocity in the vertical variable $\zeta$, using the Legendre polynomial expansion 
\begin{equation}\label{moment expansion}
	\Tilde{u}(t,x,\zeta) = u_m(t,x) + \sum_{j=1}^{N} \alpha_j(t,x) \phi_j(\zeta),
\end{equation}
where $u_m(t,x)= \displaystyle \int_{0}^{1}\Tilde{u}(t,x,\zeta) d\,\zeta $ is the mean of the horizontal velocity, and the basis functions $\phi_j(\zeta):[0,1]\rightarrow \mathbb{R} $ are the scaled Legendre polynomials of degree $j \in \mathbb{N}$ defined by
\begin{equation}\label{Legendere polynomial}
	\phi_j(\zeta) = \frac{1}{j!}\frac{d^j}{d\zeta^j}(\zeta-\zeta^2)^j.
\end{equation}
Therefore, the first scaled and shifted basis functions are
\begin{equation*}
	\phi_0(\zeta)=1,\quad
	\phi_1(\zeta)= 1-2\zeta, \quad 
	\phi_2(\zeta)= 1-6\zeta+6\zeta^2, \quad  \phi_3(\zeta)= 1-12\zeta+30\zeta^2-20\zeta^3.
\end{equation*}
These polynomials also satisfy the condition $\phi_j(0)=1$, and $\displaystyle \int_{0}^{1} \phi_j(\zeta)d\zeta = 0.$ In the moment expansion \eqrefc{moment expansion}, $N \in \mathbb{N}$ denotes both the highest degree of the Legendre polynomial and the order of the velocity expansion. Generally, a larger value of $N$ represents more complex flows. The corresponding basis coefficients $\alpha_j(t,x)$ for $j \in \left\{1,\ldots, N\right\}$ are unknown and an additional $N$ equations are derived by projecting momentum balance equation onto Legendre polynomials \cite{CiCP-25-669}. These basis coefficients represent the so-called moments and encapsulate the information of horizontal velocity in the vertical $\zeta-$direction.

\par By employing the moment approach outlined above within our specific framework, we derive the evolution equations for the extended set of variables \( (h, hu_m, h\alpha_i, h c_m)^T \) for \( i \in \{1, \ldots, N\} \). The evolution equations for the water height \( h \), the mean horizontal velocity \( u_m \), and the volumetric sediment concentration \( c_m \) are obtained through the following Galerkin projection of \eqrefc{eq:resolved system1_1}, \eqrefc{eq:resolved system1_2}, and \eqrefc{eq:resolved system1_3}, respectively
\begin{equation}\label{Galerkin projection}
	\langle \cdot, \, 1 \rangle = \int_{0}^{1} \cdot \, \, d\zeta.
\end{equation} \\
Furthermore, the \( N \) evolution equations for the moments are derived by employing a higher-order Galerkin projection of the momentum balance equation \eqrefc{eq:resolved system1_2}
\begin{equation}\label{Higher-order Galerkin projection}
	\langle \cdot, \, \phi_i \rangle = \int_{0}^{1} \cdot \, \phi_i(\zeta) \, d\zeta, \quad \quad i=1,\ldots, N.
\end{equation} 
i.e. multiplying equation \eqrefc{eq:resolved system1_2} by the test function \( \phi_i(\zeta) \) and subsequently integrating with respect to \( \zeta. \)  
\par In the three subsequent sections, we provide a detailed derivation of the depth-averaged equations for \( h \), \( hu_m \), and \( h\alpha_i \). The corresponding derivation for \( h c_m \) is presented within the morphodynamic \secref{concentration_avg} since morphodynamic change employs suspended sediment concentration.  
\subsubsection{Depth-averaging the mass balance}
To recover an explicit expression for $ \Tilde{w}$, equation \eqrefc{eq:resolved  system1_1} can be written in the following integral form 
\begin{equation}\label{vertical velocity}
	\Tilde{w}= - \partial_x \left(h \int_{0}^{\zeta} \Tilde{u} \,d\hat{\zeta}\right)+ \Tilde{u}\partial_x(\zeta h+h_b),
\end{equation}
Now, to derive the standard depth-averaged mass balance equation of the shallow water system, we depth-average the transformed mass balance equation, i.e. we apply a projection \eqrefc{Galerkin projection} to an equation \eqrefc{vertical velocity} and use the kinematic boundary conditions \eqrefc{bc_zeta_1} - \eqrefc{bc_zeta_2}, which yields the following depth-averaged mass balance equation.
\begin{equation}\label{depth-average mass balance}
	\partial_t h + \partial_x (hu_m) = F_b.
\end{equation}
Equation \eqrefc{depth-average mass balance} represents an equation for the water height $h$, and couples to the equation for the discharge $hu_m$, (\secref{momentum_discharge}) and the bed exchange rate $F_b$ (\secref{bedflux}). Without the bed exchange $F_b$, equation \eqrefc{depth-average mass balance} simplifies to the same as in the SWM \cite{CiCP-25-669}. 
\subsubsection{Depth-averaging the momentum equation}\label{momentum_discharge}
The depth-averaged momentum balance equation along the horizontal $x-$axis is derived by applying \eqrefc{Galerkin projection} to \eqrefc{eq:resolved system1_2} and is given by 	
\begin{align}\label{depth_avg_momentum}   
	\begin{split}
		\partial_t(h u_m) +\partial_x \left( h \left(u_m^2+ \sum_{j=1}^{N}  \frac{\alpha_j^2}{2j+1}\right)+ \frac{gh^2}{2} \right)
		&=  - gh \partial_x h_b + F_b u_b \\
		&- \frac{gh^2}{2 \rho}(\rho_s-\rho_w)\partial_x c_m- \epsilon |u_b|u_b.
	\end{split}
\end{align}
We provide the detailed derivation in \appref{averaged momentum balance}.
\par Here we only discuss the last three terms on the right-hand side of \eqrefc{depth_avg_momentum}, as they are new compared to the momentum balance equation in \cite{CiCP-25-669}:
\begin{enumerate}
	\item The term $F_b u_b$ quantifies momentum transfer between flow and erodible bed due to sediment exchange.
	\item The term $\displaystyle \frac{gh^2}{2 \rho}(\rho_s-\rho_w)\partial_x c_m$ represents the effects of spatial variations in sediment concentration.
	\item The term $\epsilon |u_b|u_b$ associates friction within the momentum balance equation. 
	We consider Newtonian Manning type friction at the bottom \cite{Garres} and Newtonian friction within the fluid  \parencite{ Garres, CiCP-25-669}. Since the Galerkin projection of the momentum balance is influenced only by bottom friction and remains independent of friction within the fluid, we present only the bottom friction in this section, while the treatment of friction within the fluid will be addressed in the subsequent section.\\
	Accordingly, the Manning type friction at the bottom is given by
	\begin{equation*}
		\left. \Tilde{\sigma}_{xz}(\zeta) \right\vert_{\zeta=0}= -\rho \epsilon |u_b|u_b,
	\end{equation*}
	where the bottom velocity is expressed as
	\begin{equation}\label{bottom_velocity}
		u_b=\left. \Tilde{u}\right\vert_{\zeta=0}=u_m+\displaystyle \sum_{j=1}^{N}\alpha_j(t,x)\, \phi_j(0)=u_m+\sum_{j=1}^{N}\alpha_j(t,x),
	\end{equation} 
	and $\epsilon$ is a dimensionless constant defined in \cite{Garres}.
\end{enumerate}

\subsubsection{Higher-order moments of the momentum equation}\label{Higher order averages_1}
The evolution equations for the higher-order moments of the momentum equation are derived from the higher-order Galerkin projection \eqrefc{Higher-order Galerkin projection} of \eqrefc{eq:resolved system1_2}. 
Thus, the resulting equation for $h\alpha_i$ is given by
\begin{align}\label{final_higher_average_equation}
	\begin{split}
		\partial_t (h\alpha_i)&+ \partial_x \left(2hu_m\alpha_i + h\sum_{j,k=1}^{N}\alpha_j\alpha_k A_{ijk}\right)\\
		&= u_m\partial_x (h \alpha_i) -  \sum_{j,k=1}^{N} B_{ijk}\alpha_k \partial_x (h \alpha_j)
		+F_b\left(\alpha_i+\sum_{j=1}^{N}\alpha_j (H_{ij}- G_{ij})\right)\\
		&- \sum_{j=1}^{N}\alpha_j G_{ij}\frac{1}{1-\psi}\partial_x Q_b
		+(2i+1)\frac{gh^2}{\rho} (\rho_s-\rho_w)\frac{\partial c_m}{\partial x} K_i - (2i+1)\left( \epsilon |u_b| u_b + \frac{\nu}{h}\sum_{j=1}^{N} \alpha_j C_{ij} \right),
	\end{split}
\end{align}
where $A_{ijk},\, B_{ijk},\, G_{ij},\, H_{ij}$ and $K_i$ introduced here, are constant coefficients and associated with the integrals of the scaled Legendre polynomials.
\begin{equation}\label{constant matrices}
	\begin{aligned}
		& A_{ijk}= (2i+1)\int_{0}^{1}\phi_i\phi_j\phi_k  \, d\zeta, \quad  B_{ijk}= (2i+1)\int_{0}^{1}\phi_i^\prime\left(\phi_j d\hat{\zeta} \right)\phi_k   \, d\zeta, \quad  C_{ij}= \int_{0}^{1} \phi_i^\prime \phi_j^\prime  \, d\zeta,\\
		& G_{ij}=  (2i+1)\int_{0}^{1} \phi_i\phi_j^{\prime} \,d \zeta, \quad  H_{ij}= (2i+1) \int_{0}^{1} \zeta \phi_i\phi_j^{\prime} \, d\zeta, \quad  K_i= \int_{0}^{1}\zeta \phi_i  \, d\zeta.
	\end{aligned}
\end{equation}
We provide the detailed derivation in \appref{Higher order averages}.
\par On the right-hand side of equation \eqrefc{final_higher_average_equation}, four additional terms appear compared to the classical SWM \cite{CiCP-25-669}: 
\begin{enumerate}
	\item The term $F_b\left(\alpha_i+\sum_{j=1}^{N}\alpha_j (H_{ij}- G_{ij})\right)$ represents a source term associated with erosion and deposition effects at the bottom, governed by the bed exchange rate $F_b.$
	\item The term $\sum_{j=1}^{N}\alpha_j G_{ij}\frac{1}{1-\psi}\partial_x Q_b,$ which includes the sediment discharge term $Q_b,$ represents the influence of the bedload transport, as defined in \secref{morphology}. This term appears in the equations for the higher-order moments since our model does not assume $\partial_t h_b = 0.$
	\item The term $(2i+1)\frac{gh^2}{\rho} (\rho_s-\rho_w)\frac{\partial c_m}{\partial x} K_i$ arises from the higher order Galerkin projection on the spatial variation of sediment concentration. This term signifies the interaction between suspended sediment and momentum balance.
	\item The term $(2i+1)\left( \epsilon |u_b| u_b + \frac{\nu}{h}\sum_{j=1}^{N} \alpha_j C_{ij} \right)$ results from the higher-order Galerkin projection of Manning type friction at the bottom i.e., $\left. \Tilde{\sigma}_{xz}(\zeta) \right\vert_{\zeta=0}= -\rho \epsilon |u_b|u_b,$ and Newtonian friction within the fluid i.e.,  $\left. \Tilde{\sigma}_{xz}(\zeta) \right\vert_{\zeta \in (0,1)}= - \frac{\mu}{h} \partial_\zeta \Tilde{u}(\zeta)$, and the projection reads  
	\begin{align*}
		\int_{0}^{1} \phi_i \partial_\zeta \Tilde{\sigma}_{xz}\, d \zeta  = \int_{0}^{1} \partial_\zeta (\phi_i \Tilde{\sigma}_{xz}) \,d\zeta  -  \int_{0}^{1} \Tilde{\sigma}_{xz} \partial_\zeta \phi_i \,d\zeta
		& = - \left. \phi_i(\zeta) \Tilde{\sigma}_{xz}\right\vert_{\zeta=0} - \frac{\mu}{h}\int_{0}^{1} \phi_i^\prime  \partial_\zeta \Tilde{u} \, d\zeta\\
		& =  - \rho \epsilon |u_b| u_b  - \frac{\mu}{h} \sum_{j=1}^{N} \int_{0}^{1} \alpha_j \phi_i^\prime \phi_j^\prime \,d\zeta\\
		&=- \rho \epsilon |u_b| u_b - \frac{\mu}{h}\sum_{j=1}^{N} \alpha_j C_{ij}. 
	\end{align*}    
\end{enumerate}
\begin{remark}
	We use a classical solid transport discharge laws, the Meyer-Peter$\And$Müller formula, but using the velocity moments and define it in terms of the bottom velocity \eqref{bottom_velocity} instead of the averaged velocity. A similar procedure was done in the existing literature \cite{Garres} and for the simpler friction models in \cite{CiCP-25-669}. However, in the classical literature on shallow water equations, the solid transport formula is based on a friction law which is usually calibrated for the average velocity, not for the bottom velocity. By using the bottom velocity for the solid transport discharge law, the model uses the velocity-profile information provided by the moment expansion in the transport terms. The calibration of the corresponding empirical coefficients is not addressed in this work and remains a topic for further study.
\end{remark}
\subsection{Morphodynamic equations}\label{morphology}
The morphodynamic part of the model covers the suspended and bedload transport. It consists of the depth-averaged volumetric sediment concentration equation and the Exner equation, as explained below.
\subsubsection{Depth-averaged volumetric sediment concentration}
\label{concentration equation}\label{concentration_avg}
In this section, we derive the depth-averaged equation for the suspended sediment concentration. In the mapped $\zeta$-coordinates, the vertically resolved concentration equation \eqrefc{eq:resolved  system1_3} can be written in conservative form as
\begin{equation}\label{eq:vertical_concentration_conservative}
	\partial_t(h\tilde{c})+
	\partial_x(h\tilde{c}\tilde{u})+
	\partial_\zeta J=0,
\end{equation}
where the vertical sediment flux $J$ is given by
\begin{equation}\label{eq:vertical_sediment_flux}
	J(t,x,\zeta)=\tilde{c}(t,x,\zeta)
	\left[\tilde{w}-\partial_x(\zeta h+h_b)\tilde{u}-\partial_t(\zeta h+h_b)\right].
\end{equation}
We define the depth-averaged volumetric sediment concentration $c_m$ and the depth-averaged concentration flux $	h\langle \tilde{c}\tilde{u}\rangle$, respectively, by
\begin{equation}\label{eq:depth_averaged_concentration_def}
	c_m(t,x):=\int_0^1 \tilde{c}(t,x,\zeta)\,d\zeta,
	\qquad
	h\langle \tilde{c}\tilde{u}\rangle :=
	h\int_0^1 \tilde{c}(t,x,\zeta)\tilde{u}(t,x,\zeta)\,d\zeta.
\end{equation}
Integrating \eqrefc{eq:vertical_concentration_conservative} over $\zeta\in[0,1]$ and using the kinematic boundary conditions \eqrefc{bc_zeta_1}-\eqrefc{bc_zeta_2} gives no sediment flux through the free surface, $J(1)=0$, while the bed flux is prescribed by $J(0)=E-D$. Here, the erosion and deposition rates $E$ and $D$ will be defined in \secref{bedflux}. Therefore,
\begin{equation}\label{eq:unclosed_concentration_equation}
	\partial_t(hc_m)+\partial_x\left(h\langle\tilde{c}\tilde{u}\rangle\right)=E-D.
\end{equation}
Since the present model evolves only the depth-averaged concentration $c_m$, the depth-averaged concentration flux  $h\langle \tilde{c}\tilde{u}\rangle$ is not closed. We therefore use the shallow-water closure
\begin{equation}\label{eq:concentration_flux_closure}
	\langle\tilde{c}\tilde{u}\rangle \approx c_m u_m .
\end{equation}
Thus, using \eqrefc{eq:concentration_flux_closure} the depth-averaged volumetric sediment concentration equation \eqrefc{eq:unclosed_concentration_equation} reads
\begin{equation}\label{eq:depth_averaged_concentration_equation}
	\partial_t(hc_m)+\partial_x(hc_m u_m)=E-D.
\end{equation}
The closure \eqrefc{eq:concentration_flux_closure} is a modelling simplification: it neglects the vertical covariance between concentration and velocity. A more detailed extension would require additional information on the vertical structure of the concentration, for example through separate evolution equations for concentration moments. This extension is not used in the present reduced model and is left for future work.
The near-bed concentration $c_b$ enters only the deposition closure. Following \cite{bradford1999hydrodynamics}, we prescribe $c_b$ algebraically from the depth-averaged concentration through $c_b=S_b c_m$, where $S_b$ is the empirical near-bed concentration factor defined in \eqrefc{bed_concentration}. Thus, $c_m$ represents the concentration averaged over the water column, whereas $c_b$ is the corresponding near-bed value used to evaluate deposition. Note that the suspended concentration is not treated as a passive scalar in the coupled SWEMED system. The depth-averaged concentration $c_m$ enters the mixture density $\rho=\rho(c_m)$ in \eqrefc{density average equation}, which appears in the hydrostatic pressure \eqrefc{pressure equation} and generates buoyancy contributions in the depth-averaged momentum and moment equations.

\subsubsection{Bedload mass balance equation}
We use the Exner equation  \parencite{exner1925uber,exner1920physik} to model the evolution of bedload transport, incorporating the principle of mass conservation for the sediment layer. This formulation represents sediment transport through a flux term, and the right-hand side of the equation includes a source term that describes erosion and deposition processes, which reads as 
\begin{equation}\label{Exner equation}
	\partial_t h_b  + \frac{1}{1-\psi}\partial_x Q_b = -F_b,
\end{equation}
where $h_b$ is the bed elevation, $Q_b$ is the solid sediment discharge, and $F_b$ is the bed exchange rate, which accounts for erosion and deposition processes, and defined in \secref{bedflux}. More explicitly, the Exner equation represents the rate of bed deformation. Active sediment transport and rapid bed deformation will occur if the flow entrains more sediment particles than it deposits.
\par The result of the Exner equation \eqrefc{Exner equation} is highly dependent on the choice of the sediment discharge formula, $Q_b$. We note that there is a plethora of formulas for the sediment discharge term. In this paper, we adopt the following formula \cite{Garres}
\begin{equation}\label{sediment discharge}
	Q_b = sgn(\tau)\, Q \, \Phi(\theta), 
\end{equation}
where $Q=\sqrt{\left(\frac{\rho_s}{\rho_w}-1\right)g\,d_s^3}$ is the characteristic discharge and $sgn(\cdot)$ is the sign function.	
We chose the expression \eqrefc{sediment discharge} for $Q_b$ since it depends on the bottom shear stress, used in its non-dimensional form $\theta$ and $Q_b$ can be expressed as a function of $\theta$ as $\Phi(\theta).$  The dimensionless bottom shear stress $\theta$ is also called Shields parameter \cite{cordier2011bedload} and defined by
\begin{equation}\label{Shields parameter}
	\theta = \frac{|\tau|}{g\,(\rho_s- \rho_w)\,d_s}= \frac{\rho \epsilon |u_b| u_b}{g\,(\rho_s- \rho_w)\,d_s},
\end{equation}
where $\epsilon$ is a dimensionless constant \cite{garres2023general}, $u_b$ is the velocity at bottom, as defined in  \eqrefc{bottom_velocity}, and $d_s$ is the diameter of a sediment particle.
\par There is a framework that includes many different formulas for $\Phi(\theta)$ (see  \parencite{camenen2006phase,fernandez1976erosion,meyer1948formulas,nielsen1992coastal,ribberink1994sediment}). In this work, we employ the Meyer-Peter $\And$ Müller formula \cite{meyer1948formulas}
\begin{equation}\label{eq:MPM_bed_closure}
	\Phi(\theta)= 8(\theta-\theta_{c})_{+}^{3/2},
\end{equation} 
where $(\cdot)_+$ is the positive part and $\theta_c$ is the critical shear stress. Notably, sediment transport is initiated only when the Shields parameter $\theta$ exceeds the critical threshold $\theta_c$.
\subsubsection{Morphological conditions}\label{bedflux}
The bed exchange rate $F_b$, used in \secref{Hydrodynamic} and \secref{morphology}, represents a dynamic equilibrium between two effects:
(1) the sediment deposition due to gravitational force and (2) the sediment erosion from the interface due to turbulence. This balance is crucial for understanding sedimentary processes and is defined by
\begin{equation}\label{bedflux_eq}
	F_b = \frac{E-D}{1-\psi},
\end{equation}  
Here $E$ and $D$ denote the erosion and deposition rates, respectively, and $\psi$ is the bedload porosity.    
\vspace{2mm}\\
To close the system, we refer to \parencite{bradford1999hydrodynamics,del2023lagrange,garcia1993experiments,gonzalez2020robust,zhang1993sedimentation} and take the following formulas for erosion and deposition.\\
The sediment erosion and deposition follow from \parencite{del2023lagrange,gonzalez2020robust}, as given by
\begin{equation}
	E = \omega_o (1-\psi) E_s,
\end{equation}
and 
\begin{equation}
	D= \omega_o c_b,
\end{equation} 
where
\begin{itemize}
	\item[-] $\omega_o$ is the settling velocity determined by \cite{zhang1993sedimentation}, 
	\begin{equation}\label{settling velocity}
		\omega_o = \sqrt{\left(\frac{13.95\nu_w}{d_s}\right)^2+1.09\, 	\rho_w \left(\frac{\rho_s}{\rho_w}-1\right) gd_s} - \frac{13.95\nu_w}{d_s},
	\end{equation}
	\item[-] $\nu_w$ is the kinematic viscosity of water and $d_s$ is the diameter of sediment particle,
	\item[-] $E_s$ is the sediment erosion coefficient and computed by \cite{garcia1993experiments},
	\begin{equation}
		E_s = \frac{1.3\times 10^{-7} \mathcal{Z}^5}{1+4.3\times 10^{-7} \mathcal{Z}^5},
	\end{equation} 
	\item[-]  $\mathcal{Z} = \displaystyle  \frac{\gamma_1\sqrt{c_D} |u_b| }{\omega_o}\mathcal{R}^{\gamma_2}_{p}$ with particle Reynolds number	$\mathcal{R}_p = \displaystyle \frac{\sqrt{(\rho_s-\rho_w) g d_s}d_s}{\nu_w}$, $u_b$ the bottom velocity, and $c_D$ the bed drag coefficient,
	\item[-] $\gamma_1,\gamma_2$ are two parameters depending on $\mathcal{R}_p$ \cite{gonzalez2020robust}
	$$
	\left(\gamma_1,\gamma_2\right)=\begin{cases}
		\left(1,0.6\right), & \text{if $\mathcal{R}_p> 2.36$}\\
		\left(0.586,1.23\right), & \text{if $\mathcal{R}_p \le 2.36$}
	\end{cases}
	$$
	\item[-] $c_b$ is the near-bed volumetric sediment concentration used in the deposition law. Since the model transports only the depth-averaged concentration $c_m$, we follow \cite{bradford1999hydrodynamics} and prescribe $c_b$ algebraically as
		\begin{equation}\label{bed_concentration}
			c_b=c_m(t,x)\,S_b,
			\qquad
			S_b=0.4\left(\frac{d_s}{D_{sg}}\right)^{1.64}+1.64,
		\end{equation}
		where $S_b$ is the empirical near-bed concentration factor and $D_{sg}$ is the geometric mean size of the suspended sediment mixture. In this work, all sediment particles have the same diameter, so that $D_{sg}=d_s$ and $S_b=2.04$. 
\end{itemize}
\subsection{Coupled hydro-morphodynamic model}
By incorporating the derivations and definitions outlined in \secref{Hydrodynamic} and \secref{morphology}, we obtain a closed Shallow Water Exner Moment model with Erosion and Deposition (SWEMED). In contrast to the existing HSWEM in \cite{Garres}, SWEMED (i) introduces a sediment concentration equation, (ii) couples the variable sediment–water mixture density with the momentum equation and higher-order moments, and (iii) includes additional source terms arising from erosion and deposition. The resulting system consists a total of $N+4$ coupled equations \eqref{complete system} for the set of convective variables $(h,hu_m,h\alpha_1,\ldots,h\alpha_N,h c_m,h_b)^T\in \mathbb{R}^{N+4}$, where $N$ denotes the highest degree of the Legendre polynomial expansions which is also called the order of the moment model.
\begin{equation*}
	\underbrace{1}_{\text{mass balance}} + \underbrace{1}_{\text{momentum equation}} +\underbrace{N}_{\text{moment equations}}+ 
	\textcolor{black}{							\underbrace{1}_{\text{sediment concentration}}}+
	\textcolor{black}{							\underbrace{1}_{\text{bedload mass balance}}}= \, N+4
\end{equation*} 
Thus the coupled model is formulated as follows		
\begin{equation}\label{complete system}
	\left\{
	\begin{aligned}
		&\partial_t h  + \partial_x (h u_m) 
		&&= \frac{E-D}{1-\psi},\\
		&\partial_t(h u_m)  +\partial_x \left( h \left(u_m^2+ \sum_{j=1}^{N}  \frac{\alpha_j^2}{2j+1}\right) + \frac{gh^2}{2} \right)
		&&= - gh \partial_x h_b - \frac{gh^2}{2 	\rho}(\rho_s-\rho_w)\partial_x c_m \\
		&&&\quad + \frac{(E-D) u_b}{1-\psi}- \epsilon |u_b|u_b,\\
		&\partial_t (h\alpha_i)+ \partial_x \left(2hu_m\alpha_i+ h\sum_{j,k=1}^{N}\alpha_j\alpha_k A_{ijk}\right)  
		&&= u_m\partial_x (h \alpha_i) -  \sum_{j,k=1}^{N} B_{ijk}\alpha_k \partial_x (h \alpha_j) \\
		&&&\quad + (2i+1) \frac{gh^2}{\rho} (\rho_s-\rho_w)\frac{\partial c_m}{\partial x} K_i \\
		&&&\quad +\frac{E-D}{1-\psi}\left(\alpha_i+\sum_{j=1}^{N}\alpha_j (H_{ij}- G_{ij})\right)\\
		&&&\quad -(2i+1)\left( \epsilon |u_b| u_b + \frac{\nu}{h}\sum_{j=1}^{N} \alpha_j C_{ij} \right)\\
		&&&\quad - \sum_{j=1}^{N}\alpha_j G_{ij}\displaystyle \partial_x  \left(\frac{Q_b}{1-\psi}\right),\\
		&\partial_t(h c_m) + \partial_x (h c_m u_m) 
		&&= E-D,\\
		&\partial_t h_b  + \displaystyle \partial_x \left(\frac{Q_b}{1-\psi}\right) 
		&&= \frac{D-E}{1-\psi}.
	\end{aligned}
	\right.
\end{equation}
We can write the complete system \eqrefc{complete system} as a first-order PDE with conservative, non-conservative, and source terms, which reads
\begin{equation}\label{first-order pde}
	\partial_t W + \frac{\partial F}{\partial W}\partial_x 	\textcolor{black}{W}\,= B\partial_x W - A_s\partial_x W+ S_F+S_{ED}.
\end{equation}
In the left-hand side of \eqrefc{first-order pde}, the jacobian of the conservative flux is given by
\begin{equation*}
	\medmath{
		\frac{\partial F}{\partial W}=\begin{pmatrix} 
			0 & 1 & 0 & \dots & 0 & 0 & 0 \\
			gh-u_m^2- \displaystyle \sum_{j=1}^{N}\frac{\alpha_j^2}{2j+1} 	& 2u_m & \frac{2\alpha_1}{2\cdot 1+1} & \dots & \frac{2\alpha_N}{2\cdot N+1} & 0 & 0 \\
			-2u_m\alpha_1- \displaystyle 	\sum_{j,k=1}^{N}A_{1jk}\alpha_j\alpha_k & 2\alpha_1 & 2u_m+2 \displaystyle \sum_{k=1}^{N}A_{11k}\alpha_k & \dots & 2 \displaystyle \sum_{k=1}^{N}A_{1Nk}\alpha_k & 0 & 0 \\
			\vdots & \vdots & \hspace{3cm} \ddots   &  & & \vdots& \vdots\\
			\vdots & \vdots & & \hspace{-1mm} \ddots  &  & \vdots& \vdots\\
			-2u_m\alpha_N-\displaystyle 	\sum_{j,k=1}^{N}A_{Njk}\alpha_j\alpha_k & 2\alpha_N & 2 \displaystyle \sum_{k=1}^{N}A_{N1k}\alpha_k & \dots & 2u_m+2 \displaystyle \sum_{k=1}^{N}A_{NNk}\alpha_k & 0 & 0 \\
			-c_m u_m & -c_m  & 0  & \ldots  & 0  & u_m  & 0  \\
			\delta_h & \delta_q & \delta_q & \dots & \delta_q & \delta_c & 	0 
		\end{pmatrix}, }
\end{equation*}
where
\begin{subequations}
	\begin{equation}\label{eq:delta_q}
		\begin{aligned}
			\delta_q 	=\frac{\partial}{\partial_{hu_m}}\left(\frac{Q_b}{1-\psi}\right) = 	\frac{24Q}{1-\psi}sgn(\tau)\frac{\rho \epsilon}{g(\rho_s-\rho_w)d_s}\left(\theta-\theta_c\right)_{+}^{1/2}\frac{u_b}{h},
		\end{aligned}
	\end{equation}    
	\begin{equation}\label{eq:delta_h}
		\delta_h= 	\frac{\partial}{\partial_h}\left(\frac{Q_b}{1-\psi}\right)= -u_b\left(1+\displaystyle \frac{c_m(t,x)\left(\rho_s-\rho_w\right)}{2\rho}\right)\delta_q,
	\end{equation}
	\begin{equation}\label{eq:delta_c}
		\delta_c=\frac{\partial}{\partial_{h 	c_m}}\left(\frac{Q_b}{1-\psi}\right)= u_b\displaystyle \left(\frac{\rho_s-\rho_w}{2\rho}\right)\delta_q,
	\end{equation}
	\quad \quad \text{and}
	\begin{equation}\label{eq:delta_alpha}
		\frac{\partial}{\partial_{h 	\alpha_i}}\left(\frac{Q_b}{1-\psi}\right) =  \delta_q \quad \quad i \in \left\{1, \ldots, N\right\},
	\end{equation}	
\end{subequations}
with bottom velocity $u_b$ and the Shields parameter $\theta$, defined in \eqrefc{bottom_velocity} and \eqrefc{Shields parameter}, respectively.
\newline
On the right-hand side of \eqrefc{first-order pde}, the non-conservative matrix reads   
\begin{equation*}
	\medmath{
		B=\begin{pmatrix} 
			0 & 0 & 0 & \dots & 0 & 0 & 0 \\
			\displaystyle	\frac{gh(\rho-\rho_w)}{2\rho} & 0 & 0 & \dots 	& 0 & \displaystyle\frac{-gh(\rho_s-\rho_w)}{2\rho} & -gh \\
			\displaystyle	\frac{gh(\rho-\rho_w)}{2\rho} & 0 & u_m- 	\displaystyle \sum_{k=1}^{N}B_{11k}\alpha_k & \dots & - \displaystyle \sum_{k=1}^{N}B_{1Nk}\alpha_k & \displaystyle \frac{-gh(\rho_s-\rho_w)}{2\rho} & 0 \\
			\vdots & \vdots & \hspace{3cm} \ddots   &  & & \vdots& \vdots\\
			\vdots & \vdots & & \hspace{-1mm} \ddots  &  & \vdots& \vdots\\
			0 & 0 & -\displaystyle\sum_{k=1}^{N}B_{N1k}\alpha_k & \dots & 	u_m- \displaystyle \sum_{k=1}^{N}B_{NNk}\alpha_k & 0 & 0 \\
			0 & 0 & 0 & \dots & 0 & 0 & 0 \\
			0 & 0 & 0 & \dots & 0 & 0 & 0 \\
		\end{pmatrix}, } 
\end{equation*}
and the sediment discharge matrix $A_s$ is resulting from $\displaystyle \sum_{j=1}^{N}\alpha_jG_{ij}\partial_x\left(\frac{Q_b}{1-\psi}\right)$ and reads 
\begin{equation*}
	\medmath{
		A_s=\begin{pmatrix} 
			0 & 0 & 0 & \dots & 0 & 0 & 0 \\
			0 & 0 & 0 & \dots & 0 & 0 & 0 \\
			\displaystyle \sum_{j=1}^{N}\alpha_jG_{1j}\delta_h & 	\displaystyle \sum_{j=1}^{N}\alpha_jG_{1j}\delta_q & \displaystyle \sum_{j=1}^{N}\alpha_jG_{1j}\delta_q & \dots &\displaystyle \sum_{j=1}^{N}\alpha_jG_{1j}\delta_q & \displaystyle \sum_{j=1}^{N}\alpha_jG_{1j\delta_c} & 0 \\
			\displaystyle \sum_{j=1}^{N}\alpha_jG_{2j}\delta_h & 	\displaystyle \sum_{j=1}^{N}\alpha_jG_{2j}\delta_q & \displaystyle \sum_{j=1}^{N}\alpha_jG_{2j}\delta_q & \dots & \displaystyle \sum_{j=1}^{N}\alpha_jG_{2j}\delta_q & \displaystyle \sum_{j=1}^{N}\alpha_jG_{2j\delta_c} & 0 \\
			\vdots & \vdots & \hspace{3cm} \ddots   &  & & \vdots& \vdots\\
			\vdots & \vdots & & \hspace{-1mm} \ddots  &  & \vdots& \vdots\\
			\displaystyle \sum_{j=1}^{N}\alpha_jG_{Nj}\delta_h & 	\displaystyle \sum_{j=1}^{N}\alpha_jG_{Nj}\delta_q & \displaystyle \sum_{j=1}^{N}\alpha_jG_{Nj}\delta_q & \dots & \displaystyle \sum_{j=1}^{N}\alpha_jG_{Nj}\delta_q & \displaystyle \sum_{j=1}^{N}\alpha_jG_{Nj\delta_c} & 0 \\
			0 & 0 & 0 & \dots & 0 & 0 & 0 \\
			0 & 0 & 0 & \dots & 0 & 0 & 0 \\    
		\end{pmatrix}. } 
\end{equation*}
The source terms $S_{F,i}$ and $S_{ED,i}$ result  from the combined effects of friction and the processes of erosion and deposition, respectively, and are given by
\begin{equation*}
	S_{F,i}=\begin{cases}
		0,						 & \quad \text{if $i=1,N+3,N+4,$}\\
		-\epsilon  |u_b| u_b,	 & \quad \text{if $i=2,$}\\
		-(2i+1)\displaystyle \left(\epsilon |u_b| u_b + \displaystyle 	\frac{\nu}{h} \displaystyle \sum_{j=1}^{N}\alpha_j C_{ij}\right),			 & \quad\text{if  $i \in \left\{3,\ldots,N+2\right\},$}
	\end{cases}    
\end{equation*} 
and
\begin{equation*}
	S_{ED,i}=\begin{cases}
		\displaystyle \frac{E-D}{1-\psi},		 & \text{if $i=1,$}\\
		\displaystyle \frac{(E-D)u_b}{1-\psi}, & \text{if $i=2,$}\\
		\displaystyle \frac{E-D}{1-\psi}\left(\alpha_i+\displaystyle 	\sum_{j=1}^N \alpha_j(H_{ij}-G_{ij})\right), & \text{if  $i \in \left\{3,\ldots,N+2\right\},$}\\
		E-D,	& \text{if $i=N+3,$}\\
		\displaystyle \frac{D-E}{1-\psi},		 & \text{if $i=N+4.$}\\
	\end{cases}    
\end{equation*} 
The resulting system \eqrefc{first-order pde} can be written in the form
\begin{equation}\label{final SWMEED}
	\partial_t\,W + A_{ED}(W) \partial_x \, W= S_F(W)+S_{ED}(W),
\end{equation}
where we define the state vector of conservative variables as $W= (h,hu_m,h\alpha_1,\ldots,h\alpha_N,h c_m,h_b)^T\in \mathbb{R}^{N+4}$. We express the transport matrix in explicit form as $A_{ED}(W) = \left(\displaystyle \frac{\partial F}{\partial W}-B+ A_s\right) \in \mathbb{R}^{(N+4)\times(N+4)}$. Additionally, the right-hand side of \eqrefc{final SWMEED} models the source term, which arises from friction and the considerations of erosion and deposition.

\subsubsection{Lower-order moments of SWEMED}
In this section, we present the structural framework of lower-order moments of the SWEMED, namely the zeroth-order $(N=0)$, first-order $(N=1)$, second-order $(N=2)$, and third-order $(N=3)$ formulations as representative examples. These formulations are used in the numerical simulations in \secref{sec4}.
\paragraph{\, SWEMED0}\label{zero order system}
Considering $N=0$, the velocity profile is uniform in the vertical direction, i.e., $u(t,x,\zeta)=u_m(t,x)$. The zeroth-order Shallow Water Exner Moment model with Erosion and Deposition (SWEMED0) is identical to the standard Shallow Water Exner model with Erosion and Deposition (SWEED). The SWEMED0 model with state vector $W=(h,hu_m,hc_m,h_b)^T$ can be written in the form of \eqrefc{final SWMEED} using the following definitions
\begin{equation*}
     \quad F= \begin{pmatrix}h\\ \displaystyle hu_m^2 + 	\frac{1}{2}gh^2\\h c_m u_m\\ \displaystyle \frac{Q_b}{1-\psi} \end{pmatrix}, \quad
	\frac{\partial F}{\partial W} = \begin{pmatrix}
		0 & 1 & 0 & 0   \\
		gh- u^{2}_m & 2u_m & 0 & 0  \\
		-c_m u_m & c_m & u_m & 0  \\
		\delta_h & \delta_q &\delta_c & 0
	\end{pmatrix},  
\end{equation*}
\begin{equation*}
	B=  \begin{pmatrix}
		0 & 0 & 0 & 0   \\
		\displaystyle	\frac{gh(\rho-\rho_w)}{2\rho} & 0 	&\displaystyle\frac{-gh(\rho_s-\rho_w)}{2\rho} & -gh  \\
		0 & 0 & 0 & 0  \\
		0 & 0 & 0 & 0
	\end{pmatrix},\quad
	S_F = \begin{pmatrix}0\\-\epsilon |u_m| u_m\\ 0 \\ 0 	\end{pmatrix}, \quad \text{and} \quad S_{ED} =     \displaystyle \frac{E-D}{1-\psi} 
	\begin{pmatrix}
		1 \\
		u_m \\ 
		1-\psi \\ 
		-1 
	\end{pmatrix},
\end{equation*}
leading to the transport matrix
\begin{equation}\label{zeroth_order_moment}
	A_{ED} = \frac{\partial F}{\partial W}-B+ A_s= \begin{pmatrix}
		0 & 1 & 0 & 0   \\
		gh- u^{2}_m-\displaystyle	\frac{gh(\rho-\rho_w)}{2\rho} & 2u_m & 	\displaystyle\frac{gh(\rho_s-\rho_w)}{2\rho} & gh  \\
		-c_m u_m & c_m & u_m & 0  \\
		\delta_h & \delta_q &\delta_c & 0
	\end{pmatrix}
\end{equation}
where the sediment discharge matrix, $A_s=\mathbf{0}\in \mathbb{R}^{ 4\times 4}$. We see from \eqrefc{A_s defination} that $A_s$ is resulting from $\displaystyle \sum_{j=1}^{N}\alpha_jG_{ij}\partial_x\left(\frac{Q_b}{1-\psi}\right)$, representing the projection associated with the vertical coupling term and only appears in higher-order averaged models, where such coupling effects are included, since our model does not assume $\partial_t h_b=0$.\\
The $N=0$ model, i.e. SWEMED0 allows the derivation of a dissipative energy balance as shown in the following theorem.
\begin{theorem}\label{SWEMED0 energy balance}
	For the SWEMED0 model, the total mechanical energy density of the mixture is defined as 	\begin{equation}
		\mathcal{E}
		= \frac{1}{2} h u_m^2
		+ \frac{g}{2}(h + h_b)^2,
	\end{equation}
	 satisfies the following dissipative balance law
	\begin{equation}
		\begin{aligned}
			\partial_t \mathcal{E}
			&+ \partial_x\Biggl(
			u_m\biggl(\frac{h u_m^2}{2} + g h^2 + \frac{g\beta}{2}h^2 c_m\biggr)
			+ g(h+h_b)\frac{Q_b}{1-\psi}
			+ g h u_m h_b
			\Biggr)
			\\
			&= -\,\left(
			\epsilon |u_m|u_m^2
			- \tfrac{1}{2}S_m u_m^2
			- g\frac{Q_b}{1-\psi}\,\partial_x(h+h_b)
			\right)
			+ \frac{g\beta}{2}\,c_m\,\partial_x(h^2 u_m).
		\end{aligned}
		\label{eq:total_energy_final_1}
	\end{equation}
   Proof: A straightforward but somewhat lengthy proof can be found in \appref{app:4}.
\end{theorem}

\paragraph{\, SWEMED1}\label{par_first_order}
With $N=1$, the first-order moment represents the linear change of velocity in the vertical direction, i.e. $u(t,x,\zeta)=u_m(t,x)+\alpha_1(t,x)\phi_1(\zeta)$, where $\phi_1(\zeta)$ is the linear Legendre basis function defined in \eqrefc{Legendere polynomial}. The SWEMED1 model with state vector $W=(h,hu_m,h\alpha_1,hc_m,h_b)^T$ can be written in the form of \eqrefc{final SWMEED} using the following definitions
\begin{equation*}
	F= \begin{pmatrix}hu_m\\hu_m^2+\frac{1}{2}gh^2+\frac{1}{3} h \alpha^2_1\\2hu_m\alpha_1\\h c_m u_m\\ \frac{Q_b}{1-\psi}\end{pmatrix}, \quad	\frac{\partial F}{\partial W} = \begin{pmatrix}
		0 & 1 & 0  & 0 & 0 \\
		gh- u^{2}_m-\frac{\alpha^{2}_1}{3} & 2u_m & \frac{2\alpha_1}{3}  & 0 & 0\\
		-2\alpha_1u_m & 2\alpha_1 & 2u_m & 0 & 0\\
		-c_m u_m & c_m  & 0  & u_m  & 0\\
		\delta_h	  & \delta_q  & \delta_q  & \delta_c & 0
	\end{pmatrix},
\end{equation*} 
\begin{equation*}
	B= \begin{pmatrix}
		0 & 0 & 0  & 0 & 0\\
		\displaystyle \frac{gh(\rho-\rho_w)}{2\rho} & 0 & 0  & -\displaystyle \frac{gh(\rho_s-\rho_w)}{2\rho} & -gh \\
		\displaystyle \frac{gh(\rho-\rho_w)}{2\rho} & 0 & u_m & -\displaystyle \frac{gh(\rho_s-\rho_w)}{2\rho} & 0 \\
		0 & 0 & 0  & 0 & 0\\
		0 & 0 & 0  & 0 & 0\\
	\end{pmatrix},
\end{equation*} 
\begin{equation*}
	S_F = \begin{pmatrix}0\\-\epsilon |u_m+\alpha_1| (u_m+\alpha_1)\\ -3\epsilon |u_m+\alpha_1|(u_m+\alpha_1)- 12 \displaystyle \frac{\nu}{h}\alpha_1 \\ 0 \\ 0 \end{pmatrix}, \quad  S_{ED} = \displaystyle \frac{E-D}{1-\psi} \begin{pmatrix}1 \\ u_m+\alpha_1 \\  2\alpha_1 \\ 1-\psi\\ -1 \end{pmatrix},
\end{equation*}
\text{leading to the transport matrix}  
\begin{equation}\label{first_order_moment}
	A_{ED}= \displaystyle \frac{\partial F}{\partial W}-B+ A_s= \begin{pmatrix}
		0 & 1 & 0  & 0 & 0 \\
		gh- u^{2}_m-\displaystyle \frac{\alpha^{2}_1}{3}- \displaystyle \frac{gh(\rho-\rho_w)}{2\rho} & 2u_m & \displaystyle \frac{2\alpha_1}{3}  & \displaystyle \frac{gh(\rho_s-\rho_w)}{2\rho} & gh\\
		-2\alpha_1u_m-\displaystyle \frac{gh(\rho-\rho_w)}{2\rho} & 2\alpha_1 & u_m &\displaystyle \frac{gh(\rho_s-\rho_w)}{2\rho} & 0\\
		-c_m u_m & c_m  & 0  & u_m  & 0\\
		\delta_h& \delta_q  & \delta_q  & \delta_c & 0
	\end{pmatrix},
\end{equation} 
where the sediment discharge matrix,  $A_s=\mathbf{0}\in \mathbb{R}^{ 5\times 5}$ as $A_s=\displaystyle \sum_{j=1}^{N}\alpha_jG_{ij}\partial_x\left(\frac{Q_b}{1-\psi}\right)$ involves $G_{ij}$ defined in \eqrefc{constant matrices} and clearly $G_{11}= \displaystyle \int_{0}^{1}\phi_1(\zeta)\, \phi_1^\prime(\zeta)\,d\zeta=0.$\\
As in the previous section, the SWEMED1 model allows the derivation of a dissipative energy balance, as stated in the following theorem
\begin{theorem}\label{SWEMED1 energy balance}
	For the SWEMED1 model, the mechanical energy reads
	\[
	\mathcal{E}_1
	= \frac{1}{2}h u_m^2 + \frac{1}{6} h \alpha_1^2
	+ \frac{g}{2}(h+h_b)^2,
	\]
	and the corresponding balance law can be written in the following form
	\begin{equation}
		\begin{aligned}
			\partial_t \mathcal{E}_1
			&+ \partial_x\Biggl(
			\frac12 h u_m^3
			+ \frac12 h u_m\alpha_1^2
			+ \frac{g\beta}{2}h^2 c_m\Bigl(u_m+\frac{\alpha_1}{3}\Bigr)
			+ g h u_m (h+h_b)
			+ g (h+h_b)\frac{Q_b}{1-\psi}
			\Biggr)
			\\
			&= \frac{g\beta}{2}c_m\,\partial_x\!\left(h^2\Bigl(u_m+\frac{\alpha_1}{3}\Bigr)\right)
			- \left(
			\epsilon |u_b|u_b^2
			- \frac12 S_m u_b^2
			- g\frac{Q_b}{1-\psi}\,\partial_x (h+h_b)
			\right)
			-\frac{4\nu}{h}\alpha_1^2,
		\end{aligned}
		\label{eq:SWEMED1_E1_final_1}
	\end{equation}
where $u_b = u_m + \alpha_1$ is the bottom velocity.\\
Proof: The proof can be found in \appref{app:4}.
\end{theorem}

\begin{remark}
	We derive a dissipative mechanical-energy balance for the lower-order models SWEMED0 and SWEMED1. 
	For higher-order moments, the un-regularized moment system is not necessarily hyperbolic, so an energy analysis of the full hierarchy is less informative and is not shown here.  
	We note that related energy-based analyses for the shallow water linearized moment model is possible \cite{koellermeier2022steady}, where hyperbolicity can be ensured for any order of $N$.	
\end{remark}
\paragraph{\, SWEMED2} 
The second-order moment $N=2$ corresponds to a quadratic velocity profile in the vertical direction,
\begin{equation}\label{second_order_profile}
	u(t,x,\zeta)=u_m(t,x)+\alpha_1(t,x)\phi_1(\zeta)+\alpha_2(t,x)\phi_2(\zeta),
\end{equation}
Here, $\phi_1(\zeta)$ and $\phi_2(\zeta)$ are the linear and quadratic Legendre basis functions defined in \eqrefc{Legendere polynomial}. The SWEMED2 model with state vector $	W=(h,hu_m,h\alpha_1,h\alpha_2,hc_m,h_b)^T$ can be written in the form of \eqrefc{final SWMEED} using the following definitions
\begin{equation*}
	F= \begin{pmatrix}hu_m\\hu_m^2+\frac{1}{2}gh^2+\frac{1}{3} h \alpha^2_1+\frac{1}{5} h \alpha^2_2\\2hu_m\alpha_1+\frac{4}{5} h \alpha_1  \alpha_2\\2hu_m\alpha_2+\frac{2}{3} h \alpha^2_1+\frac{2}{7} h \alpha^2_2\\h c_m u_m\\ \frac{Q_b}{1-\psi}\end{pmatrix}, \quad
	\frac{\partial F}{\partial W}=\begin{pmatrix}
		0 & 1 & 0 & 0 & 0 & 0 \\
		gh- u^{2}_m-\frac{\alpha^{2}_1}{3}- \frac{\alpha^{2}_2}{5} & 2u_m & \frac{2\alpha_1}{3} &\frac{2\alpha_2}{5} & 0 & 0\\
		-2\alpha_1u_m- \frac{4}{5}\alpha_1\alpha_2 & 2\alpha_1 & 2u_m + \frac{4\alpha_2}{5} & \frac{4\alpha_1}{5} & 0 & 0\\
		-\frac{2}{3}\alpha^{2}_1 -2u_m\alpha_2-\frac{2}{7}\alpha^{2}_2  & 2\alpha_2 & \frac{4\alpha_1}{3} & 2u_m + \frac{4\alpha_2}{7} & 0 & 0\\
		-c_m u_m & c_m  & 0  & 0  & u_m  & 0\\
		\delta_h	  & \delta_q & \delta_q  & \delta_q & \delta_c & 0
	\end{pmatrix},
\end{equation*}
\begin{equation*}
	B= \begin{pmatrix}
		0 & 0 & 0 & 0 & 0 & 0\\
		\displaystyle \frac{gh(\rho-\rho_w)}{2\rho} & 0 & 0 & 0 & -\displaystyle \frac{gh(\rho_s-\rho_w)}{2\rho} & -gh \\
		\displaystyle \frac{gh(\rho-\rho_w)}{2\rho} & 0 & u_m-\displaystyle \frac{\alpha_2}{5} & \displaystyle \frac{\alpha_1}{5} & -\displaystyle\frac{gh(\rho_s-\rho_w)}{2\rho} & 0 \\
		0 & 0 & \alpha_1 &  u_m+\displaystyle \frac{\alpha_2}{7} & 0 & 0\\
		0 & 0 & 0 & 0 & 0 & 0\\
		0 & 0 & 0 & 0 & 0 & 0\\
	\end{pmatrix}, 
\end{equation*}
\begin{equation*}
	A_s=  \begin{pmatrix}
		\mathbf{0}_{1 \times 6} \\
		\mathbf{r} \\
		\mathbf{0}_{4 \times 6}
	\end{pmatrix}, \quad
	\text{where} \quad 
	\mathbf{r} = -6\alpha_2 \left(\delta_h,\ \delta_q,\ \delta_q,\ \delta_q,\ \delta_c,\ 0\right),
\end{equation*} 
\begin{equation*}
	S_F = \begin{pmatrix}0\\-\epsilon |u_m+\alpha_1+\alpha_2| (u_m+\alpha_1+\alpha_2)\\ -3\epsilon |u_m+\alpha_1+\alpha_2| (u_m+\alpha_1+\alpha_2)- 12 \displaystyle \frac{\nu}{h}\alpha_1 \\
		-5\epsilon |u_m+\alpha_1+\alpha_2|(u_m+\alpha_1+\alpha_2)\,- 60 \,\, \displaystyle \frac{\nu}{h}\alpha_2 \\ 0 \\ 0 \end{pmatrix}, \quad  S_{ED} = \displaystyle \frac{E-D}{1-\psi}\begin{pmatrix} 1 \\  u_m+\alpha_1+\alpha_2 \\ 2\alpha_1+3\alpha_2 \\ 3\alpha_2
		\\  1-\psi \\ -1 \end{pmatrix},
\end{equation*} 
leading to the transport matrix  
\begin{equation}
	A_{ED}= \displaystyle \frac{\partial F}{\partial W}-B+ A_s= \begin{pmatrix}
		0 & 1 & 0 & 0 & 0 & 0 \\
		gh- u^{2}_m-\frac{\alpha^{2}_1}{3}- \frac{\alpha^{2}_2}{5}- \frac{gh(\rho-\rho_w)}{2\rho} & 2u_m & \frac{2\alpha_1}{3} &\frac{2\alpha_2}{5} & \frac{gh(\rho_s-\rho_w)}{2\rho} & gh\\
		-2\alpha_1u_m- \frac{4}{5}\alpha_1\alpha_2- \frac{gh(\rho-\rho_w)}{2\rho} & 2\alpha_1 & u_m + \alpha_2 & \frac{3\alpha_1}{5} &\frac{gh(\rho_s-\rho_w)}{2\rho} & 0\\
		-\frac{2}{3}\alpha^{2}_1 -2u_m\alpha_2-\frac{2}{7}\alpha^{2}_2  & 2\alpha_2 & \frac{\alpha_1}{3} & u_m + \frac{3\alpha_2}{7} & 0 & 0\\
		-c_m u_m & c_m  & 0  & 0  & u_m  & 0\\
		\delta_h	  & \delta_q & \delta_q  & \delta_q & \delta_c & 0
	\end{pmatrix}.
\end{equation} 
\paragraph{\, SWEMED3} 
The third-order moment $N=3$ corresponds to a cubic velocity profile in the vertical direction,
\begin{equation}\label{third_order_profile}
	u(t,x,\zeta)=u_m(t,x)+\alpha_1(t,x)\phi_1(\zeta)+\alpha_2(t,x)\phi_2(\zeta)+\alpha_3(t,x)\phi_3(\zeta).
\end{equation}
Here, $\phi_1(\zeta)$, $\phi_2(\zeta)$ and $\phi_3(\zeta)$ are the linear, quadratic, and cubic Legendre basis functions defined in \eqrefc{Legendere polynomial}. The SWEMED3 model with the state vector $W=(h,hu_m,h\alpha_1,h\alpha_2,h\alpha_3,hc_m,h_b)^T$ can be written in the form of \eqrefc{final SWMEED}. The corresponding transport matrix is
\begin{equation}
	A_{ED}= \begin{pmatrix}
		0 & 1 & 0 & 0 & 0 & 0 & 0\\
		gh- u^{2}_m-\frac{\alpha^{2}_1}{3}- \frac{\alpha^{2}_2}{5}- \frac{\alpha^{2}_3}{7}- \frac{gh(\rho-\rho_w)}{2\rho} & 2u_m & \frac{2\alpha_1}{3} &\frac{2\alpha_2}{5} &\frac{2\alpha_3}{7}& \frac{gh(\rho_s-\rho_w)}{2\rho} & gh\\
		-\frac{2}{37}(9\alpha_3\alpha_2+7\alpha_1(5u_m+2\alpha_2))- \frac{gh(\rho-\rho_w)}{2\rho} & 2\alpha_1 & u_m + \alpha_2 & \frac{3(\alpha_1+\alpha_3)}{5}&\frac{3\alpha_2}{7} &\frac{gh(\rho_s-\rho_w)}{2\rho} & 0\\
		-\frac{2}{21}(7\alpha^2_1+9\alpha_1\alpha_3+2\alpha^2_3+3\alpha_2(7u_m+\alpha_2))  & 2\alpha_2 & \frac{\alpha_1}{3}+\frac{9\alpha_3}{7} & u_m + \frac{3\alpha_2}{7} &\frac{4\alpha_1}{7}+\frac{\alpha_3}{3} & 0 & 0\\
		-\frac{2}{15}(15u_m\alpha_3+4\alpha_2\alpha_3+9\alpha_1\alpha_2)  & 2\alpha_3 & 0 & \frac{2(\alpha_1+\alpha_3)}{5}  & u_m + \frac{\alpha_2}{3} & 0 & 0\\
		-c_m u_m & c_m  & 0  & 0 & 0 & u_m  & 0\\
		\delta_h	  & \delta_q & \delta_q  & \delta_q & \delta_q & \delta_c & 0
	\end{pmatrix},
\end{equation} 
and the right-hand side source terms $S = S_F+S_{ED}$ read 
\begin{equation*}
	S_F = \begin{pmatrix}0\\-\epsilon |u_m+\alpha_1+\alpha_2+\alpha_3| (u_m+\alpha_1+\alpha_2+\alpha_3)\\ -3\epsilon |u_m+\alpha_1+\alpha_2+\alpha_3| (u_m+\alpha_1+\alpha_2+\alpha_3)- 12 \displaystyle \frac{\nu}{h}(\alpha_1+\alpha_3) \\
		-5\epsilon |u_m+\alpha_1+\alpha_2+\alpha_3|(u_m+\alpha_1+\alpha_2+\alpha_3)\,- 60 \,\, \displaystyle \frac{\nu}{h}\alpha_2\\
		-7\epsilon |u_m+\alpha_1+\alpha_2+\alpha_3| (u_m+\alpha_1+\alpha_2+\alpha_3)- 7 \displaystyle \frac{\nu}{h}(4\alpha_1+ 24 \alpha_3)  \\ 0 \\ 0 \end{pmatrix}, \quad  S_{ED} = \displaystyle \frac{E-D}{1-\psi}\begin{pmatrix} 1 \\  u_m+\alpha_1+\alpha_2+ \alpha_3 \\ 2\alpha_1+3(\alpha_2+\alpha_3) \\ 3\alpha_2 + 5\alpha_3\\ 4 \alpha_3
		\\  1-\psi \\ -1 \end{pmatrix}.
\end{equation*} 

Next, we discuss the hyperbolic regularization of SWEMED and the hyperbolicity conditions for the resulting transport matrix.

%% file: 03-Hyperbolicity_modified.tex
\section{Hyperbolicity analysis of SWEMED}\label{hyperbolicity}
    In this section, we analyze the hyperbolicity of the transport matrix
	$A_{ED}(W)\in \mathbb{R}^{(N+4)\times(N+4)}$ of the SWEMED system
	\eqrefc{final SWMEED}, where $N$ denotes the order of the moment expansion. The matrix $A_{ED}(W)$ contains the shallow water
	moment block together with the sediment concentration and Exner coupling.
	We first recall the definition of hyperbolicity. Then, since the classical
	SWM is not hyperbolic in general for $N>1$, we regularize the moment block of $A_{ED}(W)$ following \cite{koellermeier2020analysis}. This gives the regularized transport matrix $\tilde{A}_{ED}(W)$ in
	\eqrefc{Aed transport matrix}. Next, we derive the characteristic polynomial \eqrefc{CP_HSWEMED} and its factorization separates the sediment-coupled cubic factor $P_{\mathrm{sed}}(\lambda)$ \eqrefc{eq:P_sed}, the factor $(u_m - \lambda)$, and the regularized moment factor $\chi_{A_2}(\lambda-u_m)$. 
	The remaining questions are to determine when $P_{\mathrm{sed}}(\lambda)$ has three distinct real roots and to clarify how the moment factor affects strict and weak hyperbolicity for even and odd values of $N$.

	To define the hyperbolicity of a first-order partial differential equation of the form
	\begin{equation}\label{PDE}
		\partial_t W + A(W)\partial_x W = 0.
	\end{equation}
	we use the following terminology \cite{evans2022partial,koellermeier2020analysis,KreissOrtiz2002,SmithPettyPantano2016}.
	
	\begin{definition}\label{hyp_def}
		The system \eqrefc{PDE} is called hyperbolic if  $A(W)$ has only real eigenvalues and is diagonalizable. It is called strictly hyperbolic if all eigenvalues of $A(W)$ are real and distinct. It is called weakly hyperbolic if $A(W)$ has only real eigenvalues, but repeated eigenvalues occur and diagonalizability is not guaranteed.
	\end{definition}

The characteristic polynomial of the first-order $(N=1)$ Shallow Water Moment model (SWM1) without erosion, deposition, and bedload transport is
\begin{equation}
	\label{eq:cp_SWM1_final}
	\chi_A(\lambda)
	=
	(u_m-\lambda)
	\left((\lambda-u_m)^2-gh-\alpha_1^2\right),
\end{equation}
with the three real and distinct eigenvalues \parencite{CiCP-25-669}
\begin{equation}
	\label{eq:SWM1_eigs_final}
	\lambda_{1,2}=u_m\pm\sqrt{gh+\alpha_1^2},
	\qquad
	\lambda_3=u_m.
\end{equation}
Therefore, for $h>0$, the classical SWM1 is strictly hyperbolic. However,
for higher-order $(N>1)$ SWM, the transport matrix may have complex eigenvalues and is therefore not globally hyperbolic
\parencite{koellermeier2020analysis, CiCP-25-669}. To address this issue, \cite{koellermeier2020analysis} introduced a regularization technique based on a linear deviation from equilibrium, thereby ensuring hyperbolicity of the newly derived hyperbolic shallow water moment system. For our proposed SWEMED model,
we apply the same regularization strategy to the moment-transport part of $A_{ED}(W)$ \eqrefc{final SWMEED} for arbitrary order $N$ by modifying the sub-matrix associated with the fluid transport, not the sediment concentration and bedload transport. This is done by keeping the terms in the first $N+2$ rows and $N+2$ columns of the matrix $A_{ED}(W)\in \mathbb{R}^{(N+4)\times(N+4)}$ containing $\alpha_1$ and forcing all higher-order moments to be zero in the regularized transport block, i.e.,
$
\alpha_2,\alpha_3,\ldots,\alpha_N=0.
$
This regularization guarantees real eigenvalues for the moment block inherited from the shallow water moment system.
However, it does not by itself prove hyperbolicity of the full sediment-coupled matrix; the additional sediment-coupled factor has to be analysed separately.
We denote the resulting regularized transport matrix by $\tilde{A}_{ED}(W)$, which reads
\begin{equation}\label{Aed transport matrix}
	\tilde{A}_{ED}(W)=\begin{pmatrix}
		& 1 &  &  &  &  & &  \\
		gh-u_m^2-\frac{\alpha_1^2}{3}-\frac{gh(\rho-\rho_w)}{2\rho} & 2u_m & \frac{2\alpha_1}{3} & &  &  & \frac{gh(\rho_s-\rho_w)}{2\rho} & gh \\
		-2u_m\alpha_1-\frac{gh(\rho-\rho_w)}{2\rho}  & 2\alpha_1 & u_m & \frac{3\alpha_1}{5} &  &  & \frac{gh(\rho_s-\rho_w)}{2\rho} & \\
		-\frac{2\alpha_1^2}{3}&  & \frac{\alpha_1}{3} & u_m & \ddots &  &  & \\
		&  &  & \ddots & \ddots  & \frac{N+1}{2N+1}\alpha_1 &  &  \\
		&  &  &  & \frac{N-1}{2N-1}\alpha_1 & u_m & &  \\
		- c_m u_m & c_m  &   &  & &  & u_m &  \\
		\delta_h	  & \delta_q  & \delta_q  &  & \dots  & \delta_q & \delta_c &
	\end{pmatrix},
\end{equation}
	where all other entries are zero. 
	The regularization removes the explicit dependence on the higher-order coefficients from the moment block of $\tilde{A}_{ED}(W)$, except through the bedload derivatives $\delta_h,\delta_q,$ and $\delta_c$. However, the higher-order coefficient equations remain coupled with the other equations, so the complete system remains non-linear. Throughout this work, we retain the name Hyperbolic Shallow Water Exner Moment model with Erosion and Deposition (HSWEMED) for the regularized model, since it is obtained through the hyperbolic regularization framework used for HSWM and HSWEM. Nevertheless, depending on the moment order and the eigenvalue separation, the resulting transport matrix may be strictly hyperbolic or only weakly hyperbolic. We rewrite the complete system \eqrefc{final SWMEED} with the regularized transport matrix as
\begin{equation}\label{final HSWEMED}
	\partial_t W + \tilde{A}_{ED}(W) \partial_x W= S(W).
\end{equation}
Note that the source terms in $S(W)$ in \eqrefc{final HSWEMED}, including erosion, deposition, and friction, do not enter the characteristic polynomial of the transport matrix. The following theorem gives the characteristic polynomial of $\tilde{A}_{ED}(W)$ in \eqrefc{Aed transport matrix}. It extends the characteristic-polynomial analysis of HSWEM in \cite{Garres}, which is based on the regularized HSWM structure in \cite{koellermeier2020analysis}, to the additional sediment concentration equation.

\begin{theorem}\label{theorem}
	The transport matrix $\tilde{A}_{ED}(W) \in \mathbb{R}^{(N+4)\times(N+4)}$ \eqrefc{Aed transport matrix} of HSWEMED has the following characteristic polynomial
	\begin{equation}\label{CP_HSWEMED}
		\chi_{\tilde{A}_{ED}}(\lambda)
		=(u_m-\lambda)P_{\mathrm{sed}}(\lambda)\chi_{A_2}(\lambda-u_m),
	\end{equation}
	where
	\begin{equation}\label{eq:P_sed}
		P_{\mathrm{sed}}(\lambda)
		=
		(-\lambda)\left((\lambda-u_m)^2-gh-\alpha_1^2\right)
		+gh\left(\delta_h+\lambda\delta_q+2\alpha_1\delta_q+c_m\delta_c\right),
	\end{equation}
	and
	\begin{equation}
		\chi_{A_{2}}(\lambda-u_m)
		=
		\det \left(A_{2}- \left(\lambda-u_m \right)I \right),
	\end{equation}
	where $I \in \mathbb{R}^{N\times N}$ is the identity matrix and $A_{2} \in \mathbb{R}^{N \times N}$ is defined by
	\begin{equation}\label{eq:A2_matrix}
			A_{2} = \begin{pmatrix}
				& c_2 &  &    \\
				a_2 &  & \ddots &     \\
				& \ddots &   & c_N     \\
				&    & a_N &
			\end{pmatrix},
			\qquad
			c_i = \frac{i+1}{2i+1}\alpha_1,
			\quad
			a_i = \frac{i-1}{2i-1}\alpha_1.
	\end{equation}
\end{theorem}

\begin{proof}
	The proof follows \cite{Garres} and can be found in \appref{app:3}.
\end{proof}
We now identify under which conditions all factors in \eqrefc{CP_HSWEMED} provide real eigenvalues for the full HSWEMED transport matrix $\tilde{A}_{ED}(W)$ \eqrefc{Aed transport matrix}. We denote the three roots of $P_{\mathrm{sed}}(\lambda)$ in \eqrefc{eq:P_sed} by $\lambda_1,\lambda_2$, and $\lambda_3$. These roots give the first three eigenvalues of $\tilde{A}_{ED}(W)$. The factor $(u_m-\lambda)$ gives the fourth eigenvalue, $\lambda_4=u_m$, which is associated with the concentration equation. The remaining $N$ eigenvalues, $\lambda_5,\ldots,\lambda_{N+4}$, are obtained from the regularized moment block $\chi_{A_2}(\lambda-u_m)=0$. The eigenvalues associated with the regularized moment block are real \cite{Garres,koellermeier2020analysis}. Therefore, it remains to show that the cubic polynomial $P_{\mathrm{sed}}(\lambda)$ in \eqrefc{eq:P_sed} has three real roots.
A comparable cubic real-root analysis has been carried out for shallow water Exner systems with bedload transport \cite{diaz2008sediment,cordier2011bedload}, for coupled suspended-load and bedload transport \cite{gonzalez2020robust}, and for HSWEM \cite{Garres}. The cubic polynomial $P_{\mathrm{sed}}(\lambda)$ in \eqrefc{eq:P_sed} extends the corresponding cubic factor of HSWEM \cite{Garres}. The main difference is the additional term $c_m\delta_c$ in \eqrefc{eq:P_sed}, which originates from the dependence of the bedload flux on the suspended concentration variable $hc_m$. We follow the same cubic-line real-root analysis for the proposed HSWEMED model.
For the Meyer--Peter and M\"uller bedload closure \eqrefc{eq:MPM_bed_closure} used in this work, the bedload derivatives \eqrefc{eq:delta_q}--\eqrefc{eq:delta_alpha} satisfy
\begin{equation}
	\label{eq:delta_identity_hswemed}
	\delta_h+c_m\delta_c=-u_b\delta_q,
\end{equation}
where $u_b$ in \eqrefc{bottom_velocity} is the bottom velocity associated with the moment expansion.

Using \eqrefc{eq:delta_identity_hswemed}, the cubic polynomial \eqrefc{eq:P_sed} reduces to
	\begin{equation}\label{eq:P_sed_reduced}
		P_{\mathrm{sed}}(\lambda)
		=
		(-\lambda)\left((\lambda-u_m)^2-gh-\alpha_1^2\right)
		+gh\delta_q(\lambda-u_b+2\alpha_1).
	\end{equation}
	Thus, the explicit concentration term $c_m\delta_c$ cancels with the corresponding concentration-dependent part of $\delta_h$ \eqrefc{eq:delta_h}. The concentration may still influence the coefficients of \eqrefc{eq:P_sed_reduced} through $\delta_q$, since $\delta_q$ depends on the mixture density $\rho(c_m)$.

	Following the cubic-line construction in \cite{cordier2011bedload,Garres,gonzalez2020robust}, we define
	\begin{equation}\label{eq:f_d_alpha_def}
	f_{\alpha_1}(\lambda)
	:=
	\lambda\left((\lambda-u_m)^2-gh-\alpha_1^2\right),
	\qquad
	d_{\alpha_1}(\lambda)
	:=
	gh\delta_q(\lambda-u_b+2\alpha_1).
\end{equation}
	Then $P_{\mathrm{sed}}(\lambda)=0$ in \eqrefc{eq:P_sed_reduced} is equivalent to $f_{\alpha_1}(\lambda)=d_{\alpha_1}(\lambda)$. Hence, the three roots of $P_{\mathrm{sed}}(\lambda)$ are the intersection points of the cubic $f_{\alpha_1}(\lambda)$ and the line $d_{\alpha_1}(\lambda)$.
	
\begin{figure}[ht]
	\centering
	\includegraphics[width=0.86\textwidth]{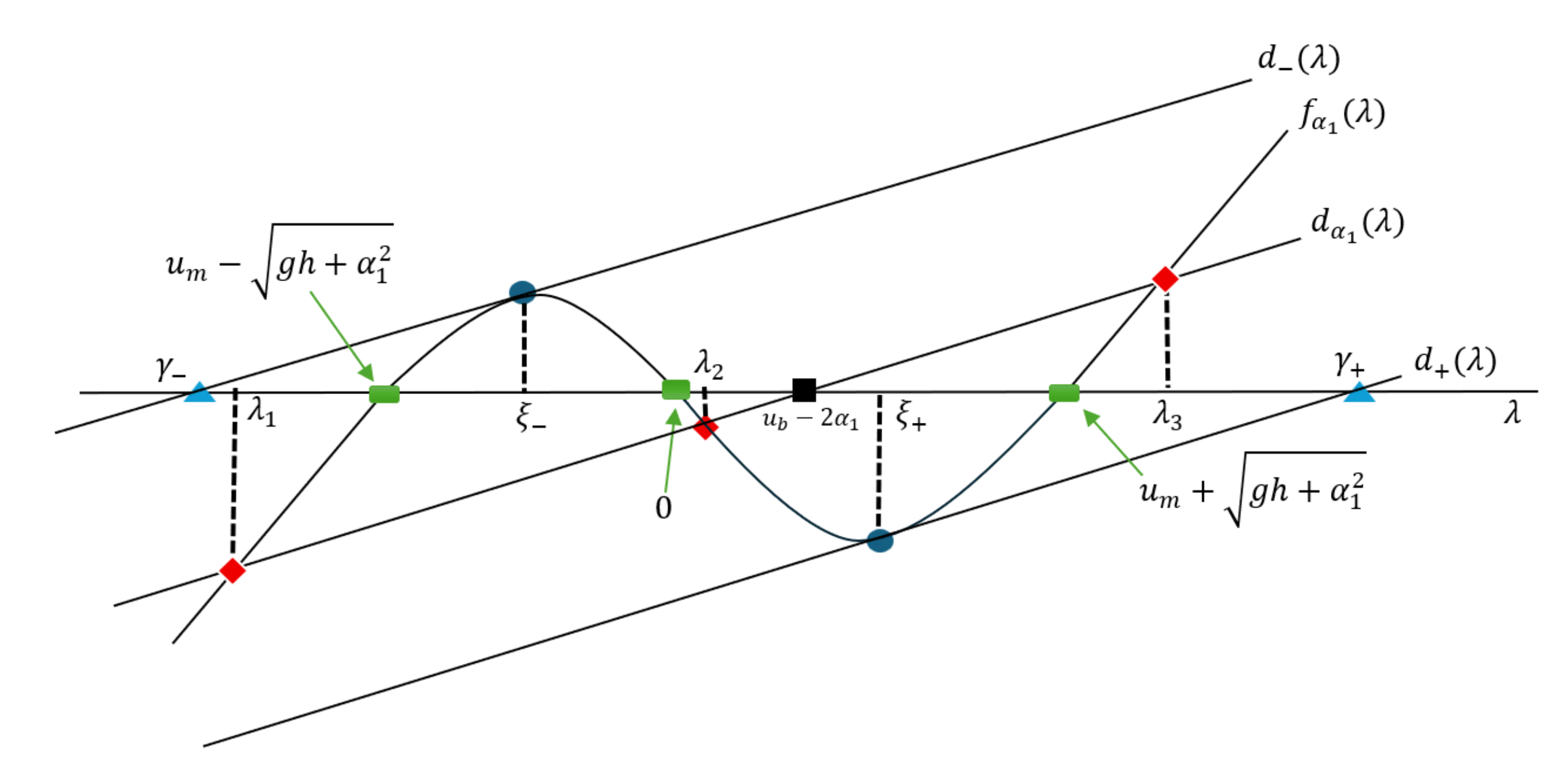}
	\caption{Cubic-line interpretation of the sediment-coupled polynomial $P_{\mathrm{sed}}(\lambda)$ in \eqrefc{eq:P_sed_reduced} \cite{cordier2011bedload,gonzalez2020robust}. The three roots of $P_{\mathrm{sed}}(\lambda)$ are $\lambda_1$, $\lambda_2$, and $\lambda_3$, which give the first three eigenvalues of $\tilde{A}_{ED}(W)$ \eqrefc{Aed transport matrix}, are the intersection points of $f_{\alpha_1}(\lambda)$ and $d_{\alpha_1}(\lambda)$. The lines $d_-(\lambda)$ and $d_+(\lambda)$ are tangent to $f_{\alpha_1}(\lambda)$ at $\xi_-$ and $\xi_+$ and have the same slope as $d_{\alpha_1}(\lambda)$. Their zeros $\gamma_-$ and $\gamma_+$ define the admissible interval for the zero of $d_{\alpha_1}(\lambda)$, namely $\lambda=u_b-2\alpha_1$. In the configuration shown, condition \eqrefc{eq:main_tangent_condition_hswemed} is satisfied and $P_{\mathrm{sed}}(\lambda)$ has three distinct real roots. The red diamonds indicate these roots, the blue circles indicate the tangency points, the green squares indicate the zeros of $f_{\alpha_1}(\lambda)$, and the blue triangles indicate $\gamma_-$ and $\gamma_+$.}
	\label{fig:hyperbolicity_figures_final}
\end{figure}

	The line $d_{\alpha_1}(\lambda)$ has slope $gh\delta_q$ and vanishes at
	\begin{equation}\label{eq:d_zero}
		\lambda=u_b-2\alpha_1
		=u_m-\alpha_1+\alpha_2+\cdots+\alpha_N.
	\end{equation}
	For a fixed local state $W\in \mathbb{R}^{N+4}$, the value $u_b-2\alpha_1$ gives the point where $d_{\alpha_1}(\lambda)$ crosses the $\lambda$-axis. In the geometric construction shown in Figure~\ref{fig:hyperbolicity_figures_final}, $h$, $u_m$, $\alpha_1$, and $\delta_q$ are fixed, while $u_b-2\alpha_1$ is varied. The higher-order moments enter this quantity through $u_b$; in the full model, they may also influence the bedload derivatives.
	The boundary of the three-real-root region is reached when two roots merge. In the cubic-line representation, this occurs when the line $d_{\alpha_1}(\lambda)$ is tangent to the cubic $f_{\alpha_1}(\lambda)$. At a tangency point $\xi$,
	\begin{equation}\label{eq:tangent_slope_hswemed}
		f'_{\alpha_1}(\xi)=d'_{\alpha_1}(\xi)=gh\delta_q,
	\end{equation}
	where
	\begin{equation}\label{eq:fprime_hswemed}
		f'_{\alpha_1}(\lambda)
		=3\lambda^2-4u_m\lambda+u_m^2-gh-\alpha_1^2,
		\qquad
		d'_{\alpha_1}(\lambda)=gh\delta_q.
	\end{equation}
	Solving \eqrefc{eq:tangent_slope_hswemed} gives
	\begin{equation}\label{eq:xi_pm_hswemed}
		\xi_\pm
		=
		\frac{2u_m\pm\sqrt{u_m^2+3\left(gh+\alpha_1^2+gh\delta_q\right)}}{3}.
	\end{equation}

Let $d_-(\lambda)$ and $d_+(\lambda)$ be the tangent lines to $f_{\alpha_1}$ \eqrefc{eq:f_d_alpha_def} at $\xi_-$ and $\xi_+$, respectively. They are given by
\begin{equation}
	\label{eq:d_pm_hswemed}
	d_\pm(\lambda)
	=
	f_{\alpha_1}(\xi_\pm)+gh\delta_q(\lambda-\xi_\pm)
	=
	gh\delta_q(\lambda-\gamma_\pm),
\end{equation}
where
\begin{equation}
	\label{eq:gamma_pm_hswemed}
	\gamma_\pm
	=
	\xi_\pm-\frac{f_{\alpha_1}(\xi_\pm)}{gh\delta_q}.
\end{equation}
The numbers $\gamma_-$ and $\gamma_+$ are the zeros of the tangent lines $d_-(\lambda)$ and $d_+(\lambda)$; they are not eigenvalues.

\begin{proposition}
	\label{prop:sediment_real_root_hswemed}
	Assume $h>0$ and $\delta_q>0$. The cubic polynomial $P_{\mathrm{sed}}(\lambda)$ \eqrefc{eq:P_sed_reduced} has three distinct real roots if and only if
	\begin{equation}
		\label{eq:main_tangent_condition_hswemed}
		\min(\gamma_-,\gamma_+)
		<
		u_b-2\alpha_1
		<
		\max(\gamma_-,\gamma_+).
	\end{equation}
\end{proposition}

\begin{proof}
	From \eqrefc{eq:P_sed_reduced} and \eqrefc{eq:f_d_alpha_def},
	\[
	P_{\mathrm{sed}}(\lambda)
	=
	d_{\alpha_1}(\lambda)-f_{\alpha_1}(\lambda).
	\]
	Its critical points are $\xi_-$ and $\xi_+$ because
	\[
	P'_{\mathrm{sed}}(\xi_\pm)
	=
	d'_{\alpha_1}(\xi_\pm)-f'_{\alpha_1}(\xi_\pm)
	=0.
	\]
	A cubic polynomial has three distinct real roots if and only if its local maximum and local minimum have opposite signs. Using \eqrefc{eq:d_pm_hswemed} and \eqrefc{eq:gamma_pm_hswemed},
	\begin{align}
		P_{\mathrm{sed}}(\xi_\pm)
		&=
		d_{\alpha_1}(\xi_\pm)-f_{\alpha_1}(\xi_\pm)\\
		&=
		gh\delta_q\left(\xi_\pm-(u_b-2\alpha_1)\right)
		-
		gh\delta_q(\xi_\pm-\gamma_\pm)\\
		&=
		gh\delta_q\left(\gamma_\pm-(u_b-2\alpha_1)\right).
	\end{align}
	Since $gh\delta_q>0$, these two values have opposite signs exactly when $u_b-2\alpha_1$ lies strictly between $\gamma_-$ and $\gamma_+$. This proves \eqrefc{eq:main_tangent_condition_hswemed}.
\end{proof}

\begin{remark}
	If $\delta_q=0$ and the identity \eqrefc{eq:delta_identity_hswemed} holds, then \eqrefc{eq:P_sed_reduced} reduces to
	\begin{equation*}
	P_{\mathrm{sed}}(\lambda)
	=
	-\lambda\left((\lambda-u_m)^2-gh-\alpha_1^2\right),
	\end{equation*}
	so the sediment-coupled factor $P_{\mathrm{sed}}(\lambda)$ has the real roots
	\begin{equation*}
	\lambda_{1,2}=u_m\pm\sqrt{gh+\alpha_1^2},\qquad \lambda_3=0.
	\end{equation*}
\end{remark}
\begin{remark}\label{rem:remaining_eigenvalues}
	In the characteristic polynomial \eqrefc{CP_HSWEMED} of the HSWEMED transport matrix $\tilde{A}_{ED}(W)\in \mathbb{R}^{(N+4)\times(N+4)}$, assume that $P_{\mathrm{sed}}(\lambda)$ has three real roots, denoted by $\lambda_1$, $\lambda_2$, and $\lambda_3$. The remaining $N+1$ eigenvalues of $\tilde{A}_{ED}(W)$ in \eqrefc{Aed transport matrix} are obtained from
	\begin{equation}\label{moment_CP}
		(u_m-\lambda)\chi_{A_2}(\lambda-u_m)=0.
	\end{equation}
	The factor $(u_m-\lambda)$ gives $\lambda_4=u_m$. The $N$ eigenvalues associated with the regularized moment block are
	\begin{equation}\label{eq:moment_block_eigenvalues}
		\lambda_{i+4}=u_m+b_i\alpha_1,
		\qquad i=1,\ldots,N,
	\end{equation}
	where $b_i$, $i=1,\ldots,N$, are the $N$ roots of the derivative of the Legendre polynomial $P_{N+1}$ \cite{Bauerle2024,koellermeier2026primitive}. Hence, for $\alpha_1\neq0$, the eigenvalues associated with the moment block are real and distinct.
	
	\tabref{tab:remaining_eigenvalues} lists the remaining eigenvalues obtained from \eqrefc{moment_CP} for $N=1,2,3,4$.  
	For $N=1$ and $N=3$, the eigenvalue $u_m$ appears repeatedly: it is obtained both from the factor $(u_m-\lambda)$ and from the regularized moment block. Therefore, the HSWEMED1 and HSWEMED3 transport matrices are not strictly hyperbolic. Under the real-root condition for $P_{\mathrm{sed}}(\lambda)$, these odd-order matrices have a real spectrum but contain the repeated eigenvalue $u_m$. In accordance with definition \ref{hyp_def}, we therefore refer to this case as weakly hyperbolic. For $N=2$ and $N=4$, the eigenvalues listed in \tabref{tab:remaining_eigenvalues} are distinct if $\alpha_1\neq0$. Strict hyperbolicity of the full HSWEMED2 and HSWEMED4 transport matrices additionally requires that the three roots of $P_{\mathrm{sed}}(\lambda)$ are distinct and do not merge with any of the eigenvalues listed in the table.
\end{remark}
\begin{table}[htbp]
	\centering
	\caption{Eigenvalues of $\tilde{A}_{ED}(W)$ obtained from the factor $(u_m-\lambda)\chi_{A_2}(\lambda-u_m)$ for $N=1,2,3,4$. The first three eigenvalues $\lambda_1$, $\lambda_2$, and $\lambda_3$, which are the roots of $P_{\mathrm{sed}}(\lambda)$, are not included in the table.}
	\label{tab:remaining_eigenvalues}
	\renewcommand{\arraystretch}{1.35}
	\setlength{\tabcolsep}{5pt}
	\resizebox{\textwidth}{!}{%
		\begin{tabular}{|c|c|c|c|c|c|}
			\hline
			$N$ & $\lambda_4$ & $\lambda_5$ & $\lambda_6$ & $\lambda_7$ & $\lambda_8$ \\
			\hline
			$1$ & $u_m$ & $u_m$ & -- & -- & -- \\
			\hline
			$2$ & $u_m$ & $u_m+\sqrt{\dfrac{1}{5}}\,\alpha_1$ & $u_m-\sqrt{\dfrac{1}{5}}\,\alpha_1$ & -- & -- \\
			\hline
			$3$ & $u_m$ & $u_m+\sqrt{\dfrac{3}{7}}\,\alpha_1$ & $u_m-\sqrt{\dfrac{3}{7}}\,\alpha_1$ & $u_m$ & -- \\
			\hline
			$4$ & $u_m$ & $u_m+\sqrt{\dfrac{1}{3}-\dfrac{2}{3\sqrt{7}}}\,\alpha_1$ & $u_m-\sqrt{\dfrac{1}{3}-\dfrac{2}{3\sqrt{7}}}\,\alpha_1$ & $u_m+\sqrt{\dfrac{1}{3}+\dfrac{2}{3\sqrt{7}}}\,\alpha_1$ & $u_m-\sqrt{\dfrac{1}{3}+\dfrac{2}{3\sqrt{7}}}\,\alpha_1$ \\
			\hline
		\end{tabular}%
	}
\end{table}

\begin{remark}
	For $N=0$, the HSWEMED model reduces to the
	Shallow Water Exner model with Erosion and Deposition (SWEED) \eqrefc{zeroth_order_moment}. In this case, hyperbolicity requires the three roots of $P_{\mathrm{sed}}(\lambda)$ to be real
	\cite{diaz2008sediment,cordier2011bedload,gonzalez2020robust}.
	Strict hyperbolicity additionally requires these roots to be distinct
	and different from the eigenvalue $u_m$. For $N\geq1$, setting $\alpha_1=0$ does not remove the moment equations. The sediment coupled cubic $P_{\mathrm{sed}}(\lambda)$ reduce to the
	corresponding SWEED cubic, while the eigenvalues obtained from
	$\chi_{A_2}(\lambda-u_m)$ reduces to $u_m$. Together with the
	eigenvalue $u_m$ obtained from the factor $(u_m-\lambda)$, this gives
	a repeated eigenvalue $u_m$ of algebraic multiplicity at least $N+1$.
	Therefore, the HSWEMED transport matrix is not strictly hyperbolic when $\alpha_1=0$. Under the corresponding real-root condition for the sediment-coupled cubic, the matrix has a real spectrum but repeated eigenvalues; hence this case is weakly hyperbolic.

\end{remark}

To summarize, the characteristic polynomial of the regularized HSWEMED transport matrix $\tilde{A}_{ED}(W)$ in \eqrefc{CP_HSWEMED} consists of the sediment-coupled cubic polynomial $P_{\mathrm{sed}}(\lambda)$, the factor $(u_m-\lambda)$, and the regularized moment factor $\chi_{A_2}(\lambda-u_m)$. Proposition~\ref{prop:sediment_real_root_hswemed} gives a condition under which $P_{\mathrm{sed}}(\lambda)$ has three distinct real roots. Under this condition, the regularized HSWEMED transport matrix has a real spectrum. For even values of the moment order $N$ and non-zero $\alpha_1$, strict hyperbolicity holds provided that the three roots of $P_{\mathrm{sed}}(\lambda)$ remain distinct from $u_m$ and from the moment-block eigenvalues. For odd values of $N$ and non-zero $\alpha_1$, the moment block contains the eigenvalue $u_m$, which is also obtained from the factor $(u_m-\lambda)$. Hence, $u_m$ is repeated and strict hyperbolicity is lost. In this case, the regularized HSWEMED transport matrix is weakly hyperbolic under the real-root condition above.

In the following section, the numerical results obtained with HSWEMED are compared with the laboratory dam-break experiments available in \cite{Spinewine01122007}.

%% file: 04-Numerical_simulation_modified.tex
\section{Numerical simulations}\label{sec4}
We numerically solve the HSWEMED system \eqrefc{final HSWEMED} using a second-order path-conservative Monotonic Upstream-Centered Scheme for Conservation Laws (MUSCL) \cite{vanLeer1979MUSCL}, combined with a local Rusanov scheme \cite{castro2012class,rusanov1962}. Since \eqrefc{final HSWEMED} contains non-conservative products, we evaluate the corresponding contributions along straight paths connecting the reconstructed states at each interface \cite{castroGallardoPares2006,pares2006}. We also retain the internal path contribution associated with the linear reconstruction inside each cell \cite{castroGallardoPares2006}. We use the generalized minmod limiter with parameter $\theta=1.3$ for the reconstruction \cite{kurganovLevy2002} and a third-order, four-stage Runge-Kutta method for the time integration.

\begin{figure}[ht]
	\centering
	\includegraphics[width=0.68\textwidth]{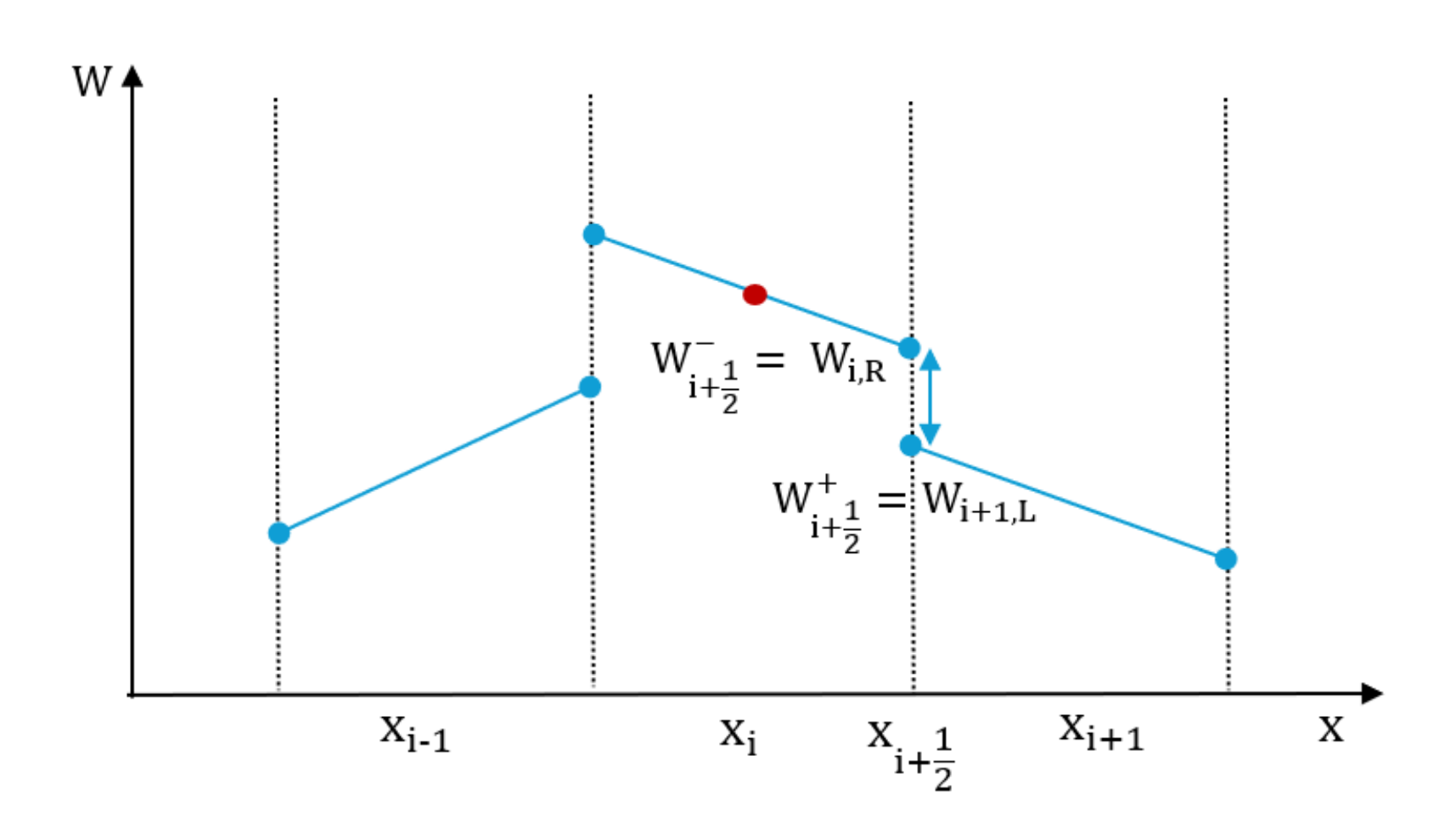}
	\caption{MUSCL reconstruction on three neighbouring cells. The reconstructed states $W^-_{i+1/2}$ and $W^+_{i+1/2}$ define the left and right states at the interface $x_{i+1/2}$.}
	\label{fig:MUSCL-reconstruction-standard}
\end{figure}

Figure~\ref{fig:MUSCL-reconstruction-standard} illustrates the MUSCL reconstruction on three neighbouring cells. At the interface $x_{i+1/2}$, the Rusanov artificial viscosity contribution is given by \cite{castro2012class,rusanov1962}. 
\begin{equation}\label{Rusanov artificial viscosity}
	a_{i+1/2}\Delta W_{i+1/2},
	\qquad
	\Delta W_{i+1/2}
	=
	W_{i+1/2}^{+}-W_{i+1/2}^{-}
\end{equation}
where
\[
a_{i+1/2}
=
\max\left\{
\max_k\left|\lambda_k\!\left(A(W_{i+1/2}^{-})\right)\right|,
\max_k\left|\lambda_k\!\left(A(W_{i+1/2}^{+})\right)\right|
\right\}.
\]
Here, $W_{i+1/2}^{-}$ and $W_{i+1/2}^{+}$ denote the reconstructed states at the interface.

The standard MUSCL-Rusanov scheme does not preserve the water-at-rest equilibrium over a non-flat bottom. 
In well-balanced finite-volume methods, stationary equilibria are commonly preserved by separating the equilibrium contribution from the perturbation that is advanced by the numerical scheme  \cite{audusse2004fast,kurganovLevy2002,koellermeier2022steady}. Motivated by this principle, we remove only the water-at-rest contribution from the Rusanov viscosity in the present HSWEMED discretization.

At water-at-rest equilibrium,
\[
u_m=\alpha_1=\cdots=\alpha_N=c_m=0,
\qquad
h+h_b=\mathrm{constant},
\]
and therefore
\[
\Delta h=-\Delta h_b.
\]
Although this jump corresponds to a stationary solution, the Rusanov viscosity \eqrefc{Rusanov artificial viscosity} acts on the complete jump of the state vector $W$ and perturbs the free surface $h+h_b$. To preserve the equilibrium, we remove the corresponding equilibrium contribution
\[
\Delta W_{i+1/2}^{\mathrm{eq}}
=
\left(
-\Delta h_{b,i+1/2},
0,\ldots,0,
\Delta h_{b,i+1/2}
\right)^{T},
\]
from the artificial viscosity term \eqrefc{Rusanov artificial viscosity}. The corrected contribution is
\begin{equation}\label{corrected Rusanov}
	a_{i+1/2}
	\left(
	\Delta W_{i+1/2}
	-
	\Delta W_{i+1/2}^{\mathrm{eq}}
	\right).
\end{equation}
The correction is applied at every interface. We modify only the Rusanov viscosity, while the reconstructed states, fluxes, path contributions, and source terms remain unchanged. We refer to the resulting scheme as the well-balanced MUSCL-Rusanov scheme for the water-at-rest equilibrium $(\text{MUSCL}_{\mathrm{WB}})$.

We consider four main test settings: first, we verify the well-balanced property for smooth and discontinuous bottom configurations in Section~\ref{water-at-rest-tests}; second, we verify the accuracy of the second-order MUSCL-Rusanov scheme in Section~\ref{accuracy_test}; third, we apply the scheme to the academic dam-break problem in Section~\ref{academic dam-break1}; and fourth, we compare the numerical results for different models with laboratory dam-break experiments in Section~\ref{laboratory experiments}.

\subsection{Well-balanced verification for the water-at-rest equilibrium}
\label{water-at-rest-tests}

In this section, we verify that the MUSCL-Rusanov scheme with the proposed water-at-rest correction \eqrefc{corrected Rusanov} is well-balanced for the water-at-rest equilibrium. We consider two bottom configurations: a smooth Gaussian bottom \eqrefc{war-smooth-ic} and a discontinuous rectangular bottom \eqrefc{war-rectangle-ic}.
For both configurations, we consider a constant initial free surface, $h(0,x)+h_b(0,x)=1.$

For the smooth Gaussian bottom $h_b$, the initial condition is
\begin{equation}
	\label{war-smooth-ic}
	h_b(0,x)
	=
	0.2\exp\!\left[
	-\left(\frac{x-5}{0.65}\right)^2
	\right],
	\qquad
	h(0,x)=1-h_b(0,x).
\end{equation}
For the rectangular bottom $h_b$, we consider
\begin{equation}
	\label{war-rectangle-ic}
	h_b(0,x)=
	\begin{cases}
		0.2, & 4\leq x\leq 6,\\
		0, & \text{otherwise},
	\end{cases}
	\qquad
	h(0,x)=1-h_b(0,x).
\end{equation}
In both cases, the water is fully settled and at rest, so that the remaining initial variables are
\begin{equation*}
	u_m(0,x)=0,
	\qquad
	\alpha_i(0,x)=0,\quad i=1,\ldots,N,
	\qquad
	c_m(0,x)=0.
\end{equation*}
The properties of the sediment particles considered in both tests are shown in \tabref{table:water-at-rest}. The computational domain is $\Omega=[0,10]$, divided into $N_x=400$ cells, and the solution is computed until $t=0.5$. We consider the third-order HSWEMED, i.e., $N=3$, for this test.
\begin{table}[H]
	\centering
	\begin{tabular}{|c|c|c|c|c|c|}
		\hline
		$\rho_w$ & $\rho_s$ &$\theta_c$ &$ \psi $& $d_s $ & $ \epsilon $ \\ \hline
		1000 & 1580 & 0.047 & 0.47 & 3.9 & 0.0324 \\ \hline
	\end{tabular}
	\caption{Sediment properties of water-at-rest equilibrium test.}
	\label{table:water-at-rest}
\end{table}
However, these parameters do not affect the stationary solutions considered in this test. Open boundary conditions are imposed at both ends of the computational domain.
\begin{figure}[H]
	\centering
	\begin{subfigure}[b]{0.49\textwidth}
		\includegraphics[width=\textwidth]{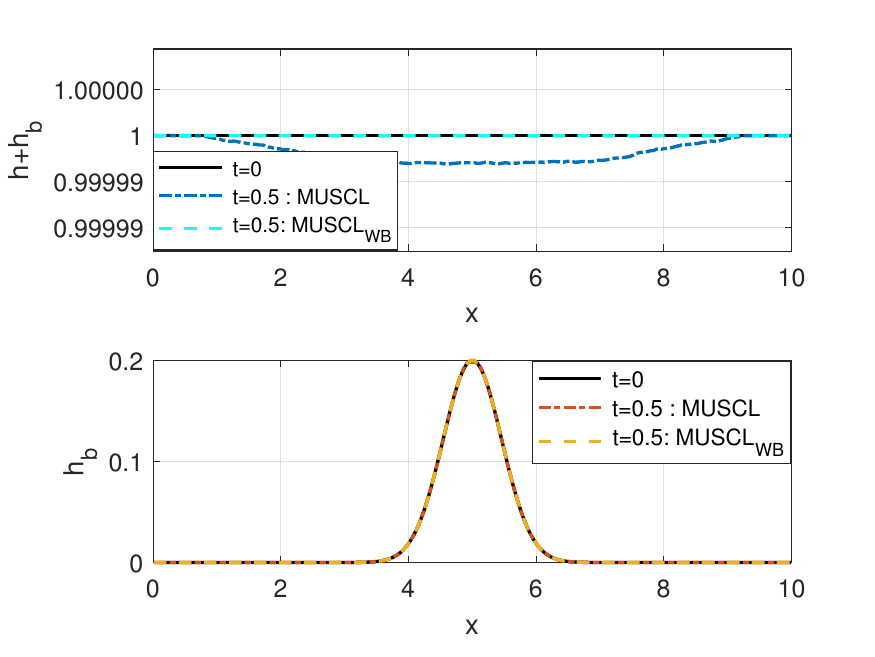}
		\caption*{(a)}
	\end{subfigure}
	\hfill
	\begin{subfigure}[b]{0.49\textwidth}
		\includegraphics[width=\textwidth]{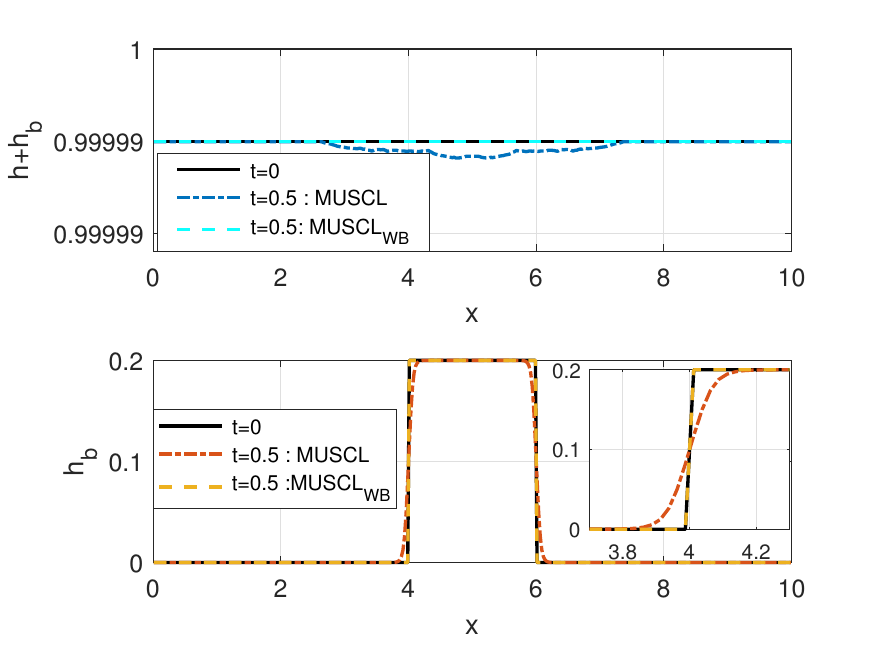}
		\caption*{(b)}
	\end{subfigure}
	\caption{Well-balanced verification at $t=0.5$: (a) smooth Gaussian bottom \eqref{war-smooth-ic} and (b) discontinuous rectangular bottom \eqref{war-rectangle-ic}. The upper panels show the free surface $h+h_b$, while the lower panels show the sediment bottom $h_b$. The corrected scheme $(\text{MUSCL}_{\mathrm{WB}})$ preserves the constant free surface and the stationary bottom, while the standard MUSCL-Rusanov scheme perturbs the equilibrium.}
	\label{Water-at-rest-cor}
\end{figure}
Figure~\ref{Water-at-rest-cor}(a) shows the free surface $h+h_b$ and sediment bottom $h_b$ at $t=0$  and $t=0.5$ for the smooth Gaussian bottom. Both the standard MUSCL-Rusanov scheme and the corrected well-balanced scheme preserve the bottom profile. However, from Figure~\ref{Water-at-rest-cor}(a) the free-surface shows that the standard MUSCL-Rusanov scheme introduces a small deviation from the constant free surface $h+h_b$. Therefore, although the error remains small for the smooth bottom, the standard MUSCL-Rusanov scheme does not preserve the water-at-rest equilibrium exactly.

Figure~\ref{Water-at-rest-cor}(b) shows the free surface $h+h_b$ and sediment bottom $h_b$ at $t=0$  and $t=0.5$ for the discontinuous rectangular bottom. At $x=4$ and $x=6$, the generalized minmod limiter reduces the reconstruction to piecewise constant states in the cells adjacent to the discontinuities. Consequently, the complete stationary jump $\Delta h=-\Delta h_b$ remains between the reconstructed states at these interfaces. A zoom at the sediment bottom shows that the non well-balanced Rusanov artificial viscosity \eqrefc{Rusanov artificial viscosity} acts on this non-zero jump and produces diffusion of the sediment bottom and also a perturbation of the free surface.
The proposed well-balanced correction \eqrefc{corrected Rusanov} removes the water-at-rest component from the jump used in the artificial viscosity. The corrected MUSCL-Rusanov scheme $(\text{MUSCL}_{\text{WB}})$ therefore preserves both the constant free surface and the discontinuous rectangular bottom. This test confirms the necessity of the correction and verifies that the resulting scheme is well-balanced for the water-at-rest equilibrium, including stationary discontinuities of the bottom.

\subsection{Accuracy test of the MUSCL-Rusanov scheme}\label{accuracy_test}
Here we perform an accuracy test for the second-order MUSCL-Rusanov scheme used in this work. Following the smooth accuracy test introduced in \cite{Garres}, we consider the computational domain $\Omega=[0,10]$ with periodic boundary conditions. The initial sediment bottom $h_b$, water height $h$, and discharge $hu_m$ are given by
\begin{equation}\label{convergence_initial_hhb}
	h_b(0,x)=\frac{1}{5}\exp\!\left[-2(x-5)^2\right],
	\qquad
	h(0,x)=2-\frac{6}{5}h_b(0,x),
	\qquad
	hu_m(0,x)=1.5.
\end{equation}
We consider the initial velocity profile in the same quadratic form as in \cite{Garres},
\begin{equation}\label{convergence_initial_velocity}
	u(0,x,\xi)=a_0+a_1\xi+a_2\xi^2,
\end{equation}
where
\begin{equation*}
	a_2=0.35,
	\qquad
	a_1=\frac{11}{15},
	\qquad
	a_0=\frac{1.5}{h(0,x)}-\frac{a_2}{3}-\frac{a_1}{2}.
\end{equation*}
Expressing the quadratic profile \eqrefc{convergence_initial_velocity} in the Legendre basis used for the moment expansion gives
\begin{equation*}
	u_m(0,x)=\frac{1.5}{h(0,x)},
	\qquad
	\alpha_1(0,x)=-\frac{13}{24},
	\qquad
	\alpha_2(0,x)=\frac{7}{120},
	\qquad
	\alpha_3(0,x)=0.
\end{equation*}

To extend the test to the complete HSWEMED system \eqrefc{final HSWEMED}, we consider a non-zero initial suspended sediment concentration,
\begin{equation}\label{convergence_initial_concentration}
	c_m(0,x)
	=
	5\times10^{-3}
	+
	10^{-3}\exp\!\left[-2(x-5)^2\right].
\end{equation}
The sediment properties used in this test are presented in Table~\ref{tab:accuracy_sediment_properties}.
\begin{table}[H]
	\centering
	\begin{tabular}{|c|c|c|c|c|c|}
		\hline
		$\rho_w$
		&
		$\rho_s$
		&
		$d_s$
		&
		$\theta_c$
		&
		$\psi$
		&
		$\epsilon$
		\\
		\hline
		$1000$
		&
		$2941.18$
		&
		$1$
		&
		$0.047$
		&
		$0.95$
		&
		$0.09$
		\\
		\hline
	\end{tabular}
	\caption{Sediment properties used in the smooth accuracy test \cite{Garres}.}
	\label{tab:accuracy_sediment_properties}
\end{table}

\begin{figure}[H]
	\centering
	\begin{subfigure}[b]{0.49\textwidth}
		\includegraphics[width=\textwidth]
		{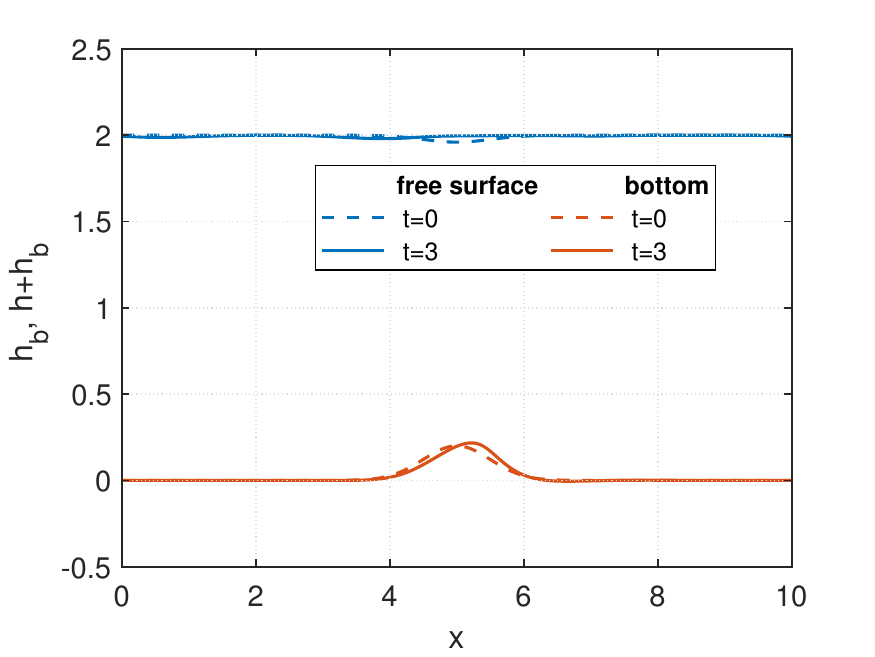}
		\caption*{(a)}
	\end{subfigure}
	\hfill
	\begin{subfigure}[b]{0.49\textwidth}
		\includegraphics[width=\textwidth]
		{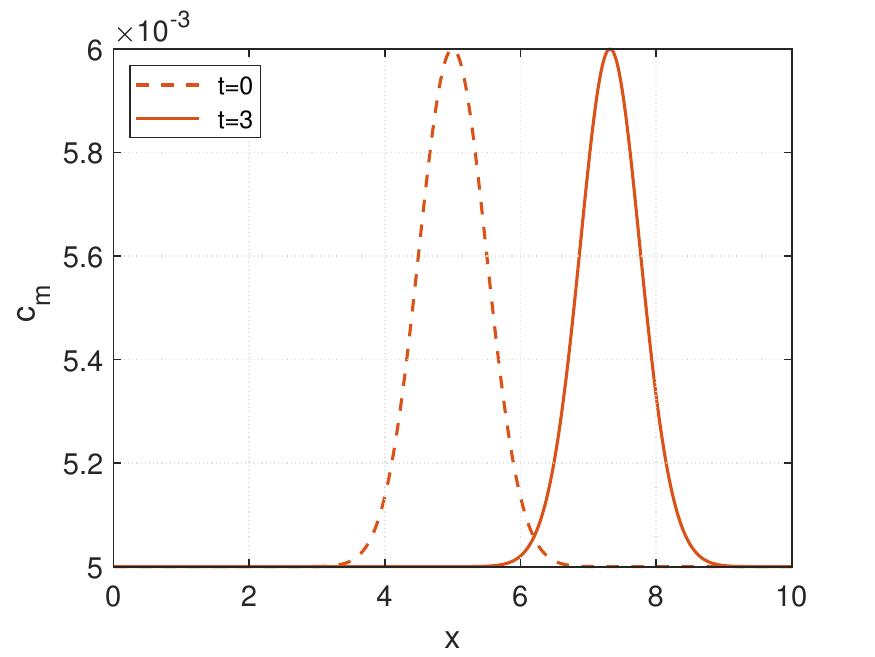}
		\caption*{(b)}
	\end{subfigure}
	\caption{(a) Free surface $h+h_b$ and sediment bottom $h_b$, and (b) Suspended sediment concentration $c_m$ at time $t=0$ and $t=3$.}
	\label{fig:convergence_results}
\end{figure}

We consider the third-order HSWEMED, i.e., $N=3$, for this test. The friction, erosion, and deposition source terms are omitted so that the test measures the accuracy of the path-conservative MUSCL-Rusanov transport discretization independently of the source treatment. We compute the solutions until $t=3$ using grids
$N_x=25,50,100,200,400,800,\text{ and }1600$. We compute a fixed numerical reference solution using $N_{x,\mathrm{ref}}=3200$.
Figure~\ref{fig:convergence_results}(a) shows the initial and final free surface $h+h_b$ and sediment bottom $h_b$ on the reference grid $N_{x,\mathrm{ref}}=3200$ at time $t=3$. We observe that the initially smooth bottom profile evolves with time and remains smooth. Figure~\ref{fig:convergence_results}(b) shows that the initial concentration perturbation is transported downstream at $t=3$.

The $L^1$ errors and computed convergence orders are presented in Table~\ref{tab:convergence_l1}. The coarser grids remain outside the asymptotic convergence range. From $N_x=200$ onward, the errors decrease consistently and the computed convergence orders are generally close to two for the water height $h$, discharge $hu_m$, three moment variables $(h\alpha_1,h\alpha_2,h\alpha_3)$, suspended sediment concentration $hc_m$, and sediment bottom $h_b$.
Overall, the results are consistent with the expected second-order accuracy of the path-conservative MUSCL-Rusanov scheme for the complete HSWEMED state vector $W$.

\begin{table}[H]
	\centering
	\setlength{\tabcolsep}{3.2pt}
	\renewcommand{\arraystretch}{1.12}
	\begin{tabular}{|c|cc|cc|cc|cc|}
		\hline
		$N_x$
		&
		\multicolumn{2}{c|}{$h$}
		&
		\multicolumn{2}{c|}{$hu_m$}
		&
		\multicolumn{2}{c|}{$hc_m$}
		&
		\multicolumn{2}{c|}{$h_b$}
		\\
		\cline{2-9}
		&
		Error
		&
		Order
		&
		Error
		&
		Order
		&
		Error
		&
		Order
		&
		Error
		&
		Order
		\\
		\hline
		25
		&
		$2.023\times10^{-1}$
		&
		--
		&
		$1.044\times10^{-1}$
		&
		--
		&
		$2.572\times10^{-3}$
		&
		--
		&
		$1.967\times10^{-1}$
		&
		--
		\\
		50
		&
		$1.102\times10^{-1}$
		&
		0.88
		&
		$6.399\times10^{-2}$
		&
		0.71
		&
		$1.338\times10^{-3}$
		&
		0.94
		&
		$1.068\times10^{-1}$
		&
		0.88
		\\
		100
		&
		$4.058\times10^{-2}$
		&
		1.44
		&
		$2.128\times10^{-2}$
		&
		1.59
		&
		$5.214\times10^{-4}$
		&
		1.36
		&
		$3.609\times10^{-2}$
		&
		1.56
		\\
		200
		&
		$1.040\times10^{-2}$
		&
		1.96
		&
		$5.670\times10^{-3}$
		&
		1.91
		&
		$1.349\times10^{-4}$
		&
		1.95
		&
		$9.477\times10^{-3}$
		&
		1.93
		\\
		400
		&
		$2.205\times10^{-3}$
		&
		2.24
		&
		$1.590\times10^{-3}$
		&
		1.83
		&
		$2.771\times10^{-5}$
		&
		2.28
		&
		$1.996\times10^{-3}$
		&
		2.25
		\\
		800
		&
		$4.606\times10^{-4}$
		&
		2.26
		&
		$3.986\times10^{-4}$
		&
		2.00
		&
		$5.489\times10^{-6}$
		&
		2.34
		&
		$3.978\times10^{-4}$
		&
		2.33
		\\
		1600
		&
		$8.069\times10^{-5}$
		&
		2.51
		&
		$8.050\times10^{-5}$
		&
		2.31
		&
		$9.817\times10^{-7}$
		&
		2.48
		&
		$6.644\times10^{-5}$
		&
		2.58
		\\
		\hline
	\end{tabular}
	
	\vspace{1.5mm}
	
	\begin{tabular}{|c|cc|cc|cc|}
		\hline
		$N_x$
		&
		\multicolumn{2}{c|}{$h\alpha_1$}
		&
		\multicolumn{2}{c|}{$h\alpha_2$}
		&
		\multicolumn{2}{c|}{$h\alpha_3$}
		\\
		\cline{2-7}
		&
		Error
		&
		Order
		&
		Error
		&
		Order
		&
		Error
		&
		Order
		\\
		\hline
		25
		&
		$3.511\times10^{-1}$
		&
		--
		&
		$8.802\times10^{-2}$
		&
		--
		&
		$7.584\times10^{-2}$
		&
		--
		\\
		50
		&
		$1.886\times10^{-1}$
		&
		0.90
		&
		$5.410\times10^{-2}$
		&
		0.70
		&
		$5.975\times10^{-2}$
		&
		0.34
		\\
		100
		&
		$6.962\times10^{-2}$
		&
		1.44
		&
		$1.643\times10^{-2}$
		&
		1.72
		&
		$2.599\times10^{-2}$
		&
		1.20
		\\
		200
		&
		$1.796\times10^{-2}$
		&
		1.95
		&
		$4.002\times10^{-3}$
		&
		2.04
		&
		$6.854\times10^{-3}$
		&
		1.92
		\\
		400
		&
		$4.530\times10^{-3}$
		&
		1.99
		&
		$1.019\times10^{-3}$
		&
		1.97
		&
		$1.467\times10^{-3}$
		&
		2.22
		\\
		800
		&
		$1.082\times10^{-3}$
		&
		2.07
		&
		$1.787\times10^{-4}$
		&
		2.51
		&
		$3.833\times10^{-4}$
		&
		1.94
		\\
		1600
		&
		$2.427\times10^{-4}$
		&
		2.16
		&
		$3.306\times10^{-5}$
		&
		2.43
		&
		$1.030\times10^{-4}$
		&
		1.90
		\\
		\hline
	\end{tabular}
	\caption{$L^1$ errors and computed convergence orders for the conservative HSWEMED variables $W$ at $t=3$. Each reported grid is compared with the fixed reference solution with $N_{x,\mathrm{ref}}=3200$.}
	\label{tab:convergence_l1}
\end{table}

\subsection{Academic dam-break test}\label{academic dam-break1}
For the academic dam-break test case, we use the following initial conditions for the water height $h$ and the sediment bottom $h_b$.
\begin{equation}\label{academic dam-break_eq}
	\vcenter{\hbox{%
		\shortstack{
			\stepcounter{figure}
			\includeinkscape[width=0.38\textwidth]{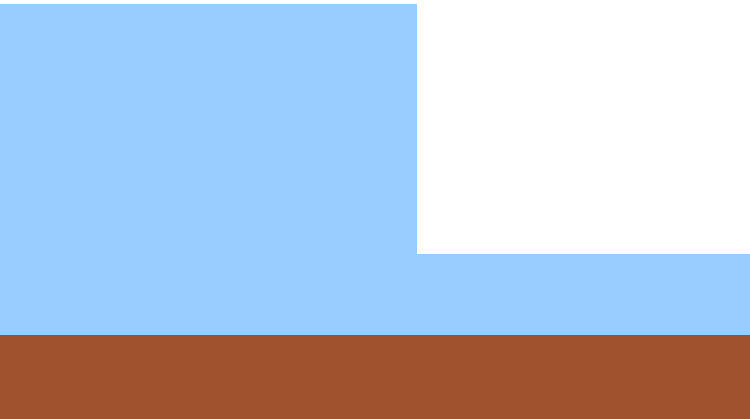}
				\\[2ex]
				\small Figure \thefigure: Initial water height and sediment bottom for the academic dam-break test \cite{Garres}.
			}
	}}
	\qquad
	\begin{aligned}
		h(0,x)
		&=
		\begin{cases}
			1,    & \text{if } x \leq 0,\\
			0.05, & \text{otherwise},
		\end{cases}\\
		h_b(0,x) &= 0.
	\end{aligned}
\end{equation}
The properties of the sediment particles considered in this test are shown in \tabref{table:ACADEMIC TEST}. The computational domain is $\Omega=[-6,6]$, divided into $N_x=1200$ cells, and the solution is computed until $t=1$.
\begin{table}[H]
	\centering
	\begin{tabular}{|c|c|c|c|c|c|}
		\hline
		$\rho_w$ & $\rho_s$ &$\theta_c$ &$ \psi $& $d_s $ & $ \epsilon $ \\ \hline
		1000 & 1580 & 0.047 & 0.47 & 3.9 & 0.0324 \\ \hline
	\end{tabular}
	\caption{Sediment properties of academic dam-break test \cite{Garres}.}
	\label{table:ACADEMIC TEST}
\end{table}
In this test case, we consider the third-order HSWEMED, i.e., $N=3$. As suggested in \cite{Garres}, this choice allows for approximating the vertical velocity profile using a polynomial of degree three. It provides a more accurate representation of the velocity distribution near the bed than a linear approximation while maintaining computational efficiency. For consistency, the same value of $N$ is used for the HSWEM, whose results are compared with those of the HSWEMED.

\begin{figure}[H]
	\centering
	\begin{subfigure}[b]{0.49\textwidth}
		\includegraphics[width=\textwidth]{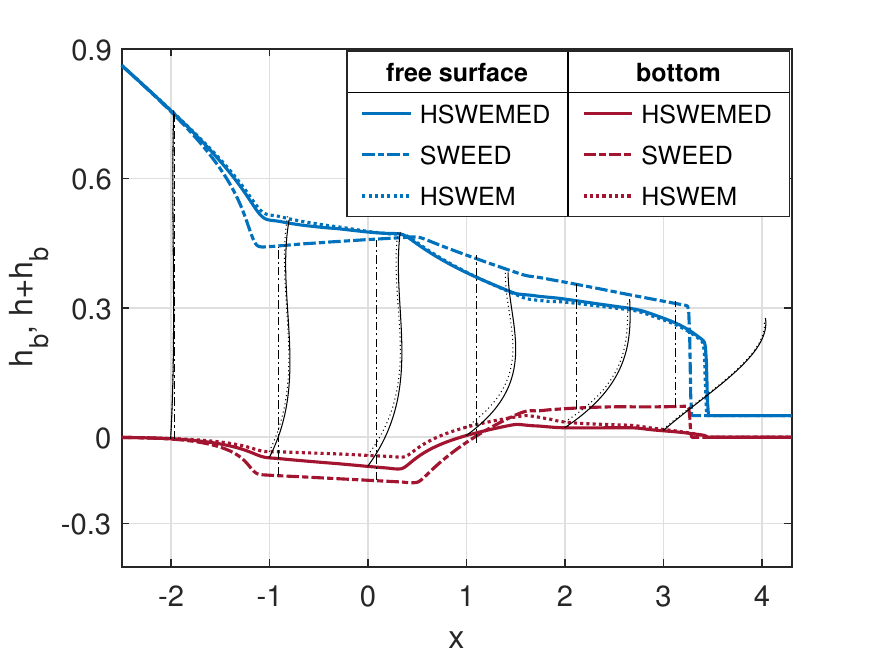}
		\caption*{(a)}
	\end{subfigure}
	\hfill
	\begin{subfigure}[b]{0.49\textwidth}
		\includegraphics[width=\textwidth]{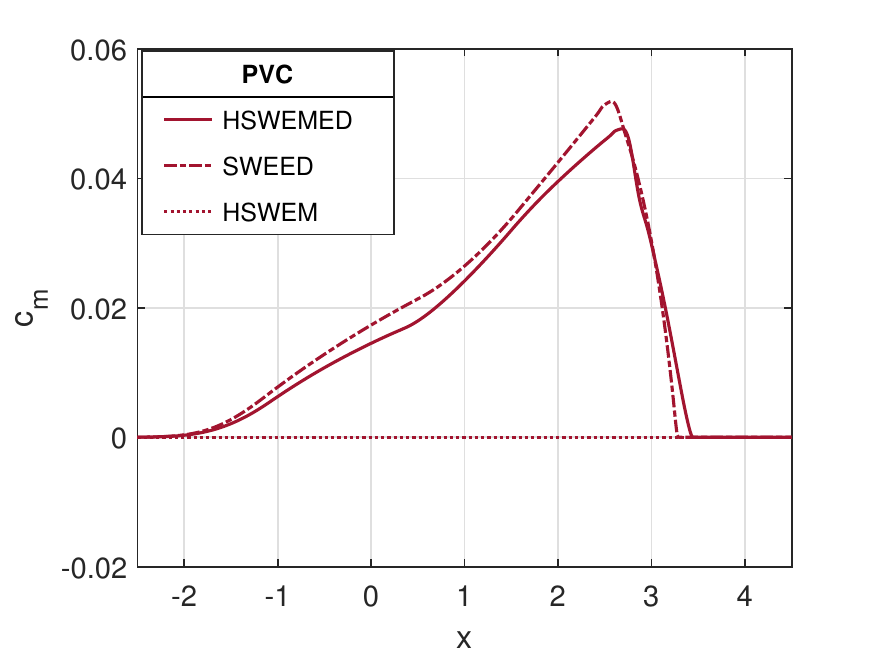}
		\caption*{(b)}
	\end{subfigure}
	\caption{Academic dam-break test \eqrefc{academic dam-break_eq} at $t=1$: (a) free surface $h+h_b$, bottom evolution $h_b$, and vertical velocity profiles $u$ (black lines) at $x=-2,-1,0,1,2,3$; and (b) volumetric sediment concentration $c_m$. Results are computed with the HSWEMED (solid), SWEED (dashed), and HSWEM (dotted).}
	\label{comparison_academic_dam_test}
\end{figure}

Figure~\ref{comparison_academic_dam_test}(a) indicates that the free surface $h+h_b$ and sediment bottom evolution $h_b$ predicted by the HSWEMED differ from those of the other two models. Specifically, when examining bottom evolution $h_b$, the HSWEMED predicts more erosion than the HSWEM and less erosion than the SWEED at the dam-break location. Furthermore, a significant fraction of the sediment remains suspended further upstream. This improvement comes from incorporating the bottom velocity $u_b$ \eqrefc{bottom_velocity}, which is lower than the depth-averaged velocity $u_m$ used in the SWEED. As expected, the SWEED tends to overestimate erosion, which explains its smaller $h_b$ around $x=0$ in the erosion domain. In contrast, the HSWEM only accounts for bedload transport without considering erosion and deposition effects, resulting in an inaccurate representation of bottom changes \cite{Garres}. The lack of erosion and deposition modeling in the HSWEM leads to underestimating the eroded material and gives the largest $h_b$ around $x=0$. Our proposed HSWEMED falls between the two, capturing a moderate level of erosion.

Regarding the free-surface evolution, Figure~\ref{comparison_academic_dam_test}(a) shows that the HSWEMED front propagates slightly farther downstream than the SWEED front. The vertical velocity profiles also show that the HSWEMED recovers an approximately parabolic velocity distribution, especially near the wavefront, with a bottom velocity $u_b$ smaller than the depth-averaged velocity $u_m$ used in the SWEED. The HSWEM result remains close to the HSWEMED near the wavefront. Since both HSWEM and HSWEMED models include the velocity moments, the small difference in the wavefront location is mainly related to the moment representation of the vertical velocity. In contrast, erosion and deposition have a stronger influence on the bed evolution behind the wavefront.

Figure~\ref{comparison_academic_dam_test}(b) illustrates the volumetric sediment concentration $c_m$. Since the HSWEM exclusively models bedload transport, its concentration remains zero. Although the HSWEMED and SWEED use the same settling velocity $\omega_o$, the SWEED yields a slightly higher suspended sediment concentration. This discrepancy arises because the SWEED induces more sediment erosion near the dam-break location due to the use of the depth-averaged velocity.

\subsection{Laboratory experiments: Dam-break flow over erodible sediment bed}\label{laboratory experiments}
In this section, we simulate one-dimensional dam-break flow over an erodible sediment bed and compare the HSWEMED results with the laboratory experiments documented in \cite{Spinewine01122007}. We consider two initial configurations:
\begin{itemize}
	\item Config 1 \eqrefc{exp_1}: dam-break flow with a wet/dry front over a flat bottom (config a in \cite{Spinewine01122007});
	\item Config 2 \eqrefc{exp_2}: dam-break flow over a discontinuous bottom with water at rest downstream (config f in \cite{Spinewine01122007}).
\end{itemize}
For both configurations, we compare the free surface $h+h_b$, the sediment-bottom evolution $h_b$, and the depth-averaged suspended concentration $c_m$ obtained with the following models:
\begin{enumerate}
	\item HSWEMED: the proposed model \eqrefc{complete system}, containing $N+4$ equations;
	\item SWEED: the standard depth-averaged model \eqrefc{zeroth_order_moment}, containing four equations;
	\item HSWEM: the moment model without erosion and deposition effects \cite{Garres}, containing $N+3$ equations.
\end{enumerate}
The comparison between HSWEMED and SWEED shows the effect of the moment approximation, whereas the comparison with HSWEM shows the effect of erosion, deposition, and suspended sediment concentration.
For both configurations, two different sediment bed materials are used: PVC pellets and uniform coarse sand. Their different densities and transport properties influence the propagation speed, suspended sediment concentration, and erosion process. The properties of the bed materials are taken from \cite{Spinewine01122007} and shown in \tabref{table:exeriments}.
\begin{table}[H]
	\centering
	\begin{tabular}{|c|c|c|c|c|c|c|}
		\hline
		Materials & $\rho_w$ & $\rho_s$ &$\theta_c$ &$ \psi $& $d_s $ & $ \epsilon $ \\ \hline
		PVC       & 1000 & 1580 & 0.047 & 0.47 & 3.9  & 0.0324 \\ \hline
		Sand      & 1000 & 2683 & 0.047 & 0.47 & 1.82 & 0.0104 \\ \hline
	\end{tabular}
	\caption{Physical parameters of laboratory experiments.}
	\label{table:exeriments}
\end{table}

\vspace{-10mm}
\subsubsection{Config 1: Dam-break flow with wet/dry front over flat bottom}\label{conf_1}
For the first configuration, we use the following initial conditions for the water height $h,$ and sediment bed $h_b$
\begin{equation}\label{exp_1}
	\vcenter{\hbox{%
			\shortstack{
				\refstepcounter{figure}
				\includeinkscape[width=0.4\textwidth]{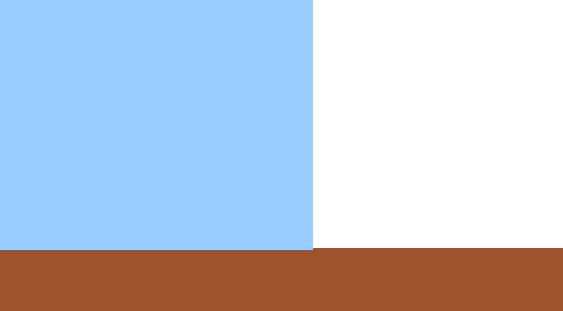}
				\\[1ex]
				\small Figure \thefigure: Initial water height and bottom for laboratory Config 1\\
				(config a in \cite{Spinewine01122007}).
			}
	}}
	\qquad
	\begin{aligned}
		h(0,x)
		&=
		\begin{cases}
			0.35, & \text{if } x \le 0,\\
			0,    & \text{otherwise,}
		\end{cases}\\
		h_b(0,x) &= 0.
	\end{aligned}
\end{equation}
The computational domain $\Omega=[-3,3]$ is divided into $N_x=1000$ grid points. Within this framework, we apply a wet/dry treatment in the HSWEMED by considering a computational cell as dry when the water height $h$ falls below a specified threshold, specifically when $h<\delta$, with $\delta=10^{-4}$, to maintain numerical stability.

We consider the first-order HSWEMED in this configuration, i.e., $N=1$. In near-dry regions, where the water height $h$ is nearly zero, the friction term in the momentum equation \eqrefc{depth_avg_momentum} and the higher-order averaged equation \eqrefc{final_higher_average_equation} becomes excessively large due to the presence of $h$ in the denominator. Such behavior introduces significant numerical stiffness, leading to instability unless extremely small time steps are considered. Therefore, incorporating higher-order moments poses computational challenges. We defer this to future research, where an asymptotic-preserving numerical scheme could be applied, for example, using projective integration \cite{amrita2022projective} or a time-splitting scheme \cite{huang2022equilibrium}.

\begin{figure}[H]
	\centering
	\begin{subfigure}[b]{0.49\textwidth}
		\includegraphics[width=\textwidth]{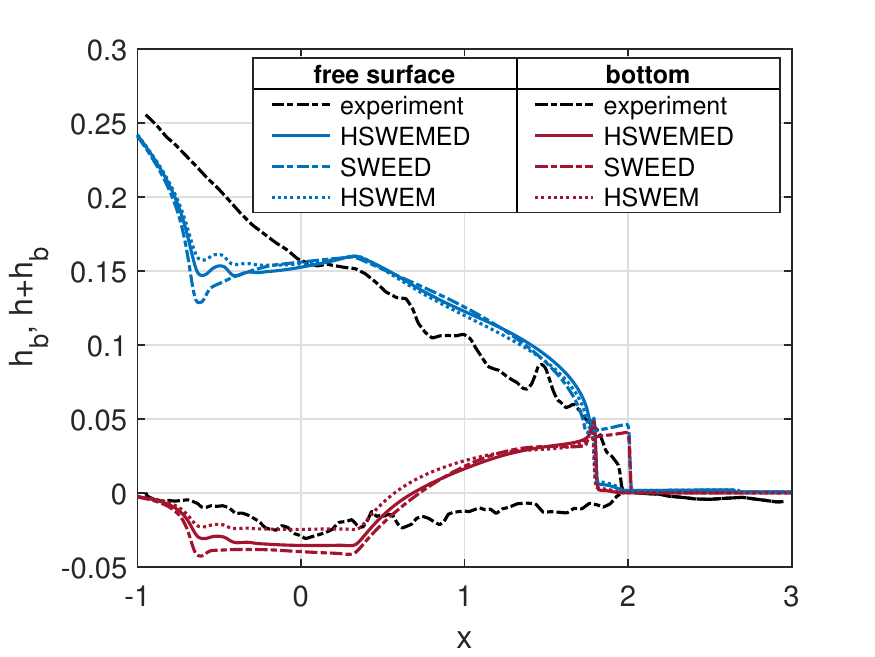}
		\caption*{(a)}
	\end{subfigure}
	\hfill
	\begin{subfigure}[b]{0.49\textwidth}
		\includegraphics[width=\textwidth]{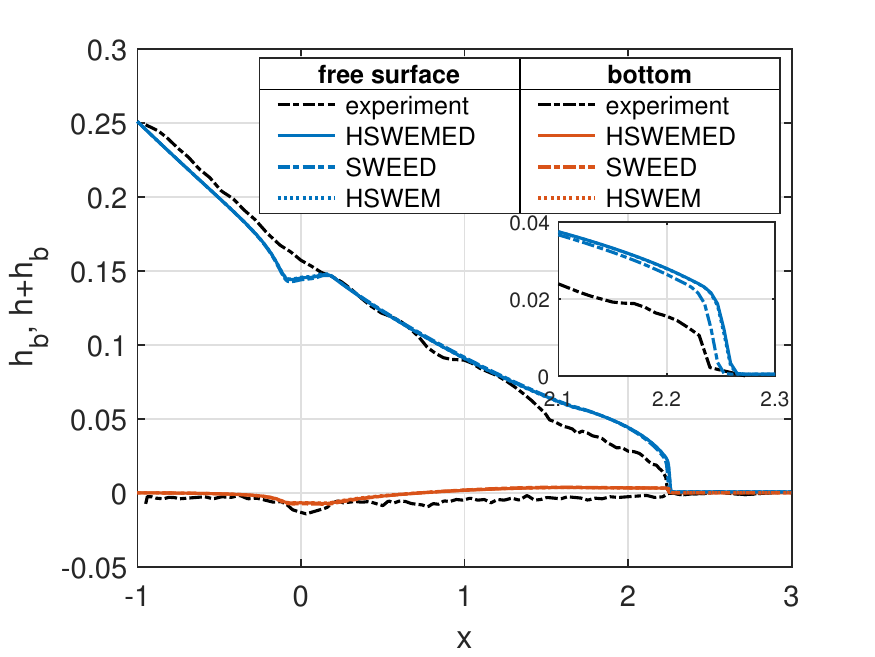}
		\caption*{(b)}
	\end{subfigure}
	\begin{subfigure}[b]{0.49\textwidth}
		\includegraphics[width=\textwidth]{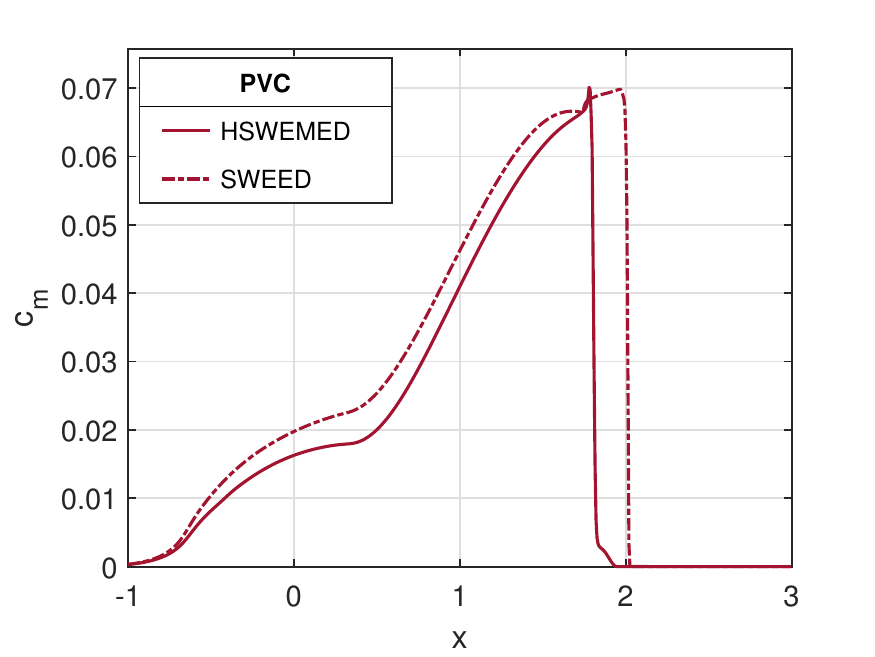}
		\caption*{(c)}
	\end{subfigure}
	\hfill
	\begin{subfigure}[b]{0.49\textwidth}
		\includegraphics[width=\textwidth]{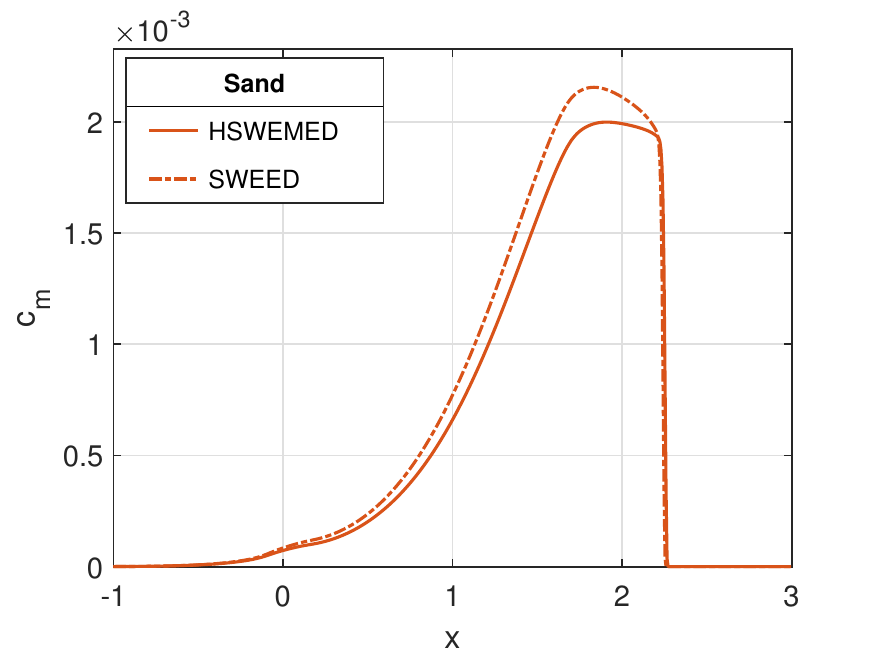}
		\caption*{(d)}
	\end{subfigure}
	\caption{Laboratory Config 1 \eqrefc{exp_1} at $t=1$: free surface $h+h_b$ and bottom $h_b$ for (a) PVC pellets and (b) uniform coarse sand, computed with the HSWEMED (solid), SWEED (dashed), and HSWEM (dotted), together with the experimental data (dash-dotted); and volumetric sediment concentration $c_m$ for (c) PVC pellets and (d) uniform coarse sand, computed with the HSWEMED (solid) and SWEED (dashed). The inset in (b) enlarges the free-surface front.}
	\label{Exp_pvc}
\end{figure}

Figure~\ref{Exp_pvc}(a) and Figure~\ref{Exp_pvc}(b) present the numerical results for the free surface $h+h_b$ and bottom evolution $h_b$ at $t=1$, as computed by the first-order HSWEMED, the SWEED, and the first-order HSWEM. These results are compared with experimental data for PVC and coarser sand beds from \cite{Spinewine01122007}. Figure~\ref{Exp_pvc}(c) and Figure~\ref{Exp_pvc}(d) show the corresponding sediment concentration $c_m$ computed with the HSWEMED and SWEED, since the HSWEM only models bedload transport.

For the PVC bed, Figure~\ref{Exp_pvc}(a) demonstrates that the sediment entrainment predicted by the HSWEMED, both upstream and downstream, aligns more closely with the experimental data compared with the SWEED. Specifically, the SWEED overestimates bottom erosion in the dam-break region. The difference between the HSWEMED and SWEED can be associated with the bottom velocity $u_b$. The erosion rate and sediment transport in the HSWEMED are computed using $u_b=u_m+\sum_{j=1}^{N}\alpha_j$, which is lower than the depth-averaged velocity $u_m$ used by the SWEED. In contrast, slight differences are observed compared with the HSWEM, which neglects erosion and deposition and only models bedload transport. Additionally, the HSWEMED accurately captures the location of the hydraulic jump around $x=0$, although it does not perfectly reproduce the experimental free surface in this region.

For the sand bed, Figure~\ref{Exp_pvc}(b) shows that the three numerical models HSWEMED, SWEED, and HSWEM produce very similar free-surface $h+h_b$ and bottom evolution $h_b$. The overall decrease of the free surface is captured, although the numerical solutions remain below the experimental data upstream and above it farther downstream. Near the dam-break location, the models reproduce the local change in the free surface, but the variation is sharper than in the experiment. The inset in Figure~\ref{Exp_pvc}(b) shows that all three models place the wavefront slightly farther downstream than the experiment. The HSWEMED and HSWEM wavefronts almost coincide and are located slightly ahead of the SWEED wavefront. This suggests that the small difference in the wavefront location is mainly related to the moment representation of the vertical velocity.

For the volumetric sediment concentration $c_m$,
Figures~\ref{Exp_pvc}(c) and~\ref{Exp_pvc}(d) show a clear difference between the PVC and sand beds. The suspended sediment concentration is considerably higher for PVC pellets than for sand particles. This difference is related to the lower density of PVC pellets, with $\rho_s/\rho_w=1.58$, compared with $\rho_s/\rho_w=2.683$ for sand. In addition, the settling velocity of PVC pellets, $\omega_o=0.1528\,\mathrm{m\,s^{-1}}$, computed from equation~\eqrefc{settling velocity}, is lower than that of sand particles, $\omega_o=0.5714\,\mathrm{m\,s^{-1}}$. Therefore, a larger amount of PVC sediment remains in suspension, with concentration values of order $10^{-2}$, whereas the sand concentration remains of order $10^{-3}$.
For both materials, the SWEED generally predicts slightly more suspended sediment than the HSWEMED over most of the suspension region. This is consistent with the stronger bottom erosion predicted by the SWEED. The difference between the HSWEMED and SWEED concentration profiles is more visible for the PVC bed than for the sand bed.
At $x\approx1.6$, the local HSWEMED concentrations are approximately $c_m=0.0640$ for PVC and $c_m=1.76\times10^{-3}$ for sand. We compare the eigenvalues at the same location using the sediment coupled factor of the characteristic polynomial \eqrefc{CP_HSWEMED}. Table~\ref{tab:hswemed_roots} shows that the largest positive eigenvalue is higher for the sand case than for the PVC case, indicating faster wave propagation for sand particles. This is consistent with the sand wavefront being located farther downstream in Figure~\ref{Exp_pvc}(b).
\begin{table}[h]
	\centering
	\renewcommand{\arraystretch}{1.5}
	\setlength{\tabcolsep}{15pt}
	\begin{tabular}{|c|c|c|}
		\hline
		\textbf{Eigenvalue}
		& \textbf{PVC} & \textbf{Sand} \\
		\hline
		$\lambda_1$ & $-0.4350$ & $-0.0619$ \\
		\hline
		$\lambda_2$ & $1.0148$ & $1.0786$ \\
		\hline
		$\lambda_3$ & $1.8333$ & $2.4437$ \\
		\hline
	\end{tabular}
	\caption{Eigenvalues obtained from the sediment coupled factor of \eqrefc{CP_HSWEMED}, evaluated around $x\approx1.6$ at $t=1$ for the HSWEMED solutions with PVC and sand beds under Config 1 \eqrefc{exp_1}. The larger positive eigenvalue for sand indicates faster wave propagation for sand particles.}
	\label{tab:hswemed_roots}
\end{table}

\subsubsection{Config 2: Dam-break flow over discontinuous bottom with water in the downstream}\label{conf_2}
For the second configuration, the initial water height $h$ and sediment bed $h_b$ are given by
\begin{equation}\label{exp_2}
	\vcenter{\hbox{%
			\shortstack{
				\refstepcounter{figure}
				\includeinkscape[width=0.4\textwidth]{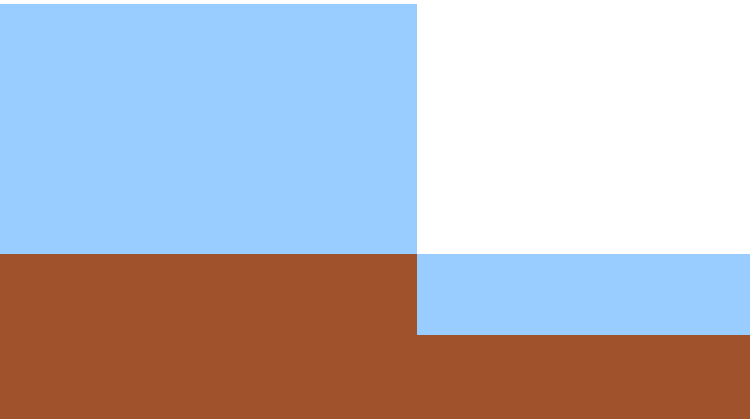}
				\\[1ex]
				\small Figure \thefigure: Initial water height and bottom for laboratory Config 2\\
				(config f in \cite{Spinewine01122007}).
			}
	}}
	\quad
	\begin{aligned}
		h(0,x)
		&=
		\begin{cases}
			0.25, & \text{if } x \le 0,\\
			0.10, & \text{otherwise,}
		\end{cases}\\
		h_b(0,x)
		&=
		\begin{cases}
			0.10, & \text{if } x \le 0,\\
			0,    & \text{otherwise.}
		\end{cases}
	\end{aligned}
\end{equation}

We consider $\Omega=[-3,3]$ as the computational domain, discretized into $N_x=1000$ equally spaced points. The main difference in this setup compared with the first configuration is the presence of a wet downstream.

In this test case, we again consider the third-order HSWEMED, i.e., $N=3$, consistent with the value used for the HSWEM in \cite{Garres}.

\begin{figure}[H]
	\centering
	\begin{subfigure}[b]{0.49\textwidth}
		\includegraphics[width=\textwidth]{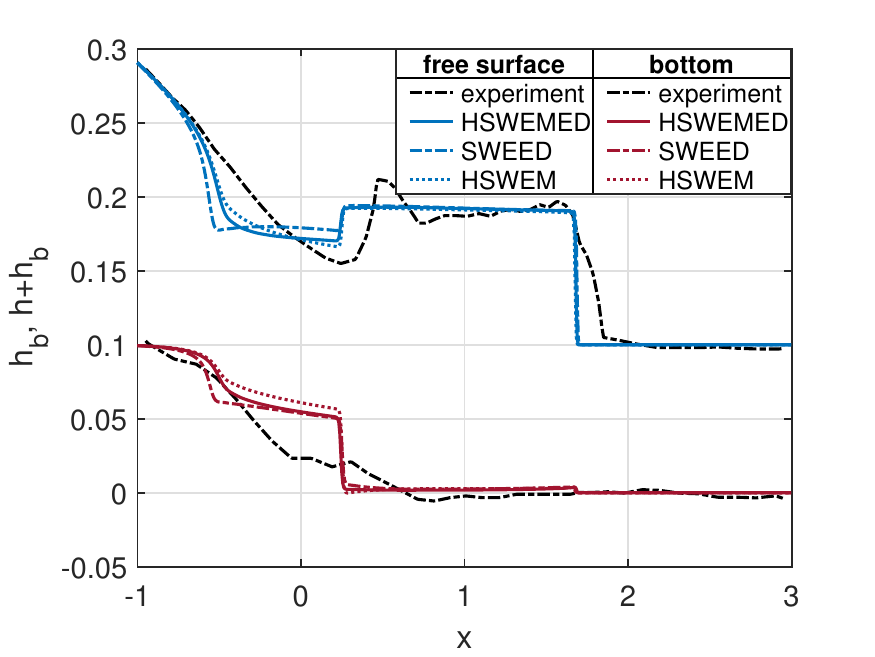}
		\caption*{(a)}
	\end{subfigure}
	\hfill
	\begin{subfigure}[b]{0.49\textwidth}
		\includegraphics[width=\textwidth]{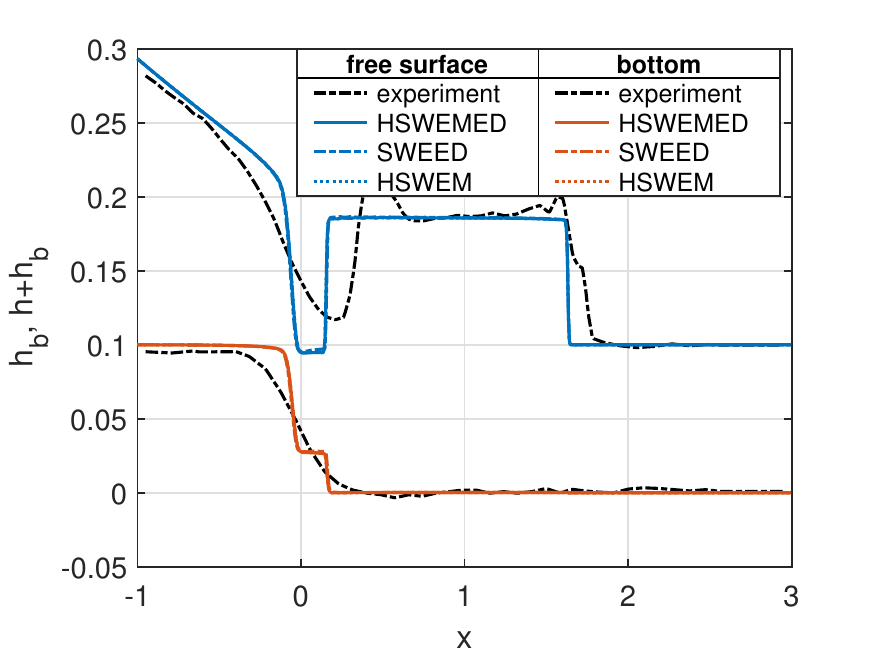}
		\caption*{(b)}
	\end{subfigure}
	\begin{subfigure}[b]{0.49\textwidth}
		\includegraphics[width=\textwidth]{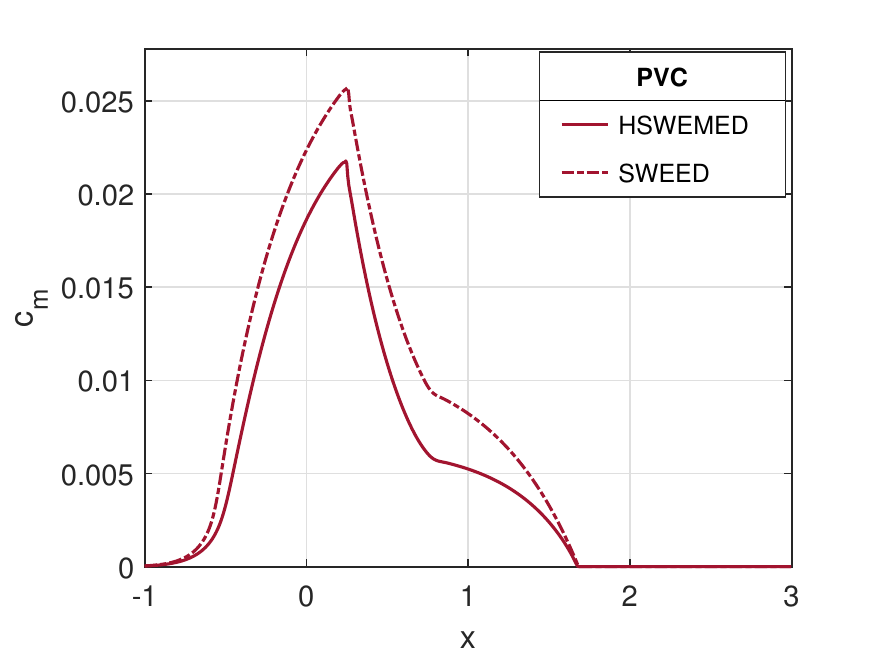}
		\caption*{(c)}
	\end{subfigure}
	\hfill
	\begin{subfigure}[b]{0.49\textwidth}
		\includegraphics[width=\textwidth]{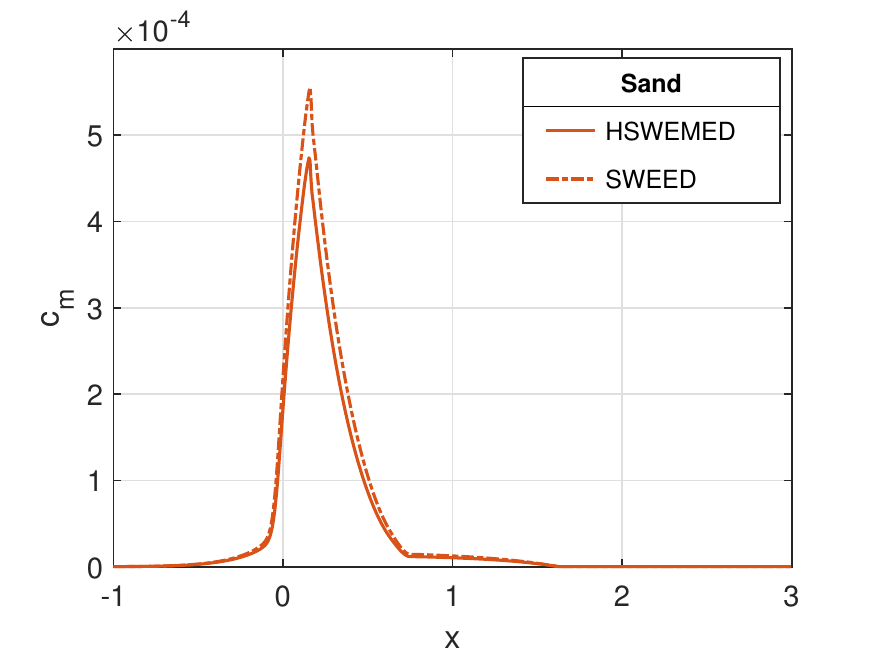}
		\caption*{(d)}
	\end{subfigure}
	\caption{Laboratory Config 2 \eqrefc{exp_2} at $t=1$: free surface $h+h_b$ and bottom $h_b$ for (a) PVC pellets and (b) uniform coarse sand, computed with the HSWEMED (solid), SWEED (dashed), and HSWEM (dotted), together with the experimental data (dash-dotted); and volumetric sediment concentration $c_m$ for (c) PVC pellets and (d) uniform coarse sand, computed with the HSWEMED (solid) and SWEED (dashed).}
	\label{exp3_pvc}
\end{figure}

Figure~\ref{exp3_pvc}(a) and Figure~\ref{exp3_pvc}(b) present the numerical results for the free surface $h+h_b$ and bottom evolution $h_b$ at $t=1$ for PVC and sand, respectively. Figure~\ref{exp3_pvc}(c) and Figure~\ref{exp3_pvc}(d) show the corresponding concentration profiles.

For the PVC bed, Figure~\ref{exp3_pvc}(a) illustrates that the free surface $h+h_b$ in the vicinity of the hydraulic jump around $x=0$ is only partially captured by the three models. This discrepancy arises from a deviation between the experimentally observed and numerically predicted positions of the free surface and may be related to the hydrostatic shallow-water assumption near the dam-break location, $x=0$. However, the models reproduce the downstream waterfront and the bottom evolution. The three numerical solutions remain close over most of the computational domain, with the largest differences observed at the hydraulic-jump region, $x= 0$.

\begin{figure}[H]
	\centering
	\begin{subfigure}[b]{0.49\textwidth}
	\includegraphics[width=\textwidth]{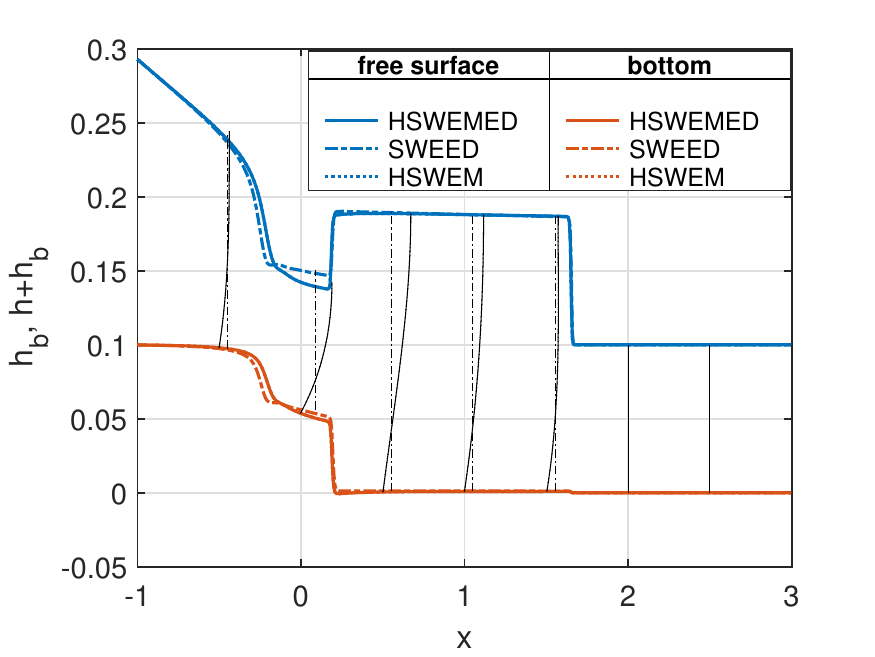}
	\caption*{(a)}
	\end{subfigure}
	\hfill
	\begin{subfigure}[b]{0.49\textwidth}
	\includegraphics[width=\textwidth]{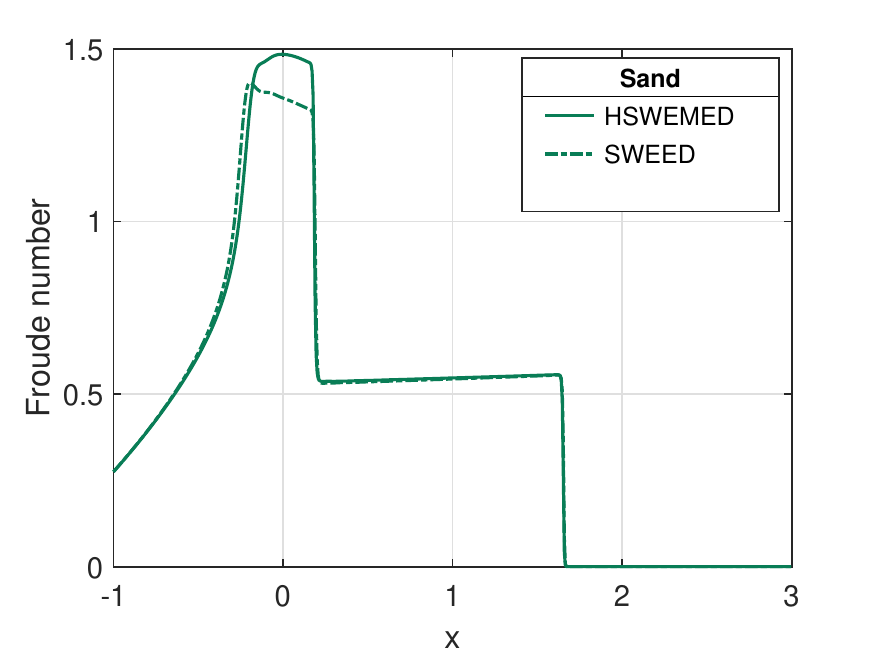}
	\caption*{(b)}
	\end{subfigure}
	\caption{For a sand bed with laboratory Config 2 \eqrefc{exp_2} at $t=1$ with the increased friction coefficient $\epsilon=0.0324$: (a) free surface, bottom evolution, and vertical velocity profiles at $x=-0.5,0,0.5,1,1.5,2,2.5$, computed with the HSWEMED (solid) and SWEED (dashed); and (b) Froude number, computed
		with the HSWEMED (solid) and the SWEED (dashed).
		}
	\label{exp3_diff_friction}
\end{figure}

For the sand bed, Figure~\ref{exp3_pvc}(b) shows a similar behavior. The hydraulic jump around $x=0$ is only partially captured, whereas the downstream free surface and the main sharp variations of the bottom are reproduced. We do not observe significant differences among the three models for the friction coefficient $\epsilon=0.0104$. This similarity is related to the small friction, which generates only a weak vertical velocity structure. When the friction coefficient is increased to $\epsilon=0.0324$ in Figure~\ref{exp3_diff_friction}(a), differences between the HSWEMED and SWEED become more visible near $x=0$, where the Froude number $\left(\dfrac{u_m}{\sqrt{gh}}\right)$ is relatively large, as shown in Figure~\ref{exp3_diff_friction}(b). The velocity profiles indicate that the vertical structure is most pronounced for $x\in[0,1]$, while the profile remains close to the depth-averaged velocity $u_m$ over the rest of the domain. A similar scenario was also analysed in \cite{Garres}.

Figures~\ref{exp3_pvc}(c) and~\ref{exp3_pvc}(d) illustrate that the peak of the sediment concentration and spatial expansion for the PVC bed are higher than those for the sand bed. The difference in concentration
computed by the HSWEMED and SWEED is related to the difference in density, settling velocity, and the low friction considered in this experiment for sand particles, as also explained in \secref{conf_1}.

%% file: 05-Conclusion.tex
\section{Conclusion}
We developed a one-dimensional coupled hydro-morphodynamic model to simulate sediment transport, incorporating bedload and suspended load dynamics with vertical velocity variations. The derived model is based on the shallow water moment framework, which governs the flow hydrodynamics, volumetric sediment concentration, and morphological evolution through the Exner equation. Using the moment approach allows us to recover the vertical structure of the horizontal velocity of the fluid, thereby approximating the velocity close to the bottom. The bottom velocity is important, as its approximation influences sediment discharge and erosion rates and, consequently, the bedload and suspended load transport. To close the system, we integrated the Meyer-Peter \& Müller formula for sediment discharge, empirical laws for erosion and deposition, non-linear Manning type friction for bed resistance, and Newtonian friction to represent fluid shear stresses. The reduced model is called the SWEMED and evolves the depth-averaged concentration $c_m$, while the near-bed concentration required in the deposition closure is prescribed algebraically from $c_m$.

\par For the hyperbolicity analysis, we applied the regularization technique of \cite{koellermeier2020analysis} to the moment block of SWEMED and obtained the HSWEMED. We derived the characteristic polynomial of the regularized transport matrix and obtained a real-root condition for the sediment-coupled cubic factor.
The analysis shows that the moment order affects the eigenvalue structure. For even moment orders, strict hyperbolicity can hold under the stated real-root and separation conditions. For odd moment orders, a repeated eigenvalue appears, so strict hyperbolicity is lost; under the same real-root condition, the regularized HSWEMED transport matrix is weakly hyperbolic.
\par The numerical results were obtained with a second-order path-conservative MUSCL-Rusanov finite-volume scheme. For the water-at-rest equilibrium, we introduced a well-balancing correction of the Rusanov viscosity. The corrected scheme preserves a constant free surface over both smooth and discontinuous bottoms, while the standard MUSCL-Rusanov scheme perturbs this equilibrium. The reference-convergence test gives rates consistent with the expected second-order behaviour of the scheme. The dam-break simulations then show the effect of the new moment-based erosion-deposition model. Compared with the standard depth-averaged SWEED, which uses only the mean velocity in the sediment transport closures, HSWEMED uses the moment-based bottom velocity. In the tested configurations, this gives a lower bottom velocity and therefore less bottom erosion, especially for the PVC bed. Compared with HSWEM, which includes the moment structure but neglects erosion, deposition, and suspended load, HSWEMED also accounts for the exchange of sediment between the water column and the erodible bed. The comparison with the laboratory experiments of \cite{Spinewine01122007} shows that HSWEMED gives a better description of the PVC bottom evolution in the tested case, while the differences between the models are smaller for sand and discrepancies remain near the hydraulic-jump region.

\par Ongoing work is directed toward equilibrium stability and steady-state analyses of the HSWEMED model. Future work should also consider improved numerical treatments for wet/dry fronts in higher-order moment models, including positivity-preserving and asymptotic-preserving schemes. Another possible extension is to introduce an independent vertical structure for the suspended concentration through additional concentration moments.

%% file: appendix_reference_system.tex
\section{Derivation of the complete reference system}\label{app:1}
\subsection{Mapping of the mass balance}
From the divergence-free conditions $\nabla \cdot \mathbf{u}=0$, we have
\begin{equation}\label{r1.1}
    \partial_x u + \partial_z w = 0.
\end{equation}
After multiplying \eqrefc{r1.1} by $h$ and applying the differential operator \eqrefc{diffop1}, the mapped mass balance yields
\begin{equation}
    \partial_x \left(h \Tilde{u}\right)-\partial_\zeta \left[\partial_x(\zeta h+h_b) \Tilde{u} \right] + \partial_\zeta \Tilde{w} = 0.
\end{equation}
\subsection{Mapping of the momentum balance}
From the reference system \eqrefc{reference system22},
\begin{equation}\label{ref:A2}
    \partial_t u + \partial_x u^2 +\partial_z (uw) = -\frac{1}{\rho} \partial_x p + \frac{1}{\rho} \partial_z \sigma_{xz}.
\end{equation}
After multiplying \eqrefc{ref:A2} by $h$ and applying the differential operator \eqrefc{diffop1}, the mapped momentum balance yields 
\begin{align}
    \begin{split}
            \partial_t (h \Tilde{u}) - \partial_\zeta \left(\Tilde{u} \partial_t\left(\zeta h+h_b\right)\right)& + \partial_x(h \Tilde{u}^2) -\partial_\zeta \left(\Tilde{u}^2 \partial_x\left(\zeta h+h_b\right)\right) + \partial_\zeta(\Tilde{u}\Tilde{w}) \\
            &+ \frac{1}{\rho} \partial_x (h \Tilde{p})- \frac{1}{\rho} \partial_\zeta \left(\Tilde{p}\partial_x(\zeta h+h_b)\right) = \frac{1}{\rho}\partial_\zeta \Tilde{\sigma}_{xz},
    \end{split}
    \\[2ex]
       \begin{split}\label{ref:A2_p}
          \Longrightarrow  \partial_t (h \Tilde{u})+ \partial_x(h \Tilde{u}^2) + & \partial_\zeta \left[\Tilde{u}\left( \underbrace{\Tilde{w} -\Tilde{u}\partial_x\left(\zeta h+h_b\right)- \partial_t\left(\zeta h+h_b\right)}_{\text{vertical coupling}}\right)\right] \\
            &+ \underbrace{ \frac{1}{\rho} \partial_x (h \Tilde{p})- \frac{1}{\rho} \partial_\zeta \left(\Tilde{p}\partial_x(\zeta h+h_b)\right)}_{\text{pressure}} = \frac{1}{\rho}\partial_\zeta \Tilde{\sigma}_{xz}.
    \end{split}
\end{align}
From the hydrostatic balance, we have $p(t,x,z)= (h_s-z)\rho g$. After mapping from $z$ to $\zeta$ we get
\begin{align}
    p(t,x,\zeta)  = (h_s-(\zeta h+ h_b)) \rho g
                  = h (1 -\zeta)\rho g.
\end{align}
Therefore, 
\begin{align}\label{pressure_der}
\frac{1}{\rho} \partial_x (h \Tilde{p})- \frac{1}{\rho} \partial_\zeta \left(\Tilde{p}\partial_x(\zeta h+h_b)\right)
&=  \frac{1}{\rho} \partial_x \left(\rho g h^2 (1-\zeta)\right) - \frac{1}{\rho} \partial_\zeta \left(\rho g h(1-\zeta) \partial_x(\zeta h+h_b)\right)\notag\\
& =  2g(1-\zeta) h \partial_x h+ \frac{gh^2}{\rho}(1-\zeta)\partial_x \rho - \left(gh\partial_x h - 2gh\zeta \partial_x h-gh \partial_x h_b\right)\notag\\
 & = gh\partial_x h + gh \partial_x h_b + \frac{gh^2}{\rho}(1-\zeta)\partial_x \rho.
\end{align}
Now we substitute \eqrefc{pressure_der} to \eqrefc{ref:A2_p}
\begin{equation}\label{ref:A2_subp}
          \partial_t (h \Tilde{u})+ \partial_x(h \Tilde{u}^2) +  \partial_\zeta \left[\Tilde{u}\left(\Tilde{w} -\Tilde{u}\partial_x\left(\zeta h+h_b\right)- \partial_t\left(\zeta h+h_b\right)\right)\right] 
            +gh\partial_x (h +  h_b) + \frac{gh^2}{\rho}(1-\zeta)\partial_x \rho = \frac{1}{\rho}\partial_\zeta \Tilde{\sigma}_{xz}.
\end{equation}

\subsection{Mapping of the sediment concentration equation}
From the reference system \eqrefc{reference system23}, we have
\begin{equation}\label{ref:A3}
    \partial_t c +\partial_x (cu) + \partial_z (cw) = 0.
\end{equation}
After multiplying \eqrefc{ref:A3} by $h$ and applying the differential operator \eqrefc{diffop1}, the mapped sediment concentration equation yields
\begin{align}
 \partial_t (h \Tilde{c})  + \partial_x \left(h \Tilde{c}\Tilde{u}\right) +\partial_\zeta \left[\Tilde{c} \left(\Tilde{w}- \partial_x(\zeta h+h_b)\Tilde{u}-\partial_t(\zeta h+h_b) \right) \right]=0.
\end{align}

%% file: appendix_depthaveraged_equation.tex
\section{Averages of the Horizontal Momentum Balances}\label{app:2}
\subsection{Depth-averaging the momentum balance}\label{averaged momentum balance}
We consider the momentum balance \eqrefc{eq:resolved  system1_2} and after integrating $\int_{0}^{1}\cdot \, d\zeta$, we get
\begin{align}\label{ref:A2_integration}
	\begin{split}
		\int_{0}^{1}\partial_t (h \Tilde{u})  \, d\zeta &+ \int_{0}^{1}\partial_x(h \Tilde{u}^2)  \, d\zeta  +  \int_{0}^{1} \partial_\zeta \left[\Tilde{u}\left(\Tilde{w} -\Tilde{u}\partial_x\left(\zeta h+h_b\right)- \partial_t\left(\zeta h+h_b\right)\right)\right] \, d\zeta  \\
		&+\int_{0}^{1}gh\partial_x h  \, d\zeta  + \int_{0}^{1}gh \partial_x h_b  \, d\zeta  + \int_{0}^{1}\frac{gh^2}{\rho}(1-\zeta)\partial_x \rho \,  \, d\zeta  = \frac{1}{\rho}\int_{0}^{1}\partial_\zeta \Tilde{\sigma}_{xz}  \, d\zeta, \notag
	\end{split}
	\\[2ex]
	\begin{split}
		\Rightarrow \partial_t \left(h \int_{0}^{1} \Tilde{u}   \, d\zeta \right) &+ \partial_x \left(h \int_{0}^{1} \Tilde{u}^2  \, d\zeta \right)   +  \partial_\zeta \left(  \int_{0}^{1}\Tilde{u}\left(\Tilde{w} -\Tilde{u}\partial_x\left(\zeta h+h_b\right)- \partial_t\left(\zeta h+h_b\right)\right) \, d\zeta\right)  \\
		&+gh\partial_x h   + gh \partial_x h_b  + \frac{gh^2}{2 \rho}\partial_x \rho  = \frac{1}{\rho}\partial_\zeta  \left(\int_{0}^{1} \Tilde{\sigma}_{xz}  \, d\zeta \right).
	\end{split}
\end{align}
After applying the kinematic boundary conditions \eqrefc{bc_zeta_1}-\eqrefc{bc_zeta_2} and substituting $\Tilde{u}(t,x,\zeta) = u_m(t,x)+ \sum_{j=1}^{N}\alpha_j(t,x)\,\phi_j(\zeta)$, we get

\begin{align}   
	\begin{split}
		\partial_t(h u_m) +  \partial_x \left( hu_m^2+ h\sum  \frac{\alpha_j^2}{2j+1}+  \frac{gh^2}{2} \right)
		=  - gh \partial_x h_b  - \frac{gh^2}{2 \rho}(\rho_s-\rho_w)\partial_x c_m +  \frac{\phi_b u_b}{1-\psi} - \epsilon |u_b|u_b.
	\end{split}
\end{align}

\textbf{Abbreviation of matrices}\label{Abbreviation}\\ \ \\
Here we introduce six matrices which are generated by the scaled Legendre polynomials
\begin{align*}
	A_{ijk} &= (2i+1)\int_{0}^{1} \phi_i \phi_j \phi_k \, d\zeta, 
	&\quad B_{ijk} &= (2i+1)\int_{0}^{1} \phi_i' \left(\phi_j d\hat{\zeta} \right) \phi_k \, d\zeta, \\
	C_{ij}   &= \int_{0}^{1} \phi_i' \phi_j' \, d\zeta, 
	&\quad H_{ij}   &= (2i+1)\int_{0}^{1} \zeta \phi_i \phi_j' \, d\zeta, \\
	G_{ij}   &= (2i+1)\int_{0}^{1} \phi_i \phi_j' \, d\zeta, 
	&\quad K_i      &= \int_{0}^{1} \zeta \phi_i \, d\zeta.
\end{align*}


\subsection{Higher-order averages}\label{Higher order averages}
To get the higher-order averages, multiply \eqrefc{eq:resolved  system1_2} with the ansatz function $\phi_i(\zeta)$ and integrate from $\zeta=0$ to $\zeta=1$, 
\begin{equation}
	\begin{split}
		&\partial_t \left(h \int_{0}^{1} \phi_i \Tilde{u}\, d\zeta\right)+ \partial_x \left(h \int_{0}^{1} \phi_i \Tilde{u}^2 \, d\zeta\right) + gh\partial_x(h+h_b)\int_{0}^{1}\phi_i \, d\zeta\\
		&+ \int_{0}^{1} \phi_i \partial_\zeta \left[\Tilde{u}\left(\Tilde{w} -\Tilde{u}\partial_x\left(\zeta h+h_b\right)- \partial_t\left(\zeta h+h_b\right)\right)\right] \, d\zeta 
		=- \frac{gh^2}{\rho}\partial_x \rho \int_{0}^{1}(1-\zeta) \phi_i \,  \, d\zeta + \frac{1}{\rho}\int_{0}^{1} \phi_i \partial_\zeta \Tilde{\sigma}_{xz}  \, d\zeta.   
	\end{split}
\end{equation}
We have, 
\begin{align}
	\int_{0}^{1} \phi_i \Tilde{u} \, d\zeta &=\int_{0}^{1} \phi_i \left(u_m+ \sum_{j=1}^{N}\alpha_j\phi_j\right) \, d\zeta
	= \frac{\alpha_i}{2i+1}.
\end{align}
and
\begin{align}
	\int_{0}^{1} \phi_i \Tilde{u}^2 \,  \, d\zeta &=\int_{0}^{1} \phi_i \left(u_m^2 + 2 u_m \sum_{j=1}^{N}\alpha_j\phi_j + \sum_{j=1}^{N}\alpha_j\phi_j \sum_{j=1}^{N}\alpha_j\phi_j  \right) \, d\zeta= \frac{1}{2i+1}\left(2u_m\alpha_i+\sum_{j,k=1}^{N}\alpha_j\alpha_kA_{ijk}\right).
\end{align}
The projection on vertical coupling term reads
\begin{align}
	&\int_{0}^{1} \phi_i \partial_\zeta \left[\Tilde{u}\left(\Tilde{w} -\Tilde{u}\partial_x\left(\zeta h+h_b\right)- \partial_t\left(\zeta h+h_b\right)\right)\right] \, d\zeta \notag\\
	= & \int_{0}^{1} \phi_i \partial_\zeta \left[\Tilde{u} \left(-\partial_x\left(h\int_{0}^{\zeta}\Tilde{u}d\hat{\zeta}\right)-\zeta\left(-\partial_x\left(hu_m\right)+\frac{F_b}{1-\psi}\right)-\partial_t\, h_b\right)\right] \, d\zeta \notag\\
	= & \int_{0}^{1} \phi_i \partial_\zeta \left[\Tilde{u} \left(\partial_x\left(h\int_{0}^{\zeta}(u_m-\Tilde{u})\, d\hat{\zeta}\right)\right) \right] \, d\zeta - \int_{0}^{1} \phi_i \partial_\zeta \left(\Tilde{u}\frac{\zeta F_b}{1-\psi}\right)\, d\zeta - \int_{0}^{1} \phi_i \partial_\zeta \left(\Tilde{u}\partial_t h_b\right) d\zeta \notag\\
	= & - \int_{0}^{1} \phi_i \partial_\zeta \left[\Tilde{u} \left(\partial_x\left(h\int_{0}^{\zeta} u_d \, d\hat{\zeta}\right)\right) \right] \, d\zeta - \frac{F_b}{1-\psi}\int_{0}^{1} \phi_i \partial_\zeta \left(\Tilde{u}\zeta\right) \, d\zeta - \partial_t h_b \int_{0}^{1} \phi_i \partial_\zeta \Tilde{u} \, d\zeta.
\end{align}

We notice the presence of the integral associated with erosion and deposition 
\begin{align}
	&\frac{F_b}{1-\psi}\int_{0}^{1} \phi_i \partial_\zeta \left(\Tilde{u}\zeta\right) \, d\zeta
	=  \frac{F_b}{1-\psi} \left( \sum_{j=1}^{N}\int_{0}^{1} \alpha_j \phi_i \phi_j \, d\zeta + \sum_{j=1}^{N}\int_{0}^{1} \alpha_j \zeta \phi_i \phi_j^{\prime}\, d\zeta \right)
	=  \frac{F_b}{1-\psi} \left( \frac{\alpha_i}{2i+1}+ \frac{\sum_{j=1}^{N} \alpha_j H_{ij}}{2i+1}\right).
\end{align}
and the projection related to the coupling with the bed evolution 
\begin{align}\label{A_s defination}
	\partial_t h_b \int_{0}^{1} \phi_i \partial_\zeta \Tilde{u} \, d\zeta  
	=  \partial_t h_b \sum_{j=1}^{N}\alpha_j \left(\int_{0}^{1}   \phi_i \phi_j^{\prime} \, d\zeta \right) 
	=  \partial_t h_b \frac{\sum_{j=1}^{N}\alpha_j G_{ij}}{2i+1}. 
\end{align}
Therefore, the projection of the vertical coupling for higher-order reads
\begin{align}
	\int_{0}^{1} \phi_i \partial_\zeta &\left[\Tilde{u}\left(\Tilde{w} -\Tilde{u}\partial_x\left(\zeta h+h_b\right)- \partial_t\left(\zeta h+h_b\right)\right)\right] \, d\zeta \notag
	=   - \frac{u_m}{2i+1}\partial_x (h \alpha_i) + \frac{1}{2i+1} \sum_{j,k=1}^{N} B_{ijk}\alpha_k \partial_x (h \alpha_j)\\
	&-\frac{1}{2i+1}\frac{F_b}{1-\psi}\left(\alpha_i+\sum_{j=1}^{N}\alpha_j (H_{ij}- G_{ij})\right)
	+\frac{1}{2i+1}  \sum_{j=1}^{N}\alpha_j G_{ij}\partial_x \frac{Q_b}{1-\psi}.
\end{align}
Similarly, we integrate the term that accounts for density variation
\begin{align}
	& \frac{gh^2}{\rho} \frac{\partial \rho}{\partial x} \int_{0}^{1} (1-\zeta) \phi_i  \, d\zeta 
	= \frac{gh^2}{\rho} (\rho_s-\rho_w)\frac{\partial c_m}{\partial x} \left(- \int_{0}^{1} \zeta \phi_i  \, d\zeta \right)
	= - \frac{gh^2}{\rho} (\rho_s-\rho_w)\frac{\partial c_m}{\partial x} K_i.
\end{align}
and the friction term for the additional moment equations,
\begin{align}
	\frac{1}{\rho}\int_{0}^{1} \phi_i \partial_\zeta \sigma_{xz} \, d\zeta= \frac{1}{\rho}\int_{0}^{1} \partial_\zeta (\phi_i \sigma_{xz})  \, d\zeta - \frac{1}{\rho} \int_{0}^{1} \sigma_{xz} \partial_\zeta \phi_i  \, d\zeta \notag
	&=  - \epsilon |u_b| u_b  - \frac{\nu}{h} \sum_{j=1}^{N} \int_{0}^{1} \alpha_j \phi_i^\prime \phi_j^\prime   \, d\zeta \notag\\
	&=- \epsilon |u_b| u_b - \frac{\nu}{h}\sum_{j=1}^{N} \alpha_j C_{ij}. 
\end{align}
After putting all together, we obtain the higher average momentum balance in $x$ direction 
\begin{align}
	\begin{split}
		& \partial_t (h\alpha_i)+ \partial_x \left(2hu_m\alpha_i+ h\sum_{j,k=1}^{N}\alpha_j\alpha_kA_{ijk}\right)= u_m\partial_x (h \alpha_i) -  \sum_{j,k=1}^{N} B_{ijk}\alpha_k \partial_x (h \alpha_j)\\
		& +\frac{F_b}{1-\psi}\left(\alpha_i+\sum_{j=1}^{N}\alpha_j (H_{ij}- G_{ij})\right)
		- \sum_{j=1}^{N}\alpha_j G_{ij}\partial_x \frac{Q_b}{1-\psi}+ (2i+1)\frac{gh^2}{\rho} (\rho_s-\rho_w)\frac{\partial c_m}{\partial x} K_i - (2i+1) \\
		&\left( \epsilon |u_b| u_b + \frac{\nu}{h}\sum_{j=1}^{N} \alpha_j C_{ij} \right).
	\end{split}
\end{align}

%% file: appendix_energy_equation.tex
\section{Energy Balance}\label{app:4}
In this section, we provide a detailed derivation of the energy balance for the
SWEMEDO (\theoremref{SWEMED0 energy balance}) and SWEMED1 (\theoremref{SWEMED1 energy balance}) model. We first
derive the potential-energy and kinetic-energy balances and then combine them
to obtain the total mechanical energy equation.

\subsection{SWEMEDO}\label{app:SWEMED0_energy}
We start from the SWEMEDO system
\begin{equation}
	\partial_t h + \partial_x(hu_m) = S_m,
	\qquad
	S_m = \frac{E-D}{1-\psi}.
	\tag{C}\label{eq:C}
\end{equation}

\begin{equation}
	\begin{aligned}
		\partial_t(hu_m)
		+ \partial_x\!\left(hu_m^2 + \frac{gh^2}{2}\right)
		&= -gh\,\partial_x h_b
		-\frac{gh^2}{2}\,\beta\,\partial_x c_m
		-\epsilon |u_m|u_m
		+ S_m u_m .
	\end{aligned}
	\tag{M}\label{eq:M}
\end{equation}

To facilitate the subsequent energy derivation, we rewrite the concentration-gradient term as
\begin{equation}
	\begin{aligned}
		h^2\,\partial_x c_m
		= h^2\,\partial_x c_m + hc_m\,\partial_x h - hc_m\,\partial_x h 
		= h\,\partial_x(hc_m) - hc_m\,\partial_x h .
	\end{aligned}
	\label{eq:h2dc_rewrite}
\end{equation}
Substituting \eqref{eq:h2dc_rewrite} into \eqref{eq:M} yields the equivalent momentum balance
\begin{equation}
	\begin{aligned}
		\partial_t(hu_m)
		+ \partial_x\!\left(hu_m^2 + \frac{gh^2}{2}\right)
		&= -gh\,\partial_x h_b
		-\frac{g\beta}{2}\Bigl(h\,\partial_x(hc_m) - hc_m\,\partial_x h\Bigr)
		-\epsilon |u_m|u_m
		+ S_m u_m .
	\end{aligned}
	\tag{M'}\label{eq:Mprime}
\end{equation}

The suspended-sediment and Exner equations are
\begin{equation}
	\partial_t(hc_m) + \partial_x(hc_m\,u_m) = E-D .
	\tag{S}\label{eq:S}
\end{equation}

\begin{equation}
	\partial_t h_b + \partial_x\!\left(\frac{Q_b}{1-\psi}\right) = -S_m .
	\tag{E}\label{eq:E}
\end{equation}

\subsection*{Potential energy}
Multiplying \eqref{eq:C} and \eqref{eq:E} by \(g(h+h_b)\) gives
\begin{equation}
	g(h+h_b)\,\partial_t h + g(h+h_b)\,\partial_x(hu_m) = g(h+h_b)\,S_m .
	\label{eq:PE_mult_C}
\end{equation}
\begin{equation}
	g(h+h_b)\,\partial_t h_b
	+ g(h+h_b)\,\partial_x\!\left(\frac{Q_b}{1-\psi}\right)
	= -g(h+h_b)\,S_m .
	\label{eq:PE_mult_E}
\end{equation}

Adding \eqref{eq:PE_mult_C} and \eqref{eq:PE_mult_E} yields
\begin{equation}
	g(h+h_b)\,\partial_t(h+h_b)
	+ g(h+h_b)\,\partial_x(hu_m)
	+ g(h+h_b)\,\partial_x\!\left(\frac{Q_b}{1-\psi}\right)
	= 0 .
	\label{eq:PE_a}
\end{equation}
Equivalently,
\begin{equation}
	\partial_t\!\left(\frac{g}{2}(h+h_b)^2\right)
	+ g(h+h_b)\,\partial_x(hu_m)
	+ g(h+h_b)\,\partial_x\!\left(\frac{Q_b}{1-\psi}\right)
	= 0 .
	\label{eq:PE_b}
\end{equation}

\subsection*{Kinetic energy}
We rewrite the momentum balance \eqref{eq:Mprime} in the form
\begin{equation}
	\begin{aligned}
		\partial_t(hu_m)
		+ \partial_x(hu_m^2)
		+ gh\,\partial_x(h+h_b)
		&= -\frac{g\beta}{2}\Bigl(h\,\partial_x(hc_m)-hc_m\,\partial_x h\Bigr)
		-\epsilon|u_m|u_m
		+ S_m u_m .
	\end{aligned}
	\tag{M1}\label{eq:M1}
\end{equation}

Next, we compute \((M1)-u_m(C)\). Using \eqref{eq:C} and \eqref{eq:M1}, we obtain
\begin{equation}
	\begin{aligned}
		\partial_t(hu_m) - u_m\,\partial_t h
		+ \partial_x(hu_m^2) - u_m\,\partial_x(hu_m)
		+ gh\,\partial_x(h+h_b)
		&= -\frac{g\beta}{2}\Bigl(h\,\partial_x(hc_m)-hc_m\,\partial_x h\Bigr)
		-\epsilon|u_m|u_m .
	\end{aligned}
	\label{eq:M1_minus_umC}
\end{equation}
Equivalently,
\begin{equation}
	\begin{aligned}
		h\,\partial_t u_m + h u_m\,\partial_x u_m + gh\,\partial_x(h+h_b)
		&= -\frac{g\beta}{2}\Bigl(h\,\partial_x(hc_m)-hc_m\,\partial_x h\Bigr)
		-\epsilon|u_m|u_m .
	\end{aligned}
	\tag{A}\label{eq:A}
\end{equation}

Averaging \eqref{eq:A} and \eqref{eq:M1} yields
\begin{equation}
	\begin{aligned}
		\frac{1}{2}\Bigl[\partial_t(hu_m)+h\,\partial_t u_m\Bigr]
		+ \frac{1}{2}\Bigl[\partial_x(hu_m^2)+h u_m\,\partial_x u_m\Bigr]
		+ gh\,\partial_x(h+h_b)
		= -\frac{g\beta}{2}\Bigl(h\,\partial_x(hc_m)-hc_m\,\partial_x h\Bigr)\\
		-\epsilon|u_m|u_m
		+ \frac{1}{2}S_m u_m .
	\end{aligned}
	\tag{S'}\label{eq:Sprime}
\end{equation}

Multiplying \eqref{eq:Sprime} by \(u_m\) gives the kinetic-energy balance
\begin{equation}
	\begin{aligned}
		\partial_t\!\left(\frac{h u_m^2}{2}\right)
		+ \partial_x\!\left(\frac{h u_m^3}{2}\right)
		+ gh\,u_m\,\partial_x(h+h_b)
		&= -\frac{g\beta}{2}\,u_m\Bigl(h\,\partial_x(hc_m)-hc_m\,\partial_x h\Bigr)
		-\epsilon u_m^2|u_m|
		+ \frac{1}{2}S_m u_m^2 .
	\end{aligned}
	\tag{K}\label{eq:K}
\end{equation}

We rewrite the coupling term as
\begin{equation}
	\begin{aligned}
		-\frac{g\beta}{2}\Bigl(u_m h\,\partial_x(hc_m) - hc_m\,u_m\,\partial_x h\Bigr)
		&= -\frac{g\beta}{2}\,h^2 u_m\,\partial_x c_m \\
		&= -\frac{g\beta}{2}\Bigl(\partial_x(h^2 c_m u_m) - c_m\,\partial_x(h^2 u_m)\Bigr),
	\end{aligned}
	\label{eq:KE_term_rewrite}
\end{equation}
since
\begin{equation}
	\partial_x(h^2 c_m u_m) = h^2 u_m\,\partial_x c_m + c_m\,\partial_x(h^2 u_m).
	\label{eq:product_rule_note}
\end{equation}

Substituting \eqref{eq:KE_term_rewrite} into \eqref{eq:K} yields
\begin{equation}
	\begin{aligned}
		\partial_t\!\left(\frac{h u_m^2}{2}\right)
		+ \partial_x\!\left(\frac{h u_m^3}{2}\right)
		+ gh\,u_m\,\partial_x(h+h_b)
		&= -\frac{g\beta}{2}\,\partial_x(h^2 c_m u_m)
		+ \frac{g\beta}{2}\,c_m\,\partial_x(h^2 u_m)
		-\epsilon |u_m|u_m^2
		+ \frac{1}{2}S_m u_m^2 .
	\end{aligned}
	\tag{K1}\label{eq:K1}
\end{equation}

Equivalently, moving the divergence term to the left gives
\begin{equation}
	\begin{aligned}
		\partial_t\!\left(\frac{h u_m^2}{2}\right)
		+ \partial_x\!\left(u_m\left(\frac{h u_m^2}{2}+\frac{g\beta}{2}h^2 c_m\right)\right)
		+ gh\,u_m\,\partial_x(h+h_b)
		&= \frac{g\beta}{2}\,c_m\,\partial_x(h^2 u_m)
		-\epsilon |u_m|u_m^2
		+ \frac{1}{2}S_m u_m^2 .
	\end{aligned}
	\label{eq:K1_fluxform}
\end{equation}

\subsection*{Total mechanical energy}
We define the total mechanical energy as
\begin{equation}
	\mathcal{E} := \frac{h u_m^2}{2} + \frac{g}{2}(h+h_b)^2 .
	\label{eq:total_energy_def}
\end{equation}

Adding the mechanical energy \eqref{eq:K1_fluxform} and the potential-energy balance \eqref{eq:PE_b} yields
\begin{equation}
	\begin{aligned}
		\partial_t \mathcal{E}
		+ \partial_x\!\left(u_m\left(\frac{h u_m^2}{2}+\frac{g\beta}{2}h^2 c_m\right)\right)
		+ gh\,u_m\,\partial_x(h+h_b)
		+ g(h+h_b)\,\partial_x(hu_m)
		+ g(h+h_b)\,\partial_x\!\left(\frac{Q_b}{1-\psi}\right)
		&\\= -\epsilon |u_m|u_m^2
		+ \frac{1}{2}S_m u_m^2
		+ \frac{g\beta}{2}\,c_m\,\partial_x(h^2 u_m).
	\end{aligned}
	\label{eq:total_energy_intermediate}
\end{equation}

Using the product rules
\[
gh\,u_m\,\partial_x(h+h_b) + g(h+h_b)\,\partial_x(hu_m)
= \partial_x\!\bigl(ghu_m(h+h_b)\bigr),
\]
and
\[
g(h+h_b)\,\partial_x\!\left(\frac{Q_b}{1-\psi}\right)
= \partial_x\!\left(g(h+h_b)\frac{Q_b}{1-\psi}\right)
- g\frac{Q_b}{1-\psi}\,\partial_x(h+h_b),
\]
we obtain the energy equation
\begin{equation}
	\begin{aligned}
		\partial_t \mathcal{E}
		+ \partial_x\!\Biggl(
		u_m\left(\frac{h u_m^2}{2} + g h^2 + \frac{g\beta}{2}h^2 c_m\right)
		+ g(h+h_b)\frac{Q_b}{1-\psi}
		+ gh\,u_m\,h_b
		\Biggr)
		= -\epsilon |u_m|u_m^2
		+ \frac{1}{2}S_m u_m^2
		\\ +\,g\frac{Q_b}{1-\psi}\,\partial_x(h+h_b)
		+ \frac{g\beta}{2}\,c_m\,\partial_x(h^2 u_m).
	\end{aligned}
	\label{eq:total_energy_final}
\end{equation}
Under the natural assumptions $\epsilon\geq 0$, $S_m\leq 0$, and $Q_b\partial_x(h+h_b)\leq 0$, the bracketed dissipation in \eqref{eq:total_energy_final} has the correct sign. Indeed, the Meyer-Peter and M\"uller flux satisfies $Q_b=\operatorname{sgn}(u_b)Q\Phi(\theta)$ with $\Phi(\theta)\geq0$, and the downslope condition $\operatorname{sgn}(u_b)\partial_x(h+h_b)\leq0$ gives
\begin{equation}
	\epsilon |u_m|u_m^2\geq0,
	\qquad
	-\frac12 S_m u_m^2\geq0,
	\qquad
	-g\frac{Q_b}{1-\psi}\partial_x(h+h_b)\geq0.
\end{equation}
The remaining density-contrast term can be written as a flux contribution and a dissipative remainder. Using $0\leq c_m\leq1$,
\begin{equation}
\begin{aligned}
\frac{g\beta}{2}c_m\partial_x(h^2u_m)
&=
\partial_x\left(\frac{g\beta}{2}\frac{1+\operatorname{sgn}(\partial_x(h^2u_m))}{2}h^2u_m\right)\\
&\quad
-\frac{g\beta}{2}\left(\frac{1+\operatorname{sgn}(\partial_x(h^2u_m))}{2}-c_m\right)\partial_x(h^2u_m).
\end{aligned}
\end{equation}
Thus, the first term is absorbed into the energy flux, while the second one contributes to the dissipation.

\subsection{SWEMED1}\label{app:SWEMED1_energy}

We consider the SWEMED1 model (mean mode $u_m$ and first moment $\alpha_1$).
Starting from
\begin{equation}
	\partial_t h + \partial_x(hu_m) = S_m,
	\qquad
	S_m = \frac{E-D}{1-\psi},
	\tag{C}\label{eq:SWEMED1_C}
\end{equation}
and the momentum balance
\begin{equation}
	\begin{aligned}
		\partial_t(hu_m)
		+ \partial_x\!\left(hu_m^2 + \frac{1}{3}h\alpha_1^2\right)
		+ gh\,\partial_x(h+h_b)
		&= -\frac{gh^2}{2}\,\beta\,\partial_x c_m
		-\epsilon\,|u_m+\alpha_1|\,(u_m+\alpha_1)
		+ S_m(u_m+\alpha_1).
	\end{aligned}
	\tag{M}\label{eq:SWEMED1_M}
\end{equation}

As before, we rewrite the concentration-gradient term as
\begin{equation}
	h^2\,\partial_x c_m
	= h\,\partial_x(hc_m) - hc_m\,\partial_x h,
	\label{eq:SWEMED1_h2dc}
\end{equation}
so that \eqref{eq:SWEMED1_M} becomes
\begin{equation}
	\begin{aligned}
		\partial_t(hu_m)
		+ \partial_x\!\left(hu_m^2 + \frac{1}{3}h\alpha_1^2\right)
		+ gh\,\partial_x(h+h_b)
		= -\frac{g\beta}{2}\Bigl(h\,\partial_x(hc_m)-hc_m\,\partial_x h\Bigr)\\
		-\epsilon\,|u_m+\alpha_1|\,(u_m+\alpha_1)
		+ S_m(u_m+\alpha_1).
	\end{aligned}
	\tag{M1}\label{eq:SWEMED1_M1}
\end{equation}

\subsubsection*{Mean kinetic-energy balance}

Compute $(M1)-u_m(C)$. Using
\(
\partial_t(hu_m)-u_m\partial_t h = h\,\partial_t u_m
\)
and
\(
\partial_x(hu_m^2)-u_m\partial_x(hu_m)=h u_m\,\partial_x u_m,
\)
we obtain
\begin{equation}
	\begin{aligned}
		h\,\partial_t u_m
		+ h u_m\,\partial_x u_m
		+ gh\,\partial_x(h+h_b)
		+ \partial_x\!\left(\frac{1}{3}h\alpha_1^2\right)
		&= -\frac{g\beta}{2}\Bigl(h\,\partial_x(hc_m)-hc_m\,\partial_x h\Bigr)
		-\epsilon\,|u_m+\alpha_1|\,(u_m+\alpha_1)
		+ S_m\alpha_1 .
	\end{aligned}
	\tag{A}\label{eq:SWEMED1_A}
\end{equation}

Averaging $(A)$ and $(M1)$ yields $(S)=\tfrac12(A)+\tfrac12(M1)$:
\begin{equation}
	\begin{aligned}
		\frac12\Bigl[h\,\partial_t u_m + \partial_t(hu_m)\Bigr]
		+\frac12\Bigl[h u_m\,\partial_x u_m + \partial_x(hu_m^2)\Bigr]
		+ gh\,\partial_x(h+h_b)
		+ \partial_x\!\left(\frac13 h\alpha_1^2\right)
		= -\frac{g\beta}{2}\Bigl(h\,\partial_x(hc_m)-hc_m\,\partial_x h\Bigr) \\
		 -\epsilon\,|u_m+\alpha_1|\,(u_m+\alpha_1)
		+ \frac12 S_m u_m + S_m\alpha_1 .
	\end{aligned}
	\tag{S}\label{eq:SWEMED1_S}
\end{equation}

Multiplying \eqref{eq:SWEMED1_S} by $u_m$ gives the kinetic-energy balance for $u_m$:
\begin{equation}
	\begin{aligned}
		\partial_t\!\left(\frac12 h u_m^2\right)
		+ \partial_x\!\left(\frac12 h u_m^3\right)
		+ gh\,u_m\,\partial_x(h+h_b)
		+ u_m\,\partial_x\!\left(\frac13 h\alpha_1^2\right)
		= -\frac{g\beta}{2}u_m\Bigl(h\,\partial_x(hc_m)-hc_m\,\partial_x h\Bigr) \\
	 -\epsilon\,u_m\,|u_m+\alpha_1|\,(u_m+\alpha_1)
		+ S_m\left(\frac12 u_m^2 + u_m\alpha_1\right).
	\end{aligned}
	\tag{\(K_u\)}\label{eq:SWEMED1_Ku}
\end{equation}

Rewrite the coupling term as
\begin{equation}
	\begin{aligned}
		-\frac{g\beta}{2}u_m\Bigl(h\,\partial_x(hc_m)-hc_m\,\partial_x h\Bigr)
		&= -\frac{g\beta}{2}h^2 u_m\,\partial_x c_m \\
		&= -\frac{g\beta}{2}\Bigl(\partial_x(h^2 c_m u_m) - c_m\,\partial_x(h^2 u_m)\Bigr),
	\end{aligned}
	\label{eq:SWEMED1_coupling_u}
\end{equation}
so \eqref{eq:SWEMED1_Ku} becomes
\begin{equation}
	\begin{aligned}
		\partial_t\!\left(\frac12 h u_m^2\right)
		+ \partial_x\!\left(\frac12 h u_m^3 + \frac{g\beta}{2}h^2 c_m\,u_m\right)
		+ gh\,u_m\,\partial_x(h+h_b)
		+ u_m\,\partial_x\!\left(\frac13 h\alpha_1^2\right)
		= \frac{g\beta}{2}c_m\,\partial_x(h^2 u_m) \\
	 -\epsilon\,u_m\,|u_m+\alpha_1|\,(u_m+\alpha_1)
		+ S_m\left(\frac12 u_m^2 + u_m\alpha_1\right).
	\end{aligned}
	\tag{\(K_u'\)}\label{eq:SWEMED1_Ku_flux}
\end{equation}

\subsubsection*{Moment equation ($N=1$) and weighted moment-energy balance}

Start from the $N=1$ moment equation
\begin{equation}
	\begin{aligned}
		\partial_t(h\alpha_1) + \partial_x(2hu_m\alpha_1)
		&= u_m\,\partial_x(h\alpha_1)
		-\frac{gh^2}{2}\beta\,\partial_x c_m
		-3\epsilon\,|u_m+\alpha_1|\,(u_m+\alpha_1)
		-\frac{12\nu}{h}\alpha_1
		+2\alpha_1 S_m ,
	\end{aligned}
	\label{eq:MN1_start}
\end{equation}
which can be rearranged as
\begin{equation}
	\begin{aligned}
		\partial_t(h\alpha_1) + \partial_x(hu_m\alpha_1) + h\alpha_1\,\partial_x u_m
		&= -\frac{g\beta}{2}\Bigl(h\,\partial_x(hc_m)-hc_m\,\partial_x h\Bigr)
		-3\epsilon\,|u_m+\alpha_1|\,(u_m+\alpha_1)
		-\frac{12\nu}{h}\alpha_1
		+2\alpha_1 S_m .
	\end{aligned}
	\tag{MN1}\label{eq:MN1}
\end{equation}

Compute $(MN1)-\alpha_1(C)$:
\begin{equation}
	\begin{aligned}
		h\,\partial_t\alpha_1 + h u_m\,\partial_x\alpha_1 + h\alpha_1\,\partial_x u_m
		&= -\frac{g\beta}{2}\Bigl(h\,\partial_x(hc_m)-hc_m\,\partial_x h\Bigr)
		-3\epsilon\,|u_m+\alpha_1|\,(u_m+\alpha_1)
		-\frac{12\nu}{h}\alpha_1
		+\alpha_1 S_m .
	\end{aligned}
	\tag{AMN1}\label{eq:AMN1}
\end{equation}

Averaging $\tfrac12(AMN1)+\tfrac12(MN1)$ yields
\begin{equation}
	\begin{aligned}
		\frac12\Bigl[\partial_t(h\alpha_1)+h\,\partial_t\alpha_1\Bigr]
		+\frac12\Bigl[\partial_x(hu_m\alpha_1)+h u_m\,\partial_x\alpha_1\Bigr]
		+h\alpha_1\,\partial_x u_m
		= -\frac{g\beta}{2}\Bigl(h\,\partial_x(hc_m)-hc_m\,\partial_x h\Bigr) \\
		 -3\epsilon\,|u_m+\alpha_1|\,(u_m+\alpha_1)
		-\frac{12\nu}{h}\alpha_1
		+\frac32\,\alpha_1 S_m .
	\end{aligned}
	\tag{SMN1}\label{eq:SMN1}
\end{equation}

Following the standard ``partial kinetic energy'' weighting for $i=1$ (division by $2i+1=3$),
we multiply \eqref{eq:SMN1} by $\alpha_1/3$ to obtain
\begin{equation}
	\begin{aligned}
		\partial_t\!\left(\frac{1}{6}h\alpha_1^2\right)
		+ \partial_x\!\left(\frac{1}{6}h u_m\alpha_1^2\right)
		+ \frac{1}{3}h\alpha_1^2\,\partial_x u_m
		&= -\frac{g\beta}{2}\frac{\alpha_1}{3}\Bigl(h\,\partial_x(hc_m)-hc_m\,\partial_x h\Bigr) \\
		&\quad -\alpha_1\epsilon\,|u_m+\alpha_1|\,(u_m+\alpha_1)
		-\frac{4\nu}{h}\alpha_1^2
		+\frac12\,\alpha_1^2 S_m .
	\end{aligned}
	\tag{\(K_{\alpha,1}\)}\label{eq:Kalpha_weighted_preflux}
\end{equation}

Rewrite the coupling term in \eqref{eq:Kalpha_weighted_preflux} as
\begin{equation}
	\begin{aligned}
		-\frac{g\beta}{2}\frac{\alpha_1}{3}\Bigl(h\,\partial_x(hc_m)-hc_m\,\partial_x h\Bigr)
		&= -\frac{g\beta}{2}h^2\frac{\alpha_1}{3}\,\partial_x c_m \\
		&= -\frac{g\beta}{2}\Bigl(\partial_x(h^2 c_m\,\tfrac{\alpha_1}{3})
		- c_m\,\partial_x(h^2\tfrac{\alpha_1}{3})\Bigr),
	\end{aligned}
	\label{eq:SWEMED1_coupling_a1}
\end{equation}
and move the divergence term to the left:
\begin{equation}
	\begin{aligned}
		\partial_t\!\left(\frac{1}{6}h\alpha_1^2\right)
		+ \partial_x\!\left(\frac{1}{6}h u_m\alpha_1^2
		+ \frac{g\beta}{2}h^2 c_m\,\frac{\alpha_1}{3}\right)
		+ \frac{1}{3}h\alpha_1^2\,\partial_x u_m
		= \frac{g\beta}{2}c_m\,\partial_x\!\left(h^2\frac{\alpha_1}{3}\right) \\
	 -\alpha_1\epsilon\,|u_m+\alpha_1|\,(u_m+\alpha_1)
		-\frac{4\nu}{h}\alpha_1^2
		+\frac12\,\alpha_1^2 S_m .
	\end{aligned}
	\tag{\(K_{\alpha,1}'\)}\label{eq:Kalpha_weighted}
\end{equation}

\subsubsection*{Total kinetic energy (mean + first moment)}

Define the total kinetic energy (mean + first moment) as
\begin{equation}
	\mathcal{K}_1
	:= \frac12 h u_m^2 + \frac{1}{6}h\alpha_1^2 .
	\label{eq:K1_def}
\end{equation}

Adding \eqref{eq:SWEMED1_Ku_flux} and \eqref{eq:Kalpha_weighted} gives
\begin{equation}
	\begin{aligned}
		\partial_t \mathcal{K}_1
		+ \partial_x\!\left(
		\frac12 h u_m^3
		+ \frac{g\beta}{2}h^2 c_m\Bigl(u_m+\frac{\alpha_1}{3}\Bigr)
		+ \frac{1}{6}h u_m\alpha_1^2
		\right)
		+ gh\,u_m\,\partial_x(h+h_b)
		+ u_m\,\partial_x\!\left(\frac13 h\alpha_1^2\right)
		+ \frac13 h\alpha_1^2\,\partial_x u_m\\
		= \frac{g\beta}{2}c_m\,\partial_x\!\left(h^2\Bigl(u_m+\frac{\alpha_1}{3}\Bigr)\right)
		-\epsilon\,|u_m+\alpha_1|\,(u_m+\alpha_1)^2
		-\frac{4\nu}{h}\alpha_1^2
		+ S_m\left(\frac12 u_m^2 + u_m\alpha_1 + \frac12\alpha_1^2\right).
	\end{aligned}
	\label{eq:K1_intermediate}
\end{equation}

Finally, use the identity
\begin{equation}
	u_m\,\partial_x\!\left(\frac13 h\alpha_1^2\right)
	+ \frac13 h\alpha_1^2\,\partial_x u_m
	+ \partial_x\!\left(\frac16 h u_m\alpha_1^2\right)
	= \partial_x\!\left(\frac12 h u_m\alpha_1^2\right),
	\label{eq:triad_identity_explicit}
\end{equation}
to rewrite \eqref{eq:K1_intermediate} in the conservative form
\begin{equation}
	\begin{aligned}
		\partial_t \mathcal{K}_1
		+ \partial_x\!\left(
		\frac12 h u_m^3
		+ \frac12 h u_m\alpha_1^2
		+ \frac{g\beta}{2}h^2 c_m\Bigl(u_m+\frac{\alpha_1}{3}\Bigr)
		\right)
		+ gh\,u_m\,\partial_x(h+h_b)
		= \frac{g\beta}{2}c_m\,\partial_x\!\left(h^2\Bigl(u_m+\frac{\alpha_1}{3}\Bigr)\right) \\
		 -\epsilon\,|u_m+\alpha_1|\,(u_m+\alpha_1)^2
		-\frac{4\nu}{h}\alpha_1^2
		+ S_m\left(\frac12 u_m^2 + u_m\alpha_1 + \frac12\alpha_1^2\right).
	\end{aligned}
	\label{eq:K1_total}
\end{equation}
\subsection*{Potential energy}

Introduce the Exner equation
\begin{equation}
	\partial_t h_b + \partial_x\!\left(\frac{Q_b}{1-\psi}\right) = -S_m .
	\tag{E}\label{eq:SWEMED1_E}
\end{equation}
 Multiply \eqref{eq:SWEMED1_C} and \eqref{eq:SWEMED1_E} by $g(h+h_b)$:
\begin{subequations}\label{eq:SWEMED1_PE_mult}
	\begin{align}
		g(h+h_b)\,(C):\qquad
		g(h+h_b)\,\partial_t h + g(h+h_b)\,\partial_x(hu_m) &= g(h+h_b)\,S_m,
		\label{eq:SWEMED1_PE_mult_C}\\
		g(h+h_b)\,(E):\qquad
		g(h+h_b)\,\partial_t h_b + g(h+h_b)\,\partial_x\!\left(\frac{Q_b}{1-\psi}\right) &= -g(h+h_b)\,S_m.
		\label{eq:SWEMED1_PE_mult_E}
	\end{align}
\end{subequations}
Adding \eqref{eq:SWEMED1_PE_mult_C} and \eqref{eq:SWEMED1_PE_mult_E} yields
\begin{subequations}\label{eq:SWEMED1_PE}
	\begin{align}
		g(h+h_b)\,\partial_t (h+h_b)
		+ g(h+h_b)\,\partial_x(hu_m)
		+ g(h+h_b)\,\partial_x\!\left(\frac{Q_b}{1-\psi}\right)
		&= 0,
		\label{eq:SWEMED1_PE_a}\\
		\Rightarrow\quad
		\partial_t\!\left(\frac{g}{2}(h+h_b)^2\right)
		+ g(h+h_b)\,\partial_x(hu_m)
		+ g(h+h_b)\,\partial_x\!\left(\frac{Q_b}{1-\psi}\right)
		&= 0.
		\label{eq:SWEMED1_PE_b}
	\end{align}
\end{subequations}

\subsection*{Total mechanical energy (kinetic + potential)}

Define the weighted total kinetic energy (mean $+$ first moment) by
\begin{equation}
	\mathcal{K}_1 := \frac12 h u_m^2 + \frac16 h\alpha_1^2,
	\label{eq:SWEMED1_K1_def}
\end{equation}
and the total mechanical energy by
\begin{equation}
	\mathcal{E}_1 := \mathcal{K}_1 + \frac{g}{2}(h+h_b)^2,
	\label{eq:SWEMED1_E1_def}
\end{equation}

Assume the total kinetic-energy balance has been obtained in the form
\begin{equation}
	\begin{aligned}
		\partial_t \mathcal{K}_1
		+ \partial_x\!\left(
		\frac12 h u_m^3
		+ \frac12 h u_m\alpha_1^2
		+ \frac{g\beta}{2}h^2 c_m\Bigl(u_m+\frac{\alpha_1}{3}\Bigr)
		\right)
		+ gh\,u_m\,\partial_x (h+h_b)
		= \frac{g\beta}{2}c_m\,\partial_x\!\left(h^2\Bigl(u_m+\frac{\alpha_1}{3}\Bigr)\right) \\
		 -\epsilon\,|u_m+\alpha_1|\,(u_m+\alpha_1)^2
		-\frac{4\nu}{h}\alpha_1^2
		 + S_m\left(\frac12 u_m^2 + u_m\alpha_1 + \frac12\alpha_1^2\right).
	\end{aligned}
	\label{eq:SWEMED1_K1_total}
\end{equation}

Adding \eqref{eq:SWEMED1_PE_b} to \eqref{eq:SWEMED1_K1_total} gives
\begin{equation}
	\begin{aligned}
		\partial_t \mathcal{E}_1
		&+ \partial_x\!\left(
		\frac12 h u_m^3
		+ \frac12 h u_m\alpha_1^2
		+ \frac{g\beta}{2}h^2 c_m\Bigl(u_m+\frac{\alpha_1}{3}\Bigr)
		\right)
		+ gh\,u_m\,\partial_x (h+h_b)
		+ g(h+h_b)\,\partial_x(hu_m)
		+ g(h+h_b)\,\partial_x\!\left(\frac{Q_b}{1-\psi}\right) \\
		&= \frac{g\beta}{2}c_m\,\partial_x\!\left(h^2\Bigl(u_m+\frac{\alpha_1}{3}\Bigr)\right)
		-\epsilon\,|u_m+\alpha_1|\,(u_m+\alpha_1)^2
		-\frac{4\nu}{h}\alpha_1^2
		+ S_m\left(\frac12 u_m^2 + u_m\alpha_1 + \frac12\alpha_1^2\right).
	\end{aligned}
	\label{eq:SWEMED1_E1_intermediate}
\end{equation}

Using the product rules
\begin{equation}
	gh\,u_m\,\partial_x (h+h_b) + g(h+h_b)\,\partial_x(hu_m)
	= \partial_x\!\bigl(ghu_m(h+h_b)\bigr),
	\label{eq:SWEMED1_prod1}
\end{equation}
and
\begin{equation}
	g(h+h_b)\,\partial_x\!\left(\frac{Q_b}{1-\psi}\right)
	= \partial_x\!\left(g(h+h_b)\frac{Q_b}{1-\psi}\right)
	- g\frac{Q_b}{1-\psi}\,\partial_x (h+h_b),
	\label{eq:SWEMED1_prod2}
\end{equation}
we obtain the conservative mechanical energy balance
\begin{equation}
	\begin{aligned}
		\partial_t \mathcal{E}_1
		+ \partial_x\!\Biggl(
		\frac12 h u_m^3
		+ \frac12 h u_m\alpha_1^2
		+ \frac{g\beta}{2}h^2 c_m\Bigl(u_m+\frac{\alpha_1}{3}\Bigr)
		+ ghu_m(h+h_b)
		+ g(h+h_b)\frac{Q_b}{1-\psi}
		\Biggr)
		\\ =\frac{g\beta}{2}c_m\,\partial_x\!\left(h^2\Bigl(u_m+\frac{\alpha_1}{3}\Bigr)\right) 
		 -\epsilon\,|u_b|\,u_b^2
		-\frac{4\nu}{h}\alpha_1^2 
		 + \frac12 S_m u_b^2
		+ g\frac{Q_b}{1-\psi}\,\partial_x (h+h_b).
	\end{aligned}
	\label{eq:SWEMED1_E1_final}
\end{equation}

Under the same assumptions $\epsilon\geq0$, $S_m\leq0$, and $Q_b\partial_x(h+h_b)\leq0$, the terms
\begin{equation}
	\epsilon |u_b|u_b^2,
	\qquad
	-\frac12 S_m u_b^2,
	\qquad
	-g\frac{Q_b}{1-\psi}\partial_x(h+h_b),
	\qquad
	\frac{4\nu}{h}\alpha_1^2
\end{equation}
are non-negative in the dissipation. The density-contrast term is treated as in the SWEMED0 case, now with $v=u_m+\alpha_1/3$:
\begin{equation}
\begin{aligned}
\frac{g\beta}{2}c_m\partial_x(h^2v)
&=
\partial_x\left(\frac{g\beta}{2}\frac{1+\operatorname{sgn}(\partial_x(h^2v))}{2}h^2v\right)\\
&\quad
-\frac{g\beta}{2}\left(\frac{1+\operatorname{sgn}(\partial_x(h^2v))}{2}-c_m\right)\partial_x(h^2v).
\end{aligned}
\end{equation}
This gives a flux term and a dissipative remainder. The additional viscous term $-4\nu\alpha_1^2/h$ is strictly dissipative when $\alpha_1\neq0$.

%% file: Proof_of_characteristic_polynomial.tex
\section{Proof of HSWEMED characteristic polynomial}
\begin{theorem}
	The transport matrix $\tilde{A}_{ED} \in \mathbb{R}^{(N+4)\times(N+4)}$ of HSWEMED has the following characteristic polynomial
	\begin{equation*}
			\chi_{\tilde{A}_{ED}}(\lambda)= (u_m-\lambda) \left[(-\lambda) \left( (\lambda-u_m)^2-gh-\alpha_1^2\right)+gh(\delta_h+\lambda\delta_q + c_m\delta_c+ 2\alpha_1\delta_q)\right] \cdot \chi_{A_2}(\lambda-u_m)
		\end{equation*} 
		where $A_{2} \in \mathbb{R}^{N \times N}$ is defined as follows
		\begin{equation*}
			\hspace{-2mm}
			A_{2} = \begin{pmatrix} 
				& c_2 &  &    \\
				a_2 &  & \ddots &     \\
				& \ddots &   & c_N     \\
				&    & a_N &   
			\end{pmatrix},
		\end{equation*}
		with values $c_i = \frac{i+1}{2i+1}\alpha_1$ and $a_i = \frac{i-1}{2i-1}\alpha_1$ the values above and below the diagonal, respectively, from \eqrefc{Aed transport matrix}.
	\end{theorem}
	Proof: The proof follows the proof of the characteristic polynomial of the HSWEM
	system matrix in \cite{Garres}, which is extended to include the sediment concentration and erosion-deposition effects at the bottom to explain the suspended and bedload transport.
	\par The transport matrix of HSWEMED for any arbitrary order $N$ is given by
	\begin{equation*}
		\hspace{-5mm}
		\tilde{A}_{ED}= \begin{pmatrix} 
			& 1 &  &  &  &  & &  \\
			gh-u_m^2-\frac{\alpha_1^2}{3}-\frac{gh(\rho-\rho_w)}{2\rho} & 2u_m & \frac{2\alpha_1}{3} & &  &  & \frac{gh(\rho_s-\rho_w)}{2\rho} & gh \\
			-2u_m\alpha_1-\frac{gh(\rho-\rho_w)}{2\rho} & 2\alpha_1 & u_m & c_2 &  &  & \frac{gh(\rho_s-\rho_w)}{2\rho} & \\
			-\frac{2\alpha_1^2}{3}&  & \frac{\alpha_1}{3} & u_m & \ddots &  &  & \\
			&  &  & \ddots & \ddots  & c_N &  &  \\
			&  &  &  & a_N & u_m & &  \\
			-c_m u_m & c_m  &   &  & &  & u_m &  \\
			\delta_h & \delta_q & \delta_q &  & \dots  & \delta_q & \delta_c &  
		\end{pmatrix}  
	\end{equation*}
	We write the matrix $\tilde{A}_{ED}$ by using the following notation for conciseness;
	\vspace{4mm}
	\newline $d_1=gh-u_m^2-\frac{\alpha_1^2}{3}-\frac{gh(\rho-\rho_w)}{2\rho} ,\quad d_2=\frac{2\alpha_1}{3}, \quad d_3=-2u_m\alpha_1-\frac{gh(\rho-\rho_w)}{2\rho},$\\
	$d_4=2\alpha_1,\quad d_5=-\frac{2\alpha_1^2}{3},a_2= \frac{\alpha_1}{3}, \quad d_6=gh,\quad d_7=\frac{gh(\rho_s-\rho_w)}{2\rho},$\\
	$d_8=\frac{gh(\rho_s-\rho_w)}{2\rho},\quad d_9= -c_m u_m,\quad d_{10}=c_m$.\\
	So the new matrix is
	\begin{equation*}
		\hspace{-5mm}
		\tilde{A}_{ED}= \begin{pmatrix} 
			& 1 &  &  &  &  & &  \\
			d_1 & 2u_m & d_2 & &  &  & d_7 & d_6 \\
			d_3 & d_4 & u_m & c_2 &  &  & d_8 & \\
			d_5 &  & a_2 & u_m & \ddots &  &  & \\
			&  &  & \ddots & \ddots  & c_N &  &  \\
			&  &  &  & a_N & u_m & &  \\
			d_9 & d_{10}  &   &  & &  & u_m &  \\
			\delta_h & \delta_q & \delta_q &  & \dots  & \delta_q & \delta_c &  
		\end{pmatrix}  
	\end{equation*}
	We write $\widetilde{A}(\textbf{W})= \tilde{A}_{ED}(\textbf{W})-u_m\textbf{I}$ and $\widetilde{\lambda}=\lambda-\textbf{I}$, so that we can compute the characteristics polynomial using
	\begin{align*}
		\chi_{\tilde{A}_{ED}} (\lambda)&= det (\tilde{A}_{ED}(\textbf{W})-\lambda\textbf{I})\\
		& = det (\tilde{A}_{ED}(\textbf{W})-u_m\textbf{I}-\lambda\textbf{I}+u_m\textbf{I})\\
		& = det ( (\tilde{A}_{ED}(\textbf{W})-u_m\textbf{I})- (\lambda - u_m)\textbf{I} )\\ 
		& =  \lvert \widetilde{A}(\textbf{W}) - \widetilde{\lambda} \textbf{I} \rvert\\ 
	\end{align*}
	Therefore,	\begin{align*}
		\lvert \widetilde{A}(\textbf{W}) - \widetilde{\lambda} \textbf{I} \rvert &= \begin{vmatrix} 
			-\widetilde{\lambda}-u_m & 1 &  &  &  &  & &  \\
			d_1 & u_m-\widetilde{\lambda} & d_2 & &  &  & d_7 & d_6 \\
			d_3 & d_4 & -\widetilde{\lambda} & c_2 &  &  & d_8 & \\
			d_5 &  & a_2 & -\widetilde{\lambda} & \ddots &  &  & \\
			&  &  & \ddots & \ddots  & c_N &  &  \\
			&  &  &  & a_N & -\widetilde{\lambda} & &  \\
			d_9 & d_{10}  &   &  & &  & -\widetilde{\lambda} &  \\
			\delta_h & \delta_q & \delta_q &  & \dots  & \delta_q & \delta_c &  -\widetilde{\lambda}-u_m
		\end{vmatrix}\\
		&= (-\widetilde{\lambda}-u_m) \underbrace{\begin{vmatrix} 
				u_m-\widetilde{\lambda} & d_2 & &  &  & d_7 & d_6 \\
				d_4 & -\widetilde{\lambda} & c_2 &  &  & d_8 & \\
				& a_2 & -\widetilde{\lambda} & \ddots &  &  & \\
				&  & \ddots & \ddots  & c_N &  &  \\
				&  &  & a_N & -\widetilde{\lambda} & &  \\
				d_{10}  &   &  & &  & -\widetilde{\lambda} &  \\
				\delta_q & \delta_q &  & \dots  & \delta_q & \delta_c &  -\widetilde{\lambda}-u_m
		\end{vmatrix} }_{M_1}\underbrace { -1 \begin{vmatrix} 
				d_1 &  d_2 & &  &  & d_7 & d_6 \\
				d_3 &  -\widetilde{\lambda} & c_2 &  &  & d_8 & \\
				d_5 & a_2 & -\widetilde{\lambda} & \ddots &  &  & \\
				&  & \ddots & \ddots  & c_N &  &  \\
				&  &    & a_N & -\widetilde{\lambda} & &  \\
				d_9   &   &  & &  & -\widetilde{\lambda} &  \\
				\delta_h  & \delta_q &  & \dots  & \delta_q & \delta_c &  -\widetilde{\lambda}-u_m
		\end{vmatrix}  }_{M_2}
	\end{align*}
	We start with the matrix $M_1$ and find the determinant with respect to the first row
	\begin{equation*}
		\hspace{-5mm}
		M_1= \underbrace {(u_m-\widetilde{\lambda})\begin{vmatrix} 
				-\widetilde{\lambda} & c_2 &  &  & d_8 & \\
				a_2 & -\widetilde{\lambda} & \ddots &  &  & \\
				& \ddots & \ddots  & c_N &  &  \\
				&  & a_N & -\widetilde{\lambda} & &  \\
				&  & &  & -\widetilde{\lambda} &  \\
				\delta_q &  & \dots  & \delta_q & \delta_c &  -\widetilde{\lambda}-u_m
		\end{vmatrix} }_{\textcolor{black}{P_1} }  \underbrace{- d_2 \begin{vmatrix} 
				d_4  & c_2  & &  & d_8 & \\
				& -\widetilde{\lambda}   & \ddots& &  & \\
				&    &\ddots & c_N &  &  \\
				&   & a_N & -\widetilde{\lambda} & &  \\
				d_{10} &   & & & -\widetilde{\lambda} &  \\
				\delta_q  &   \dots & & \delta_q & \delta_c &  -\widetilde{\lambda}-u_m
		\end{vmatrix} }_{\textcolor{black}{P_2}}
	\end{equation*}
	\begin{equation*}
		\hspace{-7mm}
		\underbrace {+(-1)^{N+3}d_7\begin{vmatrix} 
				d_4 & -\widetilde{\lambda} & c_2 &  &   & \\
				& a_2 & -\widetilde{\lambda} & \ddots &   & \\
				&  & \ddots & \ddots  & c_N &  \\
				&  &  & a_N & -\widetilde{\lambda}  &  \\
				d_{10}  &   &  & &   &  \\
				\delta_q & \delta_q &  & \dots  & \delta_q  &  -\widetilde{\lambda}-u_m
		\end{vmatrix}}_{\textcolor{black}{P_3}}  +  \underbrace {(-1)^{N+4}d_6\begin{vmatrix} 
				d_4 & -\widetilde{\lambda} & c_2 &  &  & d_8  \\
				& a_2 & -\widetilde{\lambda} & \ddots &  &   \\
				&  & \ddots & \ddots  & c_N &  \\
				&  &  & a_N & -\widetilde{\lambda} &  \\
				d_{10}  &   &  & &  & -\widetilde{\lambda} \\
				\delta_q & \delta_q &  & \dots  & \delta_q & \delta_c 
		\end{vmatrix} }_{\textcolor{black}{P_4}}
	\end{equation*}
	Next we move to the matrix $M_2$ and find the determinant with respect to the first row
	\begin{equation*}
		\hspace{-5mm}
		M_2 = \underbrace{- d_1 \begin{vmatrix} 
				-\widetilde{\lambda} & c_2 &  &  & d_8 & \\
				a_2 & -\widetilde{\lambda} & \ddots &  &  & \\
				& \ddots & \ddots  & c_N &  &  \\
				&    & a_N & -\widetilde{\lambda} & &  \\
				&  & &  & -\widetilde{\lambda} &  \\
				\delta_q &  & \dots  & \delta_q & \delta_c &  -\widetilde{\lambda}-u_m
		\end{vmatrix}}_{\textcolor{black}{P_5}}   + \underbrace{ d_2 \begin{vmatrix} 
				
				d_3  & c_2 &  &  & d_8 & \\
				d_5 & -\widetilde{\lambda} & \ddots &  &  & \\
				&  & \ddots  & c_N &  &  \\
				&    & a_N & -\widetilde{\lambda} & &  \\
				d_9 &  & &  & -\widetilde{\lambda} &  \\
				\delta_h  &  & \dots  & \delta_q & \delta_c &  -\widetilde{\lambda}-u_m
		\end{vmatrix} }_ {\textcolor{black}{P_6}}
	\end{equation*}
	\begin{equation*}
		\hspace{-5mm}
		\underbrace {-(-1)^{N+3}d_7 \begin{vmatrix} 
				d_3 &  -\widetilde{\lambda} & c_2 &  &   & \\
				d_5 & a_2 & -\widetilde{\lambda} & \ddots &   & \\
				&  & \ddots & \ddots  & c_N   &  \\
				&  &    & a_N & -\widetilde{\lambda}  &  \\
				d_9   &   &  & &   &  \\
				\delta_h  & \delta_q &  & \dots  & \delta_q  &  -\widetilde{\lambda}-u_m
		\end{vmatrix} }_{\textcolor{black}{P_7}} \underbrace {-(-1)^{N+4}d_6 \begin{vmatrix} 
				d_3 &  -\widetilde{\lambda} & c_2 &  &  & d_8  \\
				d_5 & a_2 & -\widetilde{\lambda} & \ddots &  &   \\
				&  & \ddots & \ddots  & c_N &   \\
				&  &    & a_N & -\widetilde{\lambda} &   \\
				d_9   &   &  & &  & -\widetilde{\lambda}   \\
				\delta_h  & \delta_q &  & \dots  & \delta_q & \delta_c 
		\end{vmatrix} }_{\textcolor{black}{P_8}} 
	\end{equation*}
	\begin{align*}
		\text{Now,}\,\textcolor{black}{P_1} &= (u_m - \widetilde{\lambda})
		\left|
		\begin{array}{ccccc}
			-\widetilde{\lambda} & c_2 &        & d_8 & \\
			a_2 & -\widetilde{\lambda} & \ddots &        & \\
			& \ddots & \ddots  & c_N & \\
			&        & a_N     & -\widetilde{\lambda} & \\
			\delta_q & \dots & \delta_q & \delta_c &  -\widetilde{\lambda} - u_m
		\end{array}
		\right|
		= (-\widetilde{\lambda})(\widetilde{\lambda}^2 - u_m^2) \cdot \lvert A_2 \rvert \\[2ex]
		\textcolor{black}{P_2} &= - d_2
		\left|
		\begin{array}{ccccc}
			d_4 & c_2 &        & d_8 & \\
			& -\widetilde{\lambda} & \ddots &        & \\
			&        & \ddots & c_N & \\
			d_{10} &        &        & -\widetilde{\lambda} & \\
			\delta_q & \dots & \delta_q & \delta_c &  -\widetilde{\lambda} - u_m
		\end{array}
		\right|
		= \left[(\widetilde{\lambda} + u_m) d_2 \left( -d_{10} d_8 + (-\widetilde{\lambda}) d_4 \right) \right] \cdot \lvert A_3 \rvert \\[2ex]
		\textcolor{black}{P_3} &= (-1)^{N+3} d_7
		\left|
		\begin{array}{ccccc}
			d_4 & -\widetilde{\lambda} & c_2 &         & \\
			& a_2 & -\widetilde{\lambda} & \ddots & \\
			&     & \ddots & \ddots  & c_N \\
			d_{10} &        &        &        & \\
			\delta_q & \dots & \delta_q & \delta_q & -\widetilde{\lambda} - u_m
		\end{array}
		\right|
		= (\widetilde{\lambda} + u_m) d_7 d_{10} \cdot \lvert A_2 \rvert \\[2ex]
		\textcolor{black}{P_4} &= (-1)^{N+4} d_6
		\left|
		\begin{array}{ccccc}
			d_4 & -\widetilde{\lambda} & c_2 &         & d_8 \\
			& a_2 & -\widetilde{\lambda} & \ddots & \\
			&     & \ddots & \ddots  & c_N \\
			d_{10} &        &        &        & -\widetilde{\lambda} \\
			\delta_q & \dots & \delta_q & \delta_q & \delta_c
		\end{array}
		\right|
	\end{align*}
	\begin{align*}
		& \text{Now,}\, \left|
		\begin{array}{cccccc}
			d_4 & -\widetilde{\lambda} & c_2 &        &        & d_8 \\
			& a_2 & -\widetilde{\lambda} & \ddots &        &     \\
			&     & \ddots & \ddots  & c_N    &     \\
			&     &        & a_N     & -\widetilde{\lambda} & \\
			d_{10} &     &        &        &        & -\widetilde{\lambda} \\
			\delta_q & \delta_q &        & \dots  & \delta_q & \delta_c
		\end{array}
		\right|
		= \textcolor{black}{d_4} 
		\left|
		\begin{array}{ccccc}
			a_2 & -\widetilde{\lambda} & c_3 &        &         \\
			& a_3 & \ddots & \ddots &         \\
			&     & \ddots & \ddots & c_N     \\
			&     &        & a_N & -\widetilde{\lambda} \\
			\delta_q & \delta_q & \dots  & \delta_q & \delta_c
		\end{array}
		\right|\\
		& + \textcolor{black}{(-1)^{N+2} d_{10}}
		\left|
		\begin{array}{ccccc}
			-\widetilde{\lambda} & c_2 &        &        & d_8 \\
			a_2 & -\widetilde{\lambda} & c_3 &        &     \\
			& a_3 & \ddots & \ddots &     \\
			&     & \ddots & \ddots & c_N \\
			\delta_q & \delta_q & \dots & \delta_q & \delta_c
		\end{array}
		\right| 
		+ \textcolor{black}{(-1)^{N+3} \delta_q}
		\left|
		\begin{array}{ccccc}
			-\widetilde{\lambda} & c_2 &        &        & d_8 \\
			a_2 & -\widetilde{\lambda} & c_3 &        &     \\
			& a_3 & \ddots & \ddots &     \\
			&     & \ddots & \ddots & c_N \\
			&     &        & a_N & -\widetilde{\lambda}
		\end{array}
		\right| \\[2ex]
		&= d_4 \cdot \textcolor{black}{(-1)^{2N+1} (-\widetilde{\lambda})}
		\left|
		\begin{array}{ccccc}
			a_2 & -\widetilde{\lambda} & c_3 &        &         \\
			& a_3 & \ddots & \ddots &         \\
			&     & \ddots & \ddots & c_N     \\
			&     &        & a_N & -\widetilde{\lambda} \\
			\delta_q & \delta_q & \dots  & \delta_q & 
		\end{array}
		\right| 
		+ (-1)^{N+2} d_{10} \cdot \textcolor{black}{(-1)^{N+2} d_8}
		\left|
		\begin{array}{ccccc}
			a_2 & -\widetilde{\lambda} & c_3 &        &         \\
			& a_3 & \ddots & \ddots &         \\
			&     & \ddots & \ddots & c_N     \\
			&     &        & a_N & -\widetilde{\lambda} \\
			\delta_q & \delta_q & \dots  & \delta_q &
		\end{array}
		\right| \\
		&\quad + (-1)^{N+2} d_{10} \cdot \textcolor{black}{(-1)^{2N+2} \delta_c}
		\left|
		\begin{array}{ccccc}
			-\widetilde{\lambda} & c_2 &        &        &     \\
			a_2 & -\widetilde{\lambda} & c_3 &        &     \\
			& a_3 & \ddots & \ddots &     \\
			&     & \ddots & \ddots & c_N \\
			&     &        & a_N & -\widetilde{\lambda}
		\end{array}
		\right| 
		+ (-1)^{N+3} \delta_q \cdot \textcolor{black}{(-\widetilde{\lambda})}
		\left|
		\begin{array}{ccccc}
			-\widetilde{\lambda} & c_2 &        &        &     \\
			a_2 & -\widetilde{\lambda} & c_3 &        &     \\
			& a_3 & \ddots & \ddots &     \\
			&     & \ddots & \ddots & c_N \\
			&     &        & a_N & -\widetilde{\lambda}
		\end{array}
		\right| \\[2ex]
		&= (-1)^{2N+1} d_4 (-\widetilde{\lambda}) \lvert B_2 \rvert
		+ (-1)^{2N+4} d_8 d_{10} \lvert B_2 \rvert
		+ (-1)^{N+2} d_{10} \delta_c \lvert A_2 \rvert
		+ (-1)^{N+3} (-\widetilde{\lambda}) \delta_q \lvert A_2 \rvert \\[1ex]
		&= \left[ -(-\widetilde{\lambda}) d_4 + d_8 d_{10} \right] \cdot \lvert B_2 \rvert
		+ \left[ (-1)^{N+2} d_{10} \delta_c + (-1)^{N+3} (-\widetilde{\lambda}) \delta_q \right] \cdot \lvert A_2 \rvert
	\end{align*}
	Therefore,
	\begin{align*}
		\textcolor{black}{P_4}
		&= (-1)^{N+4} d_6 \left\{ \left[ -(-\widetilde{\lambda}) d_4 + d_8 d_{10} \right] \cdot \lvert B_2 \rvert 
		+ \left[ (-1)^{N+2} d_{10} \delta_c + (-1)^{N+3} (-\widetilde{\lambda}) \delta_q \right] \cdot \lvert A_2 \rvert \right\} \\
		&= \left\{ (-1)^{N+4} d_6 \left[ -(-\widetilde{\lambda}) d_4 + d_8 d_{10} \right] \right\} \cdot \lvert B_2 \rvert 
		+ \left[ d_6 d_{10} \delta_c - (-\widetilde{\lambda}) d_6 \delta_q \right] \cdot \lvert A_2 \rvert
	\end{align*}

	\begin{equation*}
		\textcolor{black}{P_5} = -d_1
		\begin{vmatrix}
			-\widetilde{\lambda} & c_2 &        &        & d_8 & \\
			a_2 & -\widetilde{\lambda} & \ddots &        &     & \\
			& \ddots & \ddots  & c_N    &     & \\
			&        & a_N     & -\widetilde{\lambda} & & \\
			&        &         &        & -\widetilde{\lambda} & \\
			\delta_q &        & \dots  & \delta_q & \delta_c & -\widetilde{\lambda} - u_m
		\end{vmatrix}
		= -d_1 (-\widetilde{\lambda} - u_m)(-\widetilde{\lambda}) \lvert A_2 \rvert
	\end{equation*}
	
	\vspace{3mm}
	
	\begin{align*}
		\textcolor{black}{P_6}
		&= d_2 (-\widetilde{\lambda} - u_m) \cdot (-1)^{N+3} d_9 
		\begin{vmatrix}
			c_2 &        &        & d_8 \\
			-\widetilde{\lambda} & c_3 &        &     \\
			a_3 & \ddots & \ddots &     \\
			& \ddots & \ddots & c_N \\
			&        & a_N    & -\widetilde{\lambda}
		\end{vmatrix}
		+ d_2 (-\widetilde{\lambda} - u_m)(-\widetilde{\lambda})
		\begin{vmatrix}
			d_3 & c_2 &        &        \\
			d_5 & -\widetilde{\lambda} & c_3 &        \\
			& a_3 & \ddots & \ddots \\
			&     & \ddots & c_N \\
			&     & a_N    & -\widetilde{\lambda}
		\end{vmatrix} \\
		&= -(-\widetilde{\lambda} - u_m) d_2 d_9 d_8 \lvert A_3 \rvert
		+ (-\widetilde{\lambda})(-\widetilde{\lambda} - u_m) d_2 d_3 \lvert A_3 \rvert 
		- (-\widetilde{\lambda})(-\widetilde{\lambda} - u_m) d_2 c_2 d_5 \lvert A_4 \rvert
	\end{align*}
	
	\begin{align*}
		&P_7 = -(-1)^{N+3} d_7
		\begin{vmatrix}
			d_3 & -\widetilde{\lambda} & c_2 & & & \\
			d_5 & a_2 & -\widetilde{\lambda} & \ddots & & \\
			& & \ddots & \ddots & c_N & \\
			& & & a_N & -\widetilde{\lambda} & \\
			d_9 & & & & & \\
			\delta_h & \delta_q & \cdots & \delta_q & \delta_q & -\widetilde{\lambda} - u_m
		\end{vmatrix}
		= (-\widetilde{\lambda} - u_m) d_7 d_9 \cdot \lvert A_2 \rvert\\
		&\textcolor{black}{P_8} = -(-1)^{N+4} d_6 
		\begin{vmatrix}
			d_3 & -\widetilde{\lambda} & c_2 & & & d_8 \\
			d_5 & a_2 & -\widetilde{\lambda} & \ddots & & \\
			& & \ddots & \ddots & c_N & \\
			& & & a_N & -\widetilde{\lambda} & \\
			d_9 & & & & & -\widetilde{\lambda} \\
			\delta_h & \delta_q & \cdots & \delta_q & \delta_q & \delta_c
		\end{vmatrix}
	\end{align*}
	\begin{align*}
		&\text{Now,}\, \left|
		\begin{array}{cccccc}
			d_3 & -\widetilde{\lambda} & c_2 &        &        & d_8 \\
			d_5 & a_2 & -\widetilde{\lambda} & \ddots &        &     \\
			&     & \ddots & \ddots & c_N    &     \\
			&     &        & a_N & -\widetilde{\lambda} &     \\
			d_9 &     &        &     &        & -\widetilde{\lambda} \\
			\delta_h & \delta_q &        & \dots  & \delta_q & \delta_c
		\end{array}
		\right|
		= \textcolor{black}{d_3}
		\left|
		\begin{array}{ccccc}
			a_2 & -\widetilde{\lambda} & c_3 &        &        \\
			& a_3 & -\widetilde{\lambda} & \ddots &        \\
			&     & \ddots & \ddots & c_N \\
			&     &        & a_N & -\widetilde{\lambda} \\
			\delta_q & \delta_q &        & \dots  & \delta_c
		\end{array}
		\right|
		- \textcolor{black}{d_5}
		\left|
		\begin{array}{ccccccc}
			-\widetilde{\lambda} & c_2 &        &        &        &        & d_8 \\
			& a_3 & -\widetilde{\lambda} & c_4 &        &        &     \\
			&     & a_4  & -\widetilde{\lambda} & \ddots &        &     \\
			&     &      & \ddots & \ddots & c_N &     \\
			&     &      &        & a_N  & -\widetilde{\lambda} &     \\
			&     &      &        &      &        & -\widetilde{\lambda} \\
			\delta_q &          & \dots  &        & \delta_q & \delta_q & \delta_c
		\end{array}
		\right| \\
		&\quad + (-1)^{N+2} \textcolor{black}{d_9}
		\left|
		\begin{array}{ccccc}
			-\widetilde{\lambda} & c_2 &        &        & d_8 \\
			a_2 & -\widetilde{\lambda} & \ddots &        &     \\
			& \ddots & \ddots  & c_N    &     \\
			&        & a_N & -\widetilde{\lambda} &     \\
			\delta_q &        & \dots  & \delta_q & \delta_c
		\end{array}
		\right|
		+ (-1)^{N+3} \textcolor{black}{\delta_h}
		\left|
		\begin{array}{ccccc}
			-\widetilde{\lambda} & c_2 &        &        & d_8 \\
			a_2 & -\widetilde{\lambda} & \ddots &        &     \\
			& \ddots & \ddots  & c_N    &     \\
			&        & a_N & -\widetilde{\lambda} &     \\
			&        &     &        & -\widetilde{\lambda}
		\end{array}
		\right| \\
		&= d_3 \cdot \textcolor{black}{(-1)^{2N+1} (-\widetilde{\lambda})}
		\left|
		\begin{array}{ccccc}
			a_2 & -\widetilde{\lambda} & c_3 &        &        \\
			& a_3 & -\widetilde{\lambda} & \ddots &        \\
			&     & \ddots & \ddots & c_N \\
			&     &        & a_N & -\widetilde{\lambda} \\
			\delta_q & \delta_q &        & \dots & \delta_q
		\end{array}
		\right|
		- d_5 \cdot \textcolor{black}{(-1)^{2N+1} (-\widetilde{\lambda})}
		\left|
		\begin{array}{cccccc}
			-\widetilde{\lambda} & c_2 &        &        &        &     \\
			& a_3 & -\widetilde{\lambda} & c_4 &        &     \\
			&     & a_4  & -\widetilde{\lambda} & \ddots &     \\
			&     &      & \ddots & \ddots & c_N \\
			&     &      &        & a_N  & -\widetilde{\lambda} \\
			\delta_q &           & \dots  &        & \delta_q & \delta_q
		\end{array}
		\right| \\
		&\quad +  \textcolor{black}{d_9 d_8}
		\left|
		\begin{array}{ccccc}
			a_2 & -\widetilde{\lambda} & c_3 &        &        \\
			& a_3 & -\widetilde{\lambda} & \ddots &        \\
			&     & \ddots & \ddots & c_N \\
			&     &        & a_N & -\widetilde{\lambda} \\
			\delta_q &        & \dots &        & \delta_q
		\end{array}
		\right| 
		+ (-1)^{3N+5} d_9 \cdot \textcolor{black}{\delta_c}
		\left|
		\begin{array}{cccc}
			-\widetilde{\lambda} & c_2 &        &     \\
			a_2 & -\widetilde{\lambda} & \ddots &     \\
			& \ddots & \ddots & c_N \\
			&        & a_N & -\widetilde{\lambda}
		\end{array}
		\right|
		+ (-1)^{N+3} \delta_h \cdot \textcolor{black}{(-\widetilde{\lambda})}
		\left|
		\begin{array}{cccc}
			-\widetilde{\lambda} & c_2 &        &     \\
			a_2 & -\widetilde{\lambda} & \ddots &     \\
			& \ddots & \ddots & c_N \\
			&        & a_N & -\widetilde{\lambda}
		\end{array}
		\right| \\
		&= -(-\widetilde{\lambda}) d_3 \lvert B_2 \rvert
		+ (-\widetilde{\lambda}) d_5 \cdot \textcolor{black}{(-\widetilde{\lambda})} \lvert B_3 \rvert
		+ (-\widetilde{\lambda}) d_5 \cdot \textcolor{black}{(-1)^{N+1} \delta_q c_2} \lvert A_4 \rvert 
		+ d_9 d_8 \lvert B_2 \rvert
		+ (-1)^{N+2} d_9 \delta_c \lvert A_2 \rvert
		+ (-1)^{N+3} \delta_h (-\widetilde{\lambda}) \lvert A_2 \rvert
	\end{align*}
	After applying the recursion formula 
	$\lvert A_4 \rvert = -\frac{1}{a_2c_2}\left[\lvert A_2 \rvert +\widetilde{\lambda} \lvert A_3 \rvert \right]$ and $\lvert B_3 \rvert = \frac{1}{a_2}\left[\lvert B_2 \rvert + (-1)^{N+2}\delta_q \lvert A_3 \rvert \right]$ to $P_8$ yields,
	\begin{equation*}
		\begin{aligned}
			P_8 & = \left[-d_6d_9 \delta_c+ (-\widetilde{\lambda}) d_6\delta_h- (-\widetilde{\lambda})\frac{d_5d_6}{a_2}\delta_q \right]\lvert A_2 \rvert-
			\left[  (-\widetilde{\lambda})(-\widetilde{\lambda})\frac{d_5d_6}{a_2}\delta_q - (-\widetilde{\lambda})(-\widetilde{\lambda})\frac{d_5d_6}{a_2}\delta_q \right]\lvert A_3 \rvert\\
			+&(-1)^{N+5}\left[-(-\widetilde{\lambda}) d_6d_3+ d_6d_9d_8+ (-\widetilde{\lambda})(-\widetilde{\lambda})\frac{d_5d_6}{a_2}        \right] \lvert B_2 \rvert.
		\end{aligned}
	\end{equation*}
	If we add all together, we get
	\begin{equation*}
		\begin{split} 
			\chi_{\tilde{A}_{ED}}(\lambda)&=  \lvert A_2\rvert \left\{ ( -\widetilde{\lambda})(-\widetilde{\lambda}-u_m) \left((\widetilde{\lambda}^2-u_m^2)-d_1+\frac{d_2d_5}{a_2}-d_7d_{10}\right)\right\} \\
			&+\lvert A_2\rvert \left\{  ( -\widetilde{\lambda})d_6\delta_h +  d_6 \delta_c ((-\widetilde{\lambda}-u_m)d_{10}-d_9)+ (-\widetilde{\lambda}) d_6\delta_q \left( (\widetilde{\lambda}+u_m)-  \frac{d_5}{a_2} \right)\right\}\\
			&+\lvert A_3\rvert \left\{  ( -\widetilde{\lambda})( -\widetilde{\lambda}-u_m)\left[(\widetilde{\lambda}+u_m) d_2d_4 +d_2d_8d_{10} +d_2d_3+\frac{d_2d_5}{a_2}\widetilde{\lambda}   \right] \right\}+\lvert B_2 \rvert ( -\widetilde{\lambda}-u_m)\\
			&  \left \{(-1)^{N+4}d_6\left[- (-\widetilde{\lambda}) d_4 +d_8 d_{10} \right] \right \}
			+\lvert B_2 \rvert(-1)^{N+5}\left[-(-\widetilde{\lambda}) d_6d_3+ d_6d_9d_8+ (-\widetilde{\lambda})(-\widetilde{\lambda})\frac{d_5d_6}{a_2} \right] \\
			&=\lvert A_2\rvert ( -\widetilde{\lambda})\left\{ -(\widetilde{\lambda}+u_m) \left((\widetilde{\lambda}^2-u_m^2)-d_1+\frac{d_2d_5}{a_2}-d_7d_{10}\right) \right\}\\
			&+ \lvert A_2\rvert ( -\widetilde{\lambda}) d_6 \left( (\widetilde{\lambda}+u_m)\delta_q-\frac{d_5}{a_2}\delta_q+c_m\delta_c+\delta_h\right) + \lvert A_3\rvert \cdot 0 + \lvert B_2\rvert \cdot 0\\
			&= \lvert A_2\rvert ( -\widetilde{\lambda}) \left\{ (-\widetilde{\lambda}-u_m) \left(\widetilde{\lambda}^2-gh-\alpha_1^2  \right) + gh(\delta_h+c_m\delta_c+ 2\alpha_1\delta_q)+ \delta_q \left( \widetilde{\lambda}+u_m \right)\right\}
		\end{split}
	\end{equation*}
	
	Finally, substitute $\widetilde{\lambda}=\lambda-u_m$ and we get the characteristics polynomial
	\begin{equation*}
		\chi_{\tilde{A}_{ED}}(\lambda)= (u_m-\lambda)\left\{(-\lambda) \left( (\lambda-u_m)^2-gh-\alpha_1^2\right)+gh(\delta_h+\lambda\delta_q+c_m\delta_c+2\alpha_1\delta_q)\right\}\cdot \chi_{A_2}(\lambda-u_m).
	\end{equation*}

%% file: refs.bib
@Article{CiCP-25-669,
	author = {Kowalski, Julia and Torrilhon, Manuel},
	title = {Moment Approximations and Model Cascades for Shallow Flow},
	journal = {Communications in Computational Physics},
	year = {2019},
	volume = {25},
	number = {3},
	pages = {669--702},
	issn = {1991-7120},
	doi = {10.4208/cicp.OA-2017-0263}
}

@article{koellermeier2020analysis,
	author = {Koellermeier, Julian and Rominger, Marvin},
	title = {Analysis and Numerical Simulation of Hyperbolic Shallow Water Moment Equations},
	journal = {Communications in Computational Physics},
	volume = {28},
	number = {3},
	pages = {1038--1084},
	year = {2020},
	doi = {10.4208/cicp.OA-2019-0065},
	%url = {https://doi.org/10.4208/cicp.OA-2019-0065}
}

@article{Garres,
	author = {Garres-Díaz, José and Castro Díaz, Manuel Jesús and Koellermeier, Julian and Morales de Luna, Tomás},
	title = {Shallow Water Moment Models for Bedload Transport Problems},
	journal = {Communications in Computational Physics},
	year = {2021},
	volume = {30},
	number = {3},
	pages = {903--941},
	abstract = {<p style="text-align: justify;">In this work a simple but accurate shallow model for bedload sediment
	transport is proposed. The model is based on applying the moment approach to the
	Shallow Water Exner model, making it possible to recover the vertical structure of the
	flow. This approach allows us to obtain a better approximation of the fluid velocity
	close to the bottom, which is the relevant velocity for the sediment transport. A general Shallow Water Exner moment model allowing for polynomial velocity profiles of
	arbitrary order is obtained. A regularization ensures hyperbolicity and easy computation of the eigenvalues. The system is solved by means of an adapted IFCP scheme
	proposed here. The improvement of this IFCP type scheme is based on the approximation of the eigenvalue associated to the sediment transport. Numerical tests are
	presented which deal with large and short time scales. The proposed model allows to
	obtain the vertical structure of the fluid, which results in a better description on the
	bedload transport of the sediment layer.</p>},
	issn = {1991-7120},
	doi = {10.4208/cicp.OA-2020-0152},
%	url = %{https://global-sci.com/article/79588/shallow-water-moment-models-for-bedload-transport-problems}
	}

@article{huang2022equilibrium,
	author = {Huang, Qian and Koellermeier, Julian and Yong, Wen-An},
	title = {Equilibrium Stability Analysis of Hyperbolic Shallow Water Moment Equations},
	journal = {Mathematical Methods in the Applied Sciences},
	volume = {45},
	number = {10},
	pages = {6459--6480},
	year = {2022},
	doi = {10.1002/mma.8180},
	%url = {https://onlinelibrary.wiley.com/doi/abs/10.1002/mma.8180}
}

@article{koellermeier2022steady,
	author = {Koellermeier, Julian and Pimentel-García, Ernesto},
	title = {Steady States and Well-balanced Schemes for Shallow Water Moment Equations with Topography},
	journal = {Journal of Computational Physics},
	volume = {429},
	pages = {109951},
	year = {2021},
	doi = {10.1016/j.jcp.2020.109951},
	%url = {https://www.sciencedirect.com/science/article/pii/S0096300322002417}
}

@article{garres2023general,
	author = {Garres-Díaz, José and Escamilla-Sánchez, Alejandro and Morales de Luna, Tomás and Castro D{\'\i}az, Manuel Jes{\'u}s},
	title = {A General Vertical Decomposition of Euler Equations: Multilayer-Moment Models},
	journal = {Applied Numerical Mathematics},
	volume = {183},
	pages = {236--262},
	year = {2023},
	doi = {10.1016/j.apnum.2022.09.004},
	%url = {https://www.sciencedirect.com/science/article/pii/S0168927422002288}
}

@article{amrita2022projective,
	author = {Amrita, Amrita and Koellermeier, Julian},
	title = {Projective Integration for Hyperbolic Shallow Water Moment Equations},
	journal = {Axioms},
	volume = {11},
	number = {5},
	pages = {235},
	year = {2022},
	doi = {10.3390/axioms11050235},
	%url = {https://www.mdpi.com/2075-1680/11/5/235}
}

@article{gonzalez2020robust,
  author    = {González-Aguirre, Juan Carlos and Castro, Manuel Jesús and Morales de Luna, Tomás},
  title     = {A robust model for rapidly varying flows over movable bottom with suspended and bedload transport: Modelling and numerical approach},
  journal   = {Advances in Water Resources},
  volume    = {140},
  pages     = {103575},
  year      = {2020},
  doi       = {10.1016/j.advwatres.2020.103575},
  %url       = {https://doi.org/10.1016/j.advwatres.2020.103575}
}

@article{steldermann2023shallow,
    title={Shallow Moments to Capture Vertical Structure in Open Curved Shallow Flow},
    author={Steldermann, Ingo and Torrilhon, Manuel and Kowalski, Julia},
    journal = {Journal of Computational and Theoretical Transport},
    volume = {52},
    number = {7},
    pages = {475--505},
    year = {2023},
    publisher = {Taylor \& Francis},
    doi = {10.1080/23324309.2023.2284202},
    %URL = {https://doi.org/10.1080/23324309.2023.2284202},
    %eprint = {https://doi.org/10.1080/23324309.2023.2284202}
}

@article{exner1920physik,
	title={Zur Physik der Dünen},
	author={Exner Ewarten, Felix Maria},
        journal   = {Sitzungsberichte der Mathematisch-Naturwissenschaftlichen Klasse der Akademie der Wissenschaften in Wien},
        volume    = {129},
        pages     = {929--952},
        year      = {1920},
}

@article{exner1925uber,
    author={Exner Ewarten, Felix Maria},
      title     = {Über die Wechselwirkung zwischen Wasser und Geschiebe in Flüssen},
      journal   = {Sitzungsberichte der Akademie der Wissenschaften in Wien, Mathematisch-Naturwissenschaftliche Klasse},
      volume    = {134},
      number    = {2a},
      pages     = {165--204},
      year      = {1925}
      %note      = {Accessed: 2025-03-20},
      %url       = {https://books.google.com/books/about/Über_die_Wechselwirkung_zwischen_Wasser.html?id=7bjqZwEACAAJ}
}

@inproceedings{meyer1948formulas,
	author={Eugen, Meyer Peter and Robert, Müller},
         title     = {Formulas for Bed-Load Transport},
         booktitle = {Proceedings of the Second Meeting of the International Association for Hydraulic Structures Research},
         year      = {1948},
        pages     = {39--64},
        address   = {Delft, Netherlands},
        %url       = {https://repository.tudelft.nl/record/uuid%3A4fda9b61-be28-4703-ab06-43cdc2a21bd7}
}

@article{cordier2011bedload,
    title = {Bedload transport in shallow water models: Why splitting (may) fail, how hyperbolicity (can) help},
    journal = {Advances in Water Resources},
    volume = {34},
    number = {8},
    pages = {980-989},
    year = {2011},
    issn = {0309-1708},
    doi = {10.1016/j.advwatres.2011.05.002},
    %url = {https://www.sciencedirect.com/science/article/pii/S0309170811000935},
    author = {Cordier, Stéphane and Le, Minh Hanh and Morales de Luna, Tomás},
    keywords = {Shallow water system, Exner equation, Sediment transport, Hyperbolicity, Stability, Splitting methods},
    abstract = {In this paper, we are concerned with sediment transport models consisting of a shallow water system coupled with the so called Exner equation to describe the evolution of the topography. We show that, for some bedload transport models like the well-known Meyer-Peter and Müller model, the system is hyperbolic and, thus, linearly stable, only under some constraint on the velocity. In practical situations, this condition is hopefully fulfilled. Numerical approximations of such system are often based on a splitting method, solving first shallow water equation on a time step and, updating afterwards the topography. It is shown that this strategy can create spurious/unphysical oscillations which are related to the study of hyperbolicity. Using an upper bound of the largest eigenvector may improve the results although the instabilities cannot be always avoided, e.g. in supercritical regions.}
}

@article{garcia1993experiments,
	title={Experiments on the entrainment of sediment into suspension by a dense bottom current},
	author={Garcia, Marcelo Horacio and Parker, Gary},
	journal={Journal of Geophysical Research: Oceans},
	volume={98},
	number={C3},
	pages={4793--4807},
	year={1993},
	publisher={Wiley Online Library},
        doi     = {10.1029/92JC02404},
        %url     = {https://doi.org/10.1029/92JC02404}
}

@book{zhang1993sedimentation,
	title={Sedimentation research in China: Systematic selections},
	author={Zhang, Ruijin and Xie, Jianheng},
          publisher = {China Water and Power Press},
          year      = {1993},
          address   = {Beijing},
          isbn      = {9787120019433},
}

@article{fernandez1976erosion,
	title={Erosion and transport of bed-load sediment},
	author={Fernandez Luque, Rafael and Van Beek, Rens},
        journal = {Journal of Hydraulic Research},
        volume = {14},
        number = {2},
        pages = {127--144},
        year = {1976},
        publisher = {IAHR Website},
        doi = {10.1080/00221687609499677},
        %URL = {https://doi.org/10.1080/00221687609499677},
        %eprint = {https://doi.org/10.1080/00221687609499677}
        }

@article{nielsen1992coastal,
	title={Coastal bottom boundary layers and sediment transport},
	author={Nielsen, Peter},
	journal={World Sci., River Edge, NJ},
	year={1992}
}

@article{ribberink1994sediment,
	   author={Ribberink, Jan Sjoerd and Al-Salem, Abdullah},
          title     = {Sediment transport in oscillatory boundary layers in cases of rippled beds and sheet flow},
          journal   = {Journal of Geophysical Research},
          volume    = {99},
          number    = {C6},
          pages     = {12,707--12,727},
          year      = {1994},
          doi       = {10.1029/94JC00380},
          %url       = {https://doi.org/10.1029/94JC00380}
}

@article{camenen2006phase,
    title = {Phase-lag effects in sheet flow transport},
    journal = {Coastal Engineering},
    volume = {53},
    number = {5},
    pages = {531-542},
    year = {2006},
    issn = {0378-3839},
    doi = {10.1016/j.coastaleng.2005.12.003},
    %url = {https://www.sciencedirect.com/science/article/pii/S0378383905001766},
    author = {Camenen, Benoît  and Larson, Magnus},
    keywords = {Inception of sheet flow, Phase-lag effects, Net sediment transport, Bed load},
    abstract = {The inception of the sheet flow regime as well as the effects of the phase lag when the sheet flow regime is established were investigated for oscillatory flows and combined steady and oscillatory flows. A new criterion for the inception of sheet flow is proposed based on around 300 oscillatory flow cases from experiments. This criterion was introduced in the Camenen and Larson [Camenen, B., Larson, M., 2005. A bedload sediment transport formula for the nearshore. Estuarine, Coastal and Shelf Science 63, 249–260.] bed load formula in order to take into account phase-lag effects in the sheet flow regime. The modification of the Camenen and Larson formula significantly improves the overall agreement with data and yields a correct behavior in relation to some of the main governing parameters, which are the median grain size d50, the orbital wave velocity Uw, and the wave period Tw. The calibration of the new formula was based on more than 200 experimental data values on the net sediment transport rate for a full wave cycle. A conceptual model was also proposed to estimate the ratio between sediment transport rate with and without phase lag, (rpl=qs,net/qs,net,ϕ=0). This simple model provides accurate results and may be used together with any quasi-steady model for bed load transport.}
}

@article{del2023lagrange,
	author={Del Grosso, Alessia and Castro D{\'\i}az, Manuel Jes{\'u}s and Chalons, Christophe and Morales de Luna, Tomás},
	title = {On Lagrange-Projection Schemes for Shallow Water Flows over Movable Bottom with Suspended and Bedload Transport},
journal = {Numerical Mathematics: Theory, Methods and Applications},
year = {2023},
volume = {16},
number = {4},
pages = {1087--1126},
abstract = {<p style="text-align: justify;">In the present work we aim to simulate shallow water flows over movable bottom with suspended and bedload transport. In order to numerically approximate such a system, we proceed step by step. We start by considering shallow
water equations with non-constant density of the mixture water-sediment. Then,
the Exner equation is included to take into account bedload sediment transport.
Finally, source terms for friction, erosion and deposition processes are considered.
Indeed, observe that the sediment particle could go in suspension into the water or
being deposited on the bottom. For the numerical scheme, we rely on well-balanced
Lagrange-projection methods. In particular, since sediment transport is generally
a slow process, we aim to develop semi-implicit schemes in order to obtain fast simulations. The Lagrange-projection splitting is well-suited for such a purpose as it
entails a decomposition of the (fast) acoustic waves and the (slow) material waves
of the model. Hence, in subsonic regimes, an implicit approximation of the acoustic
equations allows us to neglect the corresponding CFL condition and to obtain fast
numerical schemes with large time step.</p>},
issn = {2079-7338},
doi = {10.4208/nmtma.OA-2023-0082},
%url = {https://global-sci.com/article/90240/on-lagrange-projection-schemes-for-shallow-water-flows-over-movable-bottom-with-suspended-and-bedload-transport}
}

@book{evans2022partial,
       author={Evans, Lawrence Craig},
      title     = {Partial Differential Equations},
      series    = {Graduate Studies in Mathematics},
      volume    = {19},
      edition   = {2},
      publisher = {American Mathematical Society},
      year      = {2010},
      isbn      = {9780821849743},
      %url       = {https://bookstore.ams.org/gsm-19}
}

@article{Spinewine01122007,
author = {Spinewine, Benoit and Zech, Yves},
title = {Small-scale laboratory dam-break waves on movable beds},
journal = {Journal of Hydraulic Research},
volume = {45},
number = {sup1},
pages = {73--86},
year = {2007},
publisher = {IAHR Website},
doi = {10.1080/00221686.2007.9521834},
}

@article{castro2012class,
	title={A class of computationally fast first order finite volume solvers: PVM methods},
	author={Castro D{\'\i}az, Manuel Jes{\'u}s and Fernandez-Nieto, Enrique},
	journal={SIAM Journal on Scientific Computing},
	volume={34},
	number={4},
	pages={A2173--A2196},
	year={2012},
	publisher={SIAM}
}

@article{diaz2008sediment,
      author  = {Castro D{\'\i}az, Manuel Jes{\'u}s and Fern{\'a}ndez-Nieto, Enrique Domingo and Ferreiro, Ana María},
      title   = {Sediment transport models in shallow water equations and numerical approach by high order finite volume methods},
      journal = {Computers and Fluids},
	  volume = {37},
	  number = {3},
	  pages = {299-316},
  	  year = {2008},
	 issn = {0045-7930},
     doi     = {10.1016/j.compfluid.2007.07.017},
}

@article{subhasish2014fluvial,
	title={Fluvial Hydrodynamics: Hydrodynamic and Sediment Transport Phenomena},
	author={Subhasish, Dey},
	series    = {GeoPlanet: Earth and Planetary Sciences},
        publisher = {Springer},
        year      = {2014},
        isbn      = {978-3-642-19061-2},
        %url       = {https://link.springer.com/book/10.1007/978-3-642-19062-9}
}

@article{audusse2004fast,
    author = {Audusse, Emmanuel and Bouchut, Fran\c{c}ois and Bristeau, Marie-Odile and Klein, Rupert and Perthame, Beno\i{}t},
    title = {A Fast and Stable Well-Balanced Scheme with Hydrostatic Reconstruction for Shallow Water Flows},
    journal = {SIAM Journal on Scientific Computing},
    volume = {25},
    number = {6},
    pages = {2050-2065},
    year = {2004},
    doi = {10.1137/S1064827503431090},
    %URL = {https://doi.org/10.1137/S1064827503431090},
    %eprint = {https://doi.org/10.1137/S1064827503431090}
}

@article{meng2020localized,
      title={Localized exponential time differencing method for shallow water equations: Algorithms and numerical study},
      author={Meng, Xucheng and Hoang, Thi Thao Phuong and Wang, Zhu and Ju, Lili},
      journal   = {Communications in Computational Physics},
      volume    = {29},
      number    = {1},
      pages     = {80--110},
      year      = {2020},
      doi       = {10.4208/cicp.OA-2019-0214},
      %url       = {https://global-sci.org/intro/article_detail/cicp/18423.html},
}

@article{zhao2019depth,
	title={A depth-averaged non-cohesive sediment transport model with improved discretization of flux and source terms},
	author={Zhao, Jiaheng and {\"O}zgen-Xian, Ilhan and Liang, Dongfang and Wang, Tian and Hinkelmann, Reinhard},
        journal = {Journal of Hydrology},
        volume = {570},
        pages = {647-665},
        year = {2019},
        issn = {0022-1694},
        doi = {10.1016/j.jhydrol.2018.12.059},
        %url = {https://www.sciencedirect.com/science/article/pii/S002216941930037X},
}

@article{gonzalez2023numerical,
	title={Numerical simulation of bed load and suspended load sediment transport using well-balanced numerical schemes},
	author={Gonz{\'a}lez-Aguirre, Juan Carlos and Gonz{\'a}lez-V{\'a}zquez, José Antonio and Alavez-Ram{\'\i}rez, Justino and Silva, Rodolfo and V{\'a}zquez-Cend{\'o}n, María Elena},
	journal={Communications on Applied Mathematics and Computation},
	volume={5},
	number={2},
	pages={885--922},
	year={2023},
	publisher={Springer},
        doi       = {10.1007/s42967-021-00162-1},
        %url       = {https://doi.org/10.1007/s42967-021-00162-1}
}

@article{cantero2019vertically,
    title = {Vertically-averaged and moment equations for flow and sediment transport},
    journal = {Advances in Water Resources},
    volume = {132},
    pages = {103387},
    year = {2019},
    issn = {0309-1708},
    doi = {10.1016/j.advwatres.2019.103387},
    %url = {https://www.sciencedirect.com/science/article/pii/S0309170819302672},
    author = {Cantero-Chinchilla, Francisco Nicolás  and Castro-Orgaz, Oscar and Khan, Abdul Aziz},
    keywords = {Debris flows, Erosion processes, VAM model, Weighted residual method},
    abstract = {Simulation of river flow processes including sediment transport is usually conducted using the shallow water flow equations, which are based on a hydrostatic pressure distribution. To increase the accuracy of predictions in a variety of scenarios involving horizontal length scales of the order of vertical length scales, an improved representation of the vertical flow structure is necessary. The mathematical approximation to field variables like the velocity and fluid pressure must be enhanced during the depth-integrating process. Therefore, this paper presents a 1D non-hydrostatic flow and sediment transport model developed by using the method of the weighted residuals into the RANS equations. Using continuity, momentum, and moment equations, the fluid pressure distribution is modelled using a quadratic predictor with perturbation parameters to deviate the vertical momentum balance from the hydrostatic law. The flow equations are a generalized non-hydrostatic flow solver, where the fluid density variation due to suspension of sediments and the bed deformation due to erosion-sedimentation processes are accounted for. A hybrid semi-implicit finite volume-finite difference numerical scheme is developed to solve the system of conservation laws. Two different approaches are used to model the sediment transport processes: (i) Unified computation of the total-load transport and (ii) separate computation of suspended and bed loads. The accuracy of the non-hydrostatic model is demonstrated by comparison with experimental data, highlighting better results accounting for separate determinations of the suspended and bed loads in highly erosive flows.}
}

@article{cozzolino2014novel,
    author = {Cozzolino,Luca and Cimorelli ,Luigi and Covelli, Carmine  and Della Morte, Renata and Pianese, Domenico},
    title = {Novel Numerical Approach for 1D Variable Density Shallow Flows over Uneven Rigid and Erodible Beds},
    journal = {Journal of Hydraulic Engineering},
    volume = {140},
    number = {3},
    pages = {254-268},
    year = {2014},
    doi = {10.1061/(ASCE)HY.1943-7900.0000821},
    %URL = {https://ascelibrary.org/doi/abs/10.1061/%28ASCE%29HY.1943-7900.0000821},
    %eprint = {https://ascelibrary.org/doi/pdf/10.1061/%28ASCE%29HY.1943-7900.0000821}
        abstract = { The numerical modeling of hyperconcentrated shallow flows is a challenging task because they exhibit special features, such as propagation over dry beds, profound bed elevation modifications owing to erosion or deposition phenomena, and flow discontinuities. In this paper, a novel depth-positivity preserving Harten, Lax, and van Leer—contact (HLLC) Riemann solver is devised in order to approximate the solution of the Riemann problem for the 1D (one-dimensional) hyperconcentrated shallow flows equations over horizontal beds. The solver is used as a building block for the construction of hyperconcentrated shallow flows (HCSF), a well-balanced finite-volume scheme for the solution of the hyperconcentrated shallow flows equations with variable elevation. HCSF is able to handle the case of dry beds, to take into account the variability of the topography also in the presence of bed discontinuities, considering the flow resistance and the mass exchange between the flowing mixture and the mobile bed. The numerical tests carried out confirm the well-balancing property of the scheme proposed, the robustness in the presence of dry beds, the ability to approximate the analytic solution of problems with smooth or discontinuous beds, and the ability to reproduce reasonably the results of a laboratory experiment. }
}

@article{kubo2004experimental,
    title={Experimental and numerical study of topographic effects on deposition from two-dimensional, particle-driven density currents},
    journal = {Sedimentary Geology},
    volume = {164},
    number = {3},
    pages = {311-326},
    year = {2004},
    issn = {0037-0738},
    doi = {10.1016/j.sedgeo.2003.11.002},
    %url = {https://www.sciencedirect.com/science/article/pii/S0037073803003300},
    author = {Kubo, Yasufumio},
    keywords = {Turbidity current, Tank experiments, Numerical model, Topographic effect},
    abstract = {Topographic effects on deposition from particle-driven density current were investigated. The laboratory experiments were carried out on topography consisting of a ramp and a series of humps. The results show a localized increase in deposit distribution downstream of the slope break and on the upslope of a hump. A numerical model is developed to predict the topographic effects on deposit distribution. The model, based on layer-averaged Navier–Stokes equations, is applied to the experiments, and the process of the topographic influence is analyzed. The preferential deposition downstream of the slope break is interpreted to result from deceleration of the flow, which increases deposit distribution through loss of capacity of the flow and longer duration of the flow passage. The former effect is more influential than the latter, while additional effects of a hydraulic jump may be imposed. Increased deposit on the upslope of a hump is attributed to partial blocking of the flow at the hump crest, as well as to the differential deposition due to deceleration on the upslope.}
}

@article{nakajima2002laboratory,
	title={Laboratory experiments and numerical simulation of sediment-wave formation by turbidity currents},
	author={Kubo, Yasufumi and Nakajima, Takashi},
        journal = {Marine Geology},
        volume = {192},
        number = {1},
        pages = {105-121},
        year = {2002},
        issn = {0025-3227},
        doi = {10.1016/S0025-3227(02)00551-0},
        %url = {https://www.sciencedirect.com/science/article/pii/S0025322702005510},
}

@article{morales2009shallow,
	title={On a shallow water model for the simulation of turbidity currents},
	author={Morales de Luna, Tom{\'a}s and Castro D{\'\i}az, Manuel Jes{\'u}s and Par{\'e}s Madro{\~n}al, Carlos and Fern{\'a}ndez Nieto, Enrique Domingo},
	journal = {Communications in Computational Physics},
	year = {2009},
	volume = {6},
	number = {4},
	pages = {848--882},
	abstract = {<p style="text-align: justify;">We present a model for hyperpycnal plumes or turbidity currents that takes
	into account the interaction between the turbidity current and the bottom, considering
	deposition and erosion effects as well as solid transport of particles at the bed load
	due to the current. Water entrainment from the ambient water in which the turbidity
	current plunges is also considered. Motion of ambient water is neglected and the rigid
	lid assumption is considered. The model is obtained as a depth-average system of
	equations under the shallow water hypothesis describing the balance of fluid mass,
	sediment mass and mean flow. The character of the system is analyzed and numerical
	simulations are carried out using finite volume schemes and path-conservative Roe
	schemes.</p>},
	issn = {1991-7120},
	doi = { https://doi.org/2009-CiCP-7709},
	%url = %{https://global-sci.com/article/81173/on-a-shallow-water-model-for-the-simulation-of-turbidity-currents}
	}

@article{parker1986self,
	title={Self-accelerating turbidity currents},
	author={Parker, Gary and Fukushima, Yusuke and Pantin, Henry Moir},
	journal={Journal of Fluid Mechanics},
	volume={171},
	pages={145--181},
	year={1986},
	publisher={Cambridge University Press},
        doi       = {10.1017/S0022112086001404},
        %url       = {https://doi.org/10.1017/S0022112086001404},
}

@phdthesis{boittin2019modeling,
      TITLE = {{Modeling, analysis and simulation of two geophysical flows. Sediment transport and variable density flows}},
      AUTHOR = {Boittin, L{\'e}a},
      %URL = {https://theses.hal.science/tel-02935869},
      NUMBER = {2019SORUS033},
      SCHOOL = {{Sorbonne Universit{\'e}}},
      YEAR = {2019},
      MONTH = Apr,
      KEYWORDS = {Geophysical flows ; Sediment transport ; Non-local flux ; Variable density flow ; Multilayer model ; Numerical simulation ; {\'E}coulements g{\'e}ophysiques ; Transport s{\'e}dimentaire ; Flux non-local ; {\'E}coulement {\`a} densit{\'e} variable ; Mod{\`e}le multicouches ; Simulation num{\'e}rique},
      TYPE = {Theses},
      PDF = {https://theses.hal.science/tel-02935869v3/file/these_boittin_lea_2019.pdf},
      HAL_ID = {tel-02935869},
      HAL_VERSION = {v3},
}

@book{rijn1987mathematical,
	title={Mathematical modelling of morphological processes in the case of suspended sediment transport},
	author={Rijn, Leonardus Cornelis},
    school       = {Waterloopkundig Laboratorium},
  year         = {1987},
  address      = {Delft, Netherlands},
  %url          = {https://repository.tudelft.nl/record/uuid:c1c1fce6-afc6-4ca3-b707-e108f256048c},
}

@article{el2011applicability,
  author    = {Kesserwani, Ghassan and Elghobashi, Salah and Labeeb, Ahmed},
  title     = {Applicability of Sediment Transport Capacity Formulas to Dam-Break Flows over Movable Beds},
  journal   = {Journal of Hydraulic Engineering},
  volume    = {136},
  number    = {5},
  pages     = {298--307},
  year      = {2010},
  doi       = {10.1061/(ASCE)HY.1943-7900.0000298},
  %url       = {https://ascelibrary.org/doi/abs/10.1061/(ASCE)HY.1943-7900.0000298}
}

@article{zhao2017comparison,
  author    = {Zhao, Jiaheng and Özgen-Xian, Ilhan and Liang, Dongfang and Hinkelmann, Reinhard},
  title     = {Comparison of depth-averaged concentration and bed load flux sediment transport models of dam-break flow},
  journal   = {Water Science and Engineering},
  volume    = {10},
  number    = {4},
  pages     = {287--294},
  year      = {2017},
  doi       = {10.1016/j.wse.2017.12.006},
  %url       = {https://www.sciencedirect.com/science/article/pii/S1674237017301059}
}

@article{audusse2021asymptotic,
  TITLE = {{Asymptotic derivation and simulations of a non-local Exner model in large viscosity regime}},
  AUTHOR = {Audusse, Emmanuel and Boittin, L{\'e}a and Parisot, Martin},
  % URL = {https://hal.science/hal-03312836},
  JOURNAL = {{ESAIM: Mathematical Modelling and Numerical Analysis}},
  PUBLISHER = {{Soci{\'e}t{\'e} de Math{\'e}matiques Appliqu{\'e}es et Industrielles (SMAI) / EDP}},
  VOLUME = {55},
  NUMBER = {4},
  PAGES = {1635-1668},
  YEAR = {2021},
  MONTH = Jul,
  DOI = {10.1051/m2an/2021031},
  KEYWORDS = {Free surface flow ; Shallow water equations ; Sediment transport ; Entropy dissipation ; Non-local effects},
 % url = {https://hal.science/hal-03312836v1/file/m2an200222.pdf},
  HAL_ID = {hal-03312836},
  HAL_VERSION = {v1},
}

@article{bonaventura2018multilayer,
    title = {Multilayer shallow water models with locally variable number of layers and semi-implicit time discretization},
    journal = {Journal of Computational Physics},
    volume = {364},
    pages = {209-234},
    year = {2018},
    issn = {0021-9991},
    doi = {10.1016/j.jcp.2018.03.017},
    %url = {https://www.sciencedirect.com/science/article/pii/S0021999118301694},
    author = {Bonaventura, Luca  and Fernández-Nieto, Enrique Domingo and Garres-Díaz, José and Narbona-Reina, Gladys},
    keywords = {Semi-implicit method, Multilayer approach, Depth-averaged model, Mass exchange, Sediment transport},
    abstract = {We propose an extension of the discretization approaches for multilayer shallow water models, aimed at making them more flexible and efficient for realistic applications to coastal flows. A novel discretization approach is proposed, in which the number of vertical layers and their distribution are allowed to change in different regions of the computational domain. Furthermore, semi-implicit schemes are employed for the time discretization, leading to a significant efficiency improvement for subcritical regimes. We show that, in the typical regimes in which the application of multilayer shallow water models is justified, the resulting discretization does not introduce any major spurious feature and allows again to reduce substantially the computational cost in areas with complex bathymetry. As an example of the potential of the proposed technique, an application to a sediment transport problem is presented, showing a remarkable improvement with respect to standard discretization approaches.}
}

@article{audusse2011multilayer,
  author    = {Audusse, Emmanuel and Bristeau, Marie-Odile and Perthame, Benoît and Sainte-Marie, Jacques},
  title     = {A multilayer Saint-Venant system with mass exchanges for shallow water flows. Derivation and numerical validation},
  journal   = {ESAIM: Mathematical Modelling and Numerical Analysis},
  volume    = {45},
  number    = {1},
  pages     = {169--200},
  year      = {2011},
  doi       = {10.1051/m2an/2010036},
  %url       = {https://www.numdam.org/articles/10.1051/m2an/2010036/}
}

@article{fernandez2013multilayer,
	title={A multilayer shallow water approach for polydisperse sedimentation with sediment compressibility and mixture viscosity},
	author={B{\"u}rger, Raimund and Fern{\'a}ndez-Nieto, Enrique D and Osores, V{\'\i}ctor},
	journal={Journal of Scientific Computing},
	volume={85},
	pages={1--40},
	year={2020},
	publisher={Springer},
	doi= {10.1007/s10915-020-01334-6}
}

@article{scholz2024dispersion,
  author    = {Scholz, Ullika and Kowalski, Julia and Torrilhon, Manuel},
  title     = {Dispersion in Shallow Moment Equations},
  journal   = {Computational and Applied Mathematics},
  volume    = {43},
  number    = {1},
  year      = {2023},
  pages     = {1-22},
  doi       = {10.1007/s42967-023-00325-2},
  %url       = {https://link.springer.com/article/10.1007/s42967-023-00325-2},
}

@article{bradford1999hydrodynamics,
  author    = {Bradford, Scott F. and Katopodes, Nikolaos Demetrios},
  title     = {Hydrodynamics of Turbid Underflows. I: Formulation and Numerical Analysis},
  journal   = {Journal of Hydraulic Engineering},
  volume    = {125},
  number    = {10},
  pages     = {1006--1015},
  year      = {1999},
  doi       = {10.1061/(ASCE)0733-9429(1999)125:10(1006)},
  %url       = {https://ascelibrary.org/doi/10.1061/(ASCE)0733-9429(1999)125:10(1006)},
}

@article{afzalimehr2009determination,
	title={Determination of bed shear stress in gravel-bed rivers using boundary-layer parameters},
	author={Afzalimehr, Hossein and Rennie, Colin},
	journal={Hydrological sciences journal},
	volume={54},
	number={1},
	pages={147--159},
	year={2009},
	publisher={Taylor \& Francis},
	doi = {10.1623/hysj.54.1.147},
}

@article{lajeunesse2010bed,
	title     = {Bed load transport in turbulent flow at the grain scale: Experiments and modeling},
	author    = {Lajeunesse, Eric and Malverti, Line and Charru, Fran{\c{c}}ois},
	journal   = {Journal of Geophysical Research: Earth Surface},
	volume    = {115},
	number    = {F4},
	pages     = {F04001},
	year      = {2010},
	publisher = {Wiley Online Library},
	doi       = {10.1029/2009JF001628},
	%url       = {https://agupubs.onlinelibrary.wiley.com/doi/full/10.1029/2009JF001628}
}

@article{elgamal2021moment,
	author    = {Elgamal, Mohamed},
	title     = {A Moment-Based Chezy Formula for Bed Shear Stress in Varied Flow},
	journal   = {Water},
	volume    = {13},
	number    = {9},
	pages     = {1254},
	year      = {2021},
	publisher = {MDPI},
	doi       = {10.3390/w13091254},
	%url       = {https://www.mdpi.com/2073-4441/13/9/1254}
}

@article{zordan2018entrainment,
	title     = {Entrainment, transport and deposition of sediment by saline gravity currents},
	author    = {Zordan, Jessica and Juez, Carmelo and Schleiss, Anton J. and Franca, M{\'a}rio J.},
	journal = {Advances in Water Resources},
	volume = {115},
	pages = {17-32},
	year = {2018},
	issn = {0309-1708},
	doi = {10.1016/j.advwatres.2018.02.017},
}

@article{rowan2020efficient,
	title={Efficient computational models for shallow water flows over multilayer erodible beds},
	author={Rowan, Thomas and Seaid, Mohammed},
	journal={Engineering Computations},
	volume={37},
	number={2},
	pages={401--429},
	year={2020},
	publisher={Emerald Publishing Limited}
}

@article{galappatti1985depth,
	title={A depth-integrated model for suspended sediment transport},
	author={Galappatti, G and Vreugdenhil, CB},
	journal={Journal of Hydraulic Research},
	volume={23},
	number={4},
	pages={359--377},
	year={1985},
	publisher={Taylor \& Francis}
}

@article{armanini1988one,
	title={A one-dimensional model for the transport of a sediment mixture in non-equilibrium conditions},
	author={Armanini, Aronne and Di Silvio, Giampaolo},
	journal={Journal of Hydraulic Research},
	volume={26},
	number={3},
	pages={275--292},
	year={1988},
	publisher={Taylor \& Francis}
}

@article{Bauerle2024,
	author    = {Bauerle, Matthew and Christlieb, Andrew J. and Ding, Mingchang and Huang, Juntao},
	title     = {On the Rotational Invariance and Hyperbolicity of Shallow Water Moment Equations in Two Dimensions},
	journal   = {SIAM Journal on Mathematical Analysis},
	volume    = {57},
	number    = {1},
	pages     = {421--445},
	year      = {2025},
	doi       = {10.1137/23M1579789},
	url       = {https://epubs.siam.org/doi/abs/10.1137/23M1579789},
	publisher = {Society for Industrial and Applied Mathematics}
}

@article{koellermeier2026primitive,
	author    = {Koellermeier, Julian},
	title     = {Primitive Variable Regularization to Derive Novel Hyperbolic Shallow Water Moment Equations},
	journal   = {Communications in Computational Physics},
	volume    = {36},
	number    = {6},
	pages     = {1325--1348},
	year      = {2026},
	doi       = {10.4208/cicp.OA-2025-0155},
	url       = {https://global-sci.com},
	publisher = {Global Science Press}
}

@article{vanLeer1979MUSCL,
	author  = {van Leer, Bram},
	title   = {Towards the Ultimate Conservative Difference Scheme. {V}. A Second-Order Sequel to {Godunov}'s Method},
	journal = {Journal of Computational Physics},
	volume  = {32},
	number  = {1},
	pages   = {101--136},
	year    = {1979},
	doi     = {10.1016/0021-9991(79)90145-1}
}

@article{rusanov1962,
	author  = {Rusanov, V. V.},
	title   = {The Calculation of the Interaction of Non-Stationary Shock Waves and Obstacles},
	journal = {USSR Computational Mathematics and Mathematical Physics},
	volume  = {1},
	number  = {2},
	pages   = {304--320},
	year    = {1962},
	doi     = {10.1016/0041-5553(62)90062-9},
	note    = {English translation of Zh. Vychisl. Mat. Mat. Fiz. 1(2), 267--279 (1961)}
}

@article{kurganovLevy2002,
	author  = {Kurganov, Alexander and Levy, Doron},
	title   = {Central-Upwind Schemes for the {Saint-Venant} System},
	journal = {ESAIM: Mathematical Modelling and Numerical Analysis},
	volume  = {36},
	number  = {3},
	pages   = {397--425},
	year    = {2002},
	doi     = {10.1051/m2an:2002019}
}

@article{pares2006,
	author  = {Par{\'e}s, Carlos},
	title   = {Numerical Methods for Nonconservative Hyperbolic Systems: A Theoretical Framework},
	journal = {SIAM Journal on Numerical Analysis},
	volume  = {44},
	number  = {1},
	pages   = {300--321},
	year    = {2006},
	doi     = {10.1137/050628052}
}

@article{castroGallardoPares2006,
	author  = {Castro, Manuel J. and Gallardo, Jos{\'e} M. and Par{\'e}s, Carlos},
	title   = {High Order Finite Volume Schemes Based on Reconstruction of States for Solving Hyperbolic Systems with Nonconservative Products. Applications to Shallow-Water Systems},
	journal = {Mathematics of Computation},
	volume  = {75},
	number  = {255},
	pages   = {1103--1134},
	year    = {2006},
	doi     = {10.1090/S0025-5718-06-01851-5}
}

@misc{scholz2026splitting,
  author        = {Scholz, Ullika and Paar, Robin and Torrilhon, Manuel},
  title         = {A Splitting Scheme for Dispersive Shallow Moment Equations},
  year          = {2026},
  eprint        = {2607.08870},
  archivePrefix = {arXiv},
  primaryClass  = {math.NA},
  doi           = {10.48550/arXiv.2607.08870}
}

@incollection{KreissOrtiz2002,
	author    = {Kreiss, Heinz-Otto and Ortiz, Omar E.},
	title     = {Some Mathematical and Numerical Questions Connected with First and Second Order Time-Dependent Systems of Partial Differential Equations},
	booktitle = {The Conformal Structure of Space-Time},
	editor    = {Frauendiener, J{\"o}rg and Friedrich, Helmut},
	series    = {Lecture Notes in Physics},
	volume    = {604},
	pages     = {359--370},
	publisher = {Springer},
	address   = {Berlin, Heidelberg},
	year      = {2002},
	doi       = {10.1007/3-540-45818-2_21}
}

@article{SmithPettyPantano2016,
	author  = {Smith, Tyler A. and Petty, David J. and Pantano, Carlos},
	title   = {A Roe-like Numerical Method for Weakly Hyperbolic Systems of Equations in Conservation and Non-Conservation Form},
	journal = {Journal of Computational Physics},
	volume  = {316},
	pages   = {117--138},
	year    = {2016},
	doi     = {10.1016/j.jcp.2016.04.002}
}
